\documentclass{article}

\usepackage[english]{babel}
\usepackage[utf8]{inputenc}

\usepackage{imakeidx}

\usepackage{emptypage}

\usepackage[nottoc]{tocbibind}

\usepackage{enumitem}

\usepackage{geometry}

\usepackage{empheq}
\usepackage{a4wide}
\usepackage{caption}
\usepackage{amsmath}
\usepackage{xfrac}
\usepackage{nicefrac}
\usepackage{faktor}
\usepackage{hyperref}
\usepackage{comment}
\usepackage{enumerate}

\usepackage{musicography}
\usepackage{amssymb, amsthm, dsfont}
\usepackage{mathtools}
\usepackage{units}
\usepackage{cleveref}  
\usepackage{csquotes}
\usepackage{xcolor}
\usepackage{graphicx}
\usepackage{tikz}
\usepackage{extarrows}

\usepackage{tikz-cd}
\tikzcdset{scale cd/.style={every label/.append style={scale=#1},  cells={nodes={scale=#1}}}}

\usepackage{float}

\newlist{myitems}{enumerate}{2}

\setlist[myitems, 1]
{label=\Roman{myitemsi}., 
leftmargin=50pt, 
rightmargin=10pt
}

\setlist[myitems, 2]
{label=(\roman{myitemsii}), 
leftmargin=15pt,
rightmargin=15pt}

\newcommand{\NN}{\mathbb{N}} 								                        
\newcommand{\ZZ}{\mathbb{Z}}  								                        
\newcommand{\RR}{\mathbb{R}}  								                        
\newcommand{\CC}{\mathbb{C}}  								                        
\newcommand{\KK}{\mathbb{K}}  								                        
\newcommand{\DD}{\mathbb{D}}  								                        
\renewcommand{\SS}{\mathbb{S}}  							                        
\renewcommand{\a}{\alpha}                                                           
\newcommand{\GL}{\operatorname{GL}}                                                 
\newcommand{\SL}{\operatorname{SL}}                                                 
\newcommand{\topo}[1]{\tau_{\text{#1}}}                                             
\newcommand{\p}{\varphi}                                                            
\newcommand{\Aut}{\operatorname{Aut}}                                               
\newcommand{\erzeug}[2]{\operatorname{span}_{#1}\hspace{-2pt}\left( #2 \right)}     
\renewcommand{\a}{\alpha}                                                           
\newcommand{\C}{\mathcal{C}}

\renewcommand{\epsilon}{\varepsilon}
\newcommand{\Hom}{\operatorname{Hom}}
\newcommand{\Sym}{\operatorname{Sym}}

\newcommand{\N}[2]{N_{#1}\hspace{-1.5pt}\left( #2 \right)}                          

\newcommand{\EW}{\widetilde{\W}}							

\newcommand{\KMExp}{\exp_{\operatorname{KM}}}

\newcommand{\Kac}{\mathrm{Kac}}
\renewcommand{\DD}{\Delta^{\mathrm{re}}}
\newcommand{\Stab}{\operatorname{Stab}}                     
\newcommand{\Fix}{\operatorname{Fix}}                       
\newcommand{\sing}[1]{{#1^\text{sing}}}						
\newcommand{\reg}[1]{{#1^\text{reg}}}						
\renewcommand{\aa}{\mathfrak{a}}                            

\newcommand{\F}{\mathcal{F}}			                    

\renewcommand{\AA}{\mathbb{A}}                              
\newcommand{\BB}{\textbf{B}}                                
\newcommand{\Ad}{\operatorname{Ad}} 
\newcommand{\ad}{\operatorname{ad}} 
\newcommand{\G}{\mathcal{G}} 				                
\newcommand{\g}{\mathfrak{g}}								
\newcommand{\h}{\mathfrak{h}}								
\newcommand{\U}{\mathcal{U}} 								
\newcommand{\D}{\mathcal{D}}                                
\newcommand{\W}{\mathcal{W}} 								
\newcommand{\T}{\mathfrak{T}}                              
\newcommand{\Grp}{\textbf{Grp}}                       
\newcommand{\AutEff}{\Aut_{\text{eff}}}						
\newcommand{\Trans}{\operatorname{Trans}}				   
\newcommand{\eff}[1]{\left( #1 \right)_{\operatorname{eff}}}
\newcommand{\os}{\overline{\sigma}}                                     
\newcommand{\alsp}[1]{#1^{\operatorname{as}}}                           
\newcommand{\res}[2]{#2 \hspace{-3pt} \mid_{#1}}                                      
\renewcommand{\root}[1]{#1^{\musNatural}}                                
\newcommand{\rootgrp}[1]{V_{#1}^{\musNatural}}                           
\newcommand{\mbs}[1]{\operatorname{MBS}(#1)}
\newcommand{\rel}{\musNatural}

\newcommand{\SGR}{\textbf{(SGR)}}
\newcommand{\DCS}{\textbf{DCS}}

\newtheoremstyle{MyDef}
{15pt}
{10pt}
{}
{}
{\bfseries}
{}
{\newline}
{}

\newtheoremstyle{MyPef}
{15pt}
{10pt}
{}
{}
{\itshape}
{:}
{ }
{}

\newtheoremstyle{MySef}
{15pt}
{10pt}
{\itshape}
{}
{\bfseries}
{}
{\newline}
{}

\theoremstyle{MyDef}
\newtheorem{Definition}{Definition}[section]
\newtheorem{Convention}[Definition]{Convention}
\newtheorem{Example}[Definition]{Example}
\newtheorem{Remark}[Definition]{Remark}

\theoremstyle{MyPef}
\newtheorem*{Proof}{Proof}

\theoremstyle{MySef}
\newtheorem{Lemma}[Definition]{Lemma}

\newtheorem{Corollary}[Definition]{Corollary}
\newtheorem{Proposition}[Definition]{Proposition}
\newtheorem{Theorem}[Definition]{Theorem}

\usepackage[framemethod=TikZ]{mdframed}
\mdfdefinestyle{mpdframe1}{
    frametitlebackgroundcolor   =black!15,
    frametitlerule          =true,
        roundcorner     =10pt,
        backgroundcolor = yellow!65,
        middlelinewidth     =1pt,
        innermargin     =0.5cm,
        outermargin     =0.5cm,
        innerleftmargin     =0.5cm,
        innerrightmargin    =0.5cm,
        innertopmargin      =\topskip,
        innerbottommargin   =\topskip,
            }
            
\mdfdefinestyle{mpdframe2}{
    frametitlebackgroundcolor   =black!15,
    frametitlerule          =true,
        roundcorner     =10pt,
        backgroundcolor = blue!40,
        middlelinewidth     =1pt,
        innermargin     =0.5cm,
        outermargin     =0.5cm,
        innerleftmargin     =0.5cm,automorphisms.
        innerrightmargin    =0.5cm,
        innertopmargin      =\topskip,
        innerbottommargin   =\topskip,
            }
            
\mdfdefinestyle{mpdframe3}{
    frametitlebackgroundcolor   =green!15,
    frametitlerule          =true,
        roundcorner     =2.5pt,
        middlelinewidth     =1pt,
        innermargin     =0.5cm,
        outermargin     =0.5cm,
        innerleftmargin     =0.5cm,
        innerrightmargin    =0.5cm,
        innertopmargin      =\topskip,
        innerbottommargin   =\topskip,
            }
            
\mdfdefinestyle{mpdframe4}{
    frametitlebackgroundcolor   =red!40,
    frametitlerule          =true,
        roundcorner     =2.5pt,
        middlelinewidth     =1pt,
        innermargin     =0.5cm,
        outermargin     =0.5cm,
        innerleftmargin     =0.5cm,
        innerrightmargin    =0.5cm,
        innertopmargin      =\topskip,
        innerbottommargin   =\topskip,
            }
            
\mdfdefinestyle{question}{%
        style=mpdframe1,
        frametitle={Question},
    }
\newmdenv[style=question]{que}

\mdfdefinestyle{answer}{%
        style=mpdframe2,
        frametitle={Answer},
    }
\newmdenv[style=answer]{ans}

\mdfdefinestyle{solution}{%
        style=mpdframe3,
        frametitle={Solution},
    }
\newmdenv[style=solution]{sol}

\mdfdefinestyle{hint}{%
        style=mpdframe4,
        frametitle={Watch out!},
    }
\newmdenv[style=hint]{hint}

\title{Kac--Moody Symmetric Spaces: arbitrary symmetrizable complex or almost split real type}
\author{Ralf Köhl and Christian Vock}
\date{\today}

\newcommand{\Addresses}{{
  \bigskip
  \footnotesize

  Ralf Köhl, \textsc{Mathematisches Seminar, Christian-Albrechts-Universität zu Kiel,
    Heinrich-Hecht-Platz 6, 24118 Kiel}\par\nopagebreak
  E-mail address: \texttt{koehl@math.uni-kiel.de}

  \medskip

  Christian Vock, \textsc{Mathematisches Seminar, Christian-Albrechts-Universität zu Kiel,
    Heinrich-Hecht-Platz 6, 24118 Kiel}\par\nopagebreak
  E-mail address: \texttt{vock@math.uni-kiel.de}

}}

\begin{document}

\maketitle

\begin{abstract}
We revisit the theory of Kac--Moody symmetric spaces (of "non-compact type") in order to include the affine case, complex numbers, and Galois descent to almost split real Kac--Moody symmetric spaces.
\end{abstract}

\section{Introduction}

Non-compact (semi)simple Lie groups admit two natural geometries: a symmetric space of non-compact type and a Tits building in the boundary at infinity of the symmetric space; cf.\ \cite{BurnsSpatzier}. From a purely group-theoretic point of view, the symmetric space is suitable for establishing, resp.\ confirming the simplicity of the acting Lie group as an abstract group -- the simplicity is related to the maximality of the maximal compact subgroup (that is, the maximal compact subgroup actually being a maximal subgroup) with in turn is related to the primitivity of the Lie group action on the symmetric space; cf.\ \cite{Caprace_2005}.  
Kac--Moody groups have been extensively studied since the 1980s, initiated by \cite{Kac-Peterson_1985} and \cite{Tits_1987}, and it turns out that they actually admit two buildings which form a so-called twin building. However, Kac--Moody symmetric spaces have only been formally defined and investigated in 2020 in \cite{FHHK_2017} and only over the real numbers and only for indefinite type, although the heuristic idea that Kac--Moody symmetric spaces might be useful has been around in particular in mathematical physics for some time; cf.\ \cite[Section~8.3]{DHN}. Kac--Moody symmetric spaces actually admit a causal structure defining a causal partial order (\cite[Proposition~7.46]{FHHK_2017} combined with \cite[Theorem~6.3.18]{Diss2}). The twin building can be recovered from inside the causal boundary of the Kac--Moody symmetric space; cf.\ \cite[Corollary~7.41]{FHHK_2017}. 

\smallskip
It is still an open question whether split real Kac--Moody groups of indefinite type are abstractly simple. One of the goals of \cite{FHHK_2017} was to propose a symmetric space that might allow a primitive action as in \cite{Caprace_2005} -- in order to establish the abstract simplicity of the acting group. In case of non-invertible generalized Cartan matrix this lead to some modeling choices that eventually excluded affine types -- which is not surprising as affine Kac--Moody groups are well known to not be abstractly simple.

\smallskip
In the present article we make different modeling choices. For split real Kac--Moody groups with invertible generalized Cartan matrix they lead to exactly the same class of symmetric spaces -- see also the discussion in Section~\ref{comparing} below. For a non-invertible generalized Cartan matrix these modeling choices lead to a symmetric space which obviously is imprimitive, and they allow to also consider the case of an affine generalized Cartan matrix. 
Moreover, in the present article we also include the field $\mathbb{C}$ of complex numbers (see Convention~\ref{convention3.5}), and we conduct Galois descent based on \cite{Remy_2002}, in order to model Kac--Moody symmetric spaces also in the almost split real case; see Section~\ref{almostsplit}.

\smallskip
The following result is our main existence theorem. 

\smallskip
{\bf Proposition \ref{Prop:KM-SymSp}} and {\bf Theorem \ref{Thm:Galois-Sym-Sp}}\\
{\em
Let $\mathbb{K} \in \{ \mathbb{R}, \mathbb{C} \} $ and let $G$ be a split or almost split Kac--Moody $\mathbb{K}$-group (satisfying \eqref{Eq:DCS} in the almost split situation) and denote by $G(\mathbb{K})$ the $\mathbb{K}$-rational points of $G$. Let $\Theta$ be the Cartan--Chevalley involution of $G(\mathbb{K})$ over the real numbers, resp.\ the Cartan--Chevalley involution composed with the continuous field involution over the complex numbers, and let $K$ be the subgroup of $G$ of $\Theta$-fixed elements. Then the quotient 
$$
X(\mathbb{K}) \coloneqq \nicefrac{G(\mathbb{K})}{K(\mathbb{K})}
$$
is a symmetric space in the sense of Loos \cite{Loos_1969} with respect to the reflection map $$(gK(\mathbb{K}),hK(\mathbb{K})) \mapsto g\Theta(g)^{-1}\Theta(h)K(\mathbb{K}).$$ Moreover, there is a natural action of $G(\mathbb{K})$ on the symmetric space given by automorphisms
\begin{align*}
    G(\mathbb{K}) &\to \Sym(X(\mathbb{K})) \\
    g &\mapsto (hK(\mathbb{K}) \mapsto ghK(\mathbb{K})).
\end{align*}
}

\smallskip
The condition \eqref{Eq:DCS} is only relevant for the almost split situation and is central to the theory of almost split Kac--Moody groups: it guarantees that the rational points of the Kac--Moody group correspond to the fixed points of the complex Kac--Moody group under Galois involution. For details on this involution see \cref{Sec:RR-Forms} and for details on the condition see \cref{Sec:GaloisDescent}.

\smallskip 
In analogy to \cite{FHHK_2017} we also study flats and automorphisms of these symmetric spaces and the Weyl group action.
We took quite some effort to make the present article accessible independently from \cite{FHHK_2017}, but it is quite natural that we need to borrow and cite a wealth of results from \cite{FHHK_2017}; on top of that, also a detailed comparison of our results with the results from \cite{FHHK_2017} is called for.

\tableofcontents

\medskip

\textbf{Acknowledgment:}\\
The authors thank the DFG SPP 2026 priority program “Geometry at infinity” for partial financial support via KO 4323/14. During the first half of the project both authors were affiliated with and funded by JLU Gie\ss en, during the second half affiliated with and funded by CAU Kiel. The authors would also like to thank Paul Zellhofer in particular for a wealth of very helpful discussions and comments. Moreover, the authors thank an anonymous referee for literally hundreds of valuable observations and comments that tremendously helped improve the article. This article is based on the PhD thesis of the second author that has been advised by the first author.


\section{Basics on Symmetric Spaces}\label{Sec:SymSp}

In \cite{Loos_1969}, Loos analyzes Riemannian and pseudo-Riemannian symmetric spaces and states an abstract definition without any use of differential geometry. This makes the definition amenable to situations in which one tries to generalize symmetric spaces into a context without a differential structure.

\begin{Definition}[Definition 1, Chapter II in \cite{Loos_1969}]\label{Def:SymSp-Loos}
A topological space $X$ with a continuous map, which is called \emph{reflection}, $\mu \colon X \times X \to X$, $(x,y) \mapsto \mu(x,y) \coloneqq x.y$ is called a \emph{symmetric space} $(X, \mu)$ if the map satisfies the following properties:
\begin{myitems}[label=\textbf{S.\arabic*}]
    \item \label{itm:S1} $x.x = x$ for all $x \in X$,
    \item \label{itm:S2} $x.(x.y) = y$ for all $x,y \in X$,
    \item \label{itm:S3} $x.(y.z) = (x.y).(x.z)$ for all $x,y,z \in X$, and
    \item \label{Local} for every $x \in X$ there is a neighborhood $U$ of $x$ such that $x.y = y$ implies $x = y$ for all $y \in U$.
\end{myitems}
\end{Definition}

If the Riemannian symmetric space is of non-compact type, i.e.\ $X$ is non-compact, one can use the global condition instead of the local condition \ref{Local}:
\begin{align*}
\begin{matrix}
\textbf{S.4}_{\text{global}} & \forall x,y \in X \; \; x.y = y \Longrightarrow x=y.
\end{matrix}
\end{align*}

Since we are interested in generalizations of Riemannian symmetric spaces of non-compact type, the class of symmetric spaces that we consider in this article will automatically satisfy $\textbf{S.4}_{\text{global}}$ instead of just \ref{Local}. Note that this in fact allows one to model symmetric spaces on sets instead of topological spaces.   \\
A \emph{morphism of symmetric spaces} is a (continuous) map $f \colon X_1 \to X_2$ between two symmetric spaces $(X_1, \mu_1)$ and $(X_2, \mu_2)$ satisfying $f(\mu_1(x,y)) = \mu_2(f(x), f(y))$. A pair $((X,\mu),\beta)$ consisting of a symmetric space $(X,\mu)$ and a point $\beta \in X$ is called a \emph{pointed symmetric space}. A morphism between pointed symmetric spaces is a morphism between their underlying symmetric spaces with the additional property that the base point is preserved.\\
In the rest of the article, if the map $\mu$ is clear from the context, we will write $X$ for the symmetric space pair $(X,\mu)$ for short.

\begin{Example}\label{Ex:SymSpaces}
\begin{myitems}
    \item[]
    
    \item \label{itm:EuclSymSpace} One of the first examples that comes to mind of a symmetric space is the $n$-dimensional Euclidean space $\mathbb{E}^n$ for any $n \in \NN$. In particular, the symmetric space is given by the pair $\RR^n$ together with the point reflection on $x \in \RR^n$, i.e. the symmetric space map $\mu_{\text{Eucl}}$ is defined as $\mu_{\text{Eucl}}(x,y) \coloneqq 2x - y$.

    \item Considering the sphere with the corresponding geodesic reflection leads to an example, which satisfies only \ref{Local}. 

    \item An example of a Riemannian symmetric space of non-compact type is the hyperbolic space. Since there are no closed geodesics, the property $\textbf{S.4}_{\text{global}}$ is satisfied.

    \item \label{itm:GroupModel} Another example is any topological group $G$ with the reflection $\mu \colon G \times G \to G$, $\mu(g,h) \coloneqq hgh^{-1}$. Here \ref{Local} is in general not satisfied, but one can easily check that \ref{itm:S1} - \ref{itm:S3} are satisfied. For a detailed calculation of the first three axioms for topological groups, we refer to the proof of \cite[Porposition 4.2]{FHHK_2017}. There, an assumption is made that guarantees the fourth axiom, and this assumption also holds for the Kac--Moody situation, as we will see in \Cref{Lem:tau(G)-K-e}.
\end{myitems}
\end{Example}

Define the \emph{point reflection} at $x \in X$ as the map
$$
s_x \colon X \to X \, , \; y \mapsto \mu(x,y).
$$
The product of two point reflections is called \emph{transvection} and the group
$$
\Trans(X) \coloneqq \langle s_x \circ s_y \mid x,y \in X \rangle
$$
is called the \emph{transvection group} of $X$. 

\begin{Lemma}[p.64, line 15 of \cite{Loos_1969}]\label{Lem:Eigenschaft-Spiegelungen}
Let $(X, \mu)$ be a symmetric space, $x$ a point in $X$ and let $\alpha \in \Aut(X)$ be an automorphism of the symmetric space. Then
$$
\alpha \circ s_x \circ \alpha^{-1} = s_{\alpha(x)}.
$$
\end{Lemma}
\begin{Proof}
For each $p \in X$ one calculates
$$
\left(\alpha \circ s_x \circ \alpha^{-1} \right) \hspace{-2pt}(p) = \alpha\hspace{-2pt}\left( s_x\hspace{-2pt}\left( \alpha^{-1}(p) \right) \right) = \alpha\hspace{-2pt}\left( \mu\hspace{-2pt}\left(x,\alpha^{-1}(p) \right) \right) = \mu(\alpha(x), p) = s_{\alpha(x)}\hspace{-1pt}(p).
$$
\qed
\end{Proof}

By setting $\alpha = s_y$ for $y \in X$ in the previous lemma, we obtain the useful identity 
\begin{align}\label{Eq:Verkettung-SxSy}
s_y \circ s_x \circ s_y^{-1} = s_y \circ s_x \circ s_y = s_{s_y(x)}.    
\end{align}

Consider a pointed symmetric space $(X, \mu)$ with base point $b \in X$. We now want to describe the symmetric space via the \emph{quadratic representation}, see \cite[Section II.1]{Loos_1969} or \cite[Remark 2.10]{FHHK_2017}.\\
For $x \in X$ define 
$$
t_x \coloneqq s_x \circ s_b \in \Trans(X)
$$
and
$$
T(X,b) \coloneqq \{ t_x \mid x \in X \}.
$$
By defining the following map
$$
\mu_T \colon T(X,b) \times T(X,b) \to T(X,b) \, , \; (s,t) \mapsto s \circ t^{-1} \circ s,
$$
the set $T(X,b)$ is turned into a symmetric space with reflection map $\mu_T$.

\begin{Proposition}\label{Rem:QaudraticRepresentation}
The symmetric spaces $\left(T(X,b), \mu_T \right)$ and $(X,\mu)$ are isomorphic. 
\end{Proposition}
\begin{Proof}
Evidently, the map
$$
\eta \colon X \to T(X,b) \, , \; x \mapsto t_x
$$
is surjective. Since $t_x \circ s_b = s_x$ and since $s_x = s_y$ implies $x = y$ (due to axiom $\textbf{S.4}_{\text{global}}$), the map $\eta$ is also injective. A direct calculation shows that $\eta \colon (X,\mu) \to (T(X,b),\mu_T)$ is an isomorphism of symmetric spaces. Indeed, one has to check that the reflection map of $(X, \mu)$ is preserved, i.e.\ $\eta(\mu(x,y)) = \mu_T(t_x, t_y)$. On the right-hand side one can calculate
\begin{align*}
    \mu_T(t_x,t_y) &= t_x \circ t_y^{-1} \circ t_x \\
    &= (s_x \circ s_b) \circ (s_y \circ s_b)^{-1} \circ (s_x \circ s_b) \\
    &= s_x \circ s_b \circ s_b \circ s_y \circ s_x \circ s_b \\
    &= s_x \circ s_y \circ s_x \circ s_b \\
    &\overset{\ref{Eq:Verkettung-SxSy}}{=} s_{s_x(y)} \circ s_b = t_{s_x(y)}.
\end{align*}
And on the left-hand side one has 
$$
\eta(\mu(x,y)) = \eta(s_x(y)) = t_{s_x(y)}.
$$
\qed
\end{Proof}

Recall the following concepts of an abstract symmetric space $(X,\mu)$ from \cite[Section 2.15]{FHHK_2017}:
\begin{itemize}
    \item A \emph{midpoint} of two points $x,y \in X$ is a point $m \in X$, such that $\mu(m,x) = y$ and $\mu(m,y) = x$. 
    \item Call a subset $U \subseteq X$
        \begin{itemize}
            \item a \emph{reflection subspace}, if for all pair of points $x,y \in U$ the image $\mu(x,y)$ is contained in $U$.
            \item \emph{midpoint convex}, if for all $x,y \in U$ there is a midpoint of $x$ and $y$ in $U$. 
        \end{itemize}
        \item Two points $x,y \in X$ \emph{weakly commute}, if for every point $p \in X$ we have: $x.(p.(y.p)) = y.(p.(x.p))$.
        \item Two points $x,y \in X$ \emph{commute}, if for all points $p,q \in X$ we have: $ x.(p.(y.q)) = y.(p.(x.q)) $.
        \item A reflection subspace $F \subseteq X$ is a \emph{(weak) flat}, if it is closed, midpoint convex, contains at least two points and if all points $x,y \in F$ (weakly) commute.
\end{itemize}

\begin{Definition}[cf.\ Definition 2.23 from \cite{FHHK_2017}]
A  \emph{Euclidean flat of rank $n$}  is a closed symmetric subspace $F \subseteq X$ isomorphic to the $n$-dimensional Euclidean space $\mathbb{E}^n$ as a symmetric space. For $n=1$, the flat is called a \emph{geodesic}.
\end{Definition}


\section{Kac--Moody Groups}\label{Sec:Basics}
In this section, we give a brief overview of the basic notions as well as the construction of a Kac--Moody group. Most of this can be found in \cite{Kac_1990}, \cite{Tits_1987}, \cite{Remy_2002} or \cite{Marquis_2018}.

\subsection{Kac--Moody Algebra}\label{Sec:KM_Algebra}
\begin{Definition}
Let $\AA = \left( a_{ij} \right)_{1 \leq i,j \leq n} \in \ZZ^{n \times n}$ be a square matrix over the integers. The matrix $\AA$ is called a \emph{generalized Cartan matrix} if it satisfies the following properties
\begin{itemize}
    \item $a_{ii} = 2$ for all $i = 1, \ldots, n $,
    \item $a_{ij} \leq 0$ for all $i,j = 1, \ldots, n $, $i \neq j$, and
    \item $a_{ij} = 0$ if and only if $a_{ji} = 0$.
\end{itemize}
\end{Definition}

The generalized Cartan matrix $\AA$ is \emph{indecomposable} if it cannot be decomposed into a non-trivial direct sum after rearranging the indices, cf.\ \cite[§1.1]{Kac_1990}. Furthermore, $\AA$ is \emph{symmetrizable} if there is a symmetric matrix $B \in \RR^{n \times n}$ and an invertible diagonal matrix $D \in ( \RR^\times )^{n \times n}$ such that $\AA = D B$, cf.\ \cite[§2.1]{Kac_1990}.

If a generalized Cartan matrix is the Cartan matrix of a finite-dimensional semisimple complex Lie algebra, then $\AA$ is of \emph{spherical type}; otherwise, $\AA$ is \emph{non-spherical}.

In \cite{Kac_1990} the notion of a realization of $\AA$ is introduced, which assigns a triple to a generalized Cartan matrix, allowing the construction of a Kac--Moody algebra. As in the finite-dimensional situation, a Kac--Moody algebra is naturally related to several different groups; therefore, in order to define a Kac--Moody group associated to the algebra one needs a more subtle definition (see \cite[Chapitre 8]{Remy_2002} or \cite[Definition 7.9]{Marquis_2018}).

\begin{Definition}\label{Def:RootDatum}
Let $\AA = \left( a_{ij} \right)_{i,j \in I}$ be a generalized Cartan matrix with index set $I$, let $\Lambda$ be a free $\ZZ$-module of finite rank and $\Lambda^\vee$ the $\ZZ$-dual of $\Lambda$, and let $(c_i)_{i \in I} \subset \Lambda$ and $(h_i)_{i \in I} \subset \Lambda^\vee$ such that $c_j(h_i) = a_{ij}$. The quintuple $\D \coloneqq \left( I, \AA, \Lambda, (c_i)_{i \in I}, (h_i)_{i \in I} \right)$ is called a \emph{Kac--Moody root datum}.
\end{Definition}

Consider a generalized Cartan matrix $\AA$ of arbitrary type and a corresponding Kac--Moody root datum $\D$. Define the \emph{Cartan subalgebra} 
$$\h_\D \coloneqq \Lambda^\vee \otimes_\ZZ \CC$$
with elements, $\a_i^\vee \coloneqq h_i \otimes 1$ and define the dual space 
$$\h_\D^* \coloneqq \Lambda \otimes_\ZZ \CC$$
with elements $\a_i \coloneqq c_i \otimes 1$. Call the elements $\Pi = \{\a_i \mid i \in I \} \subset \h^*$ \emph{simple roots}, $\Pi^\vee = \{\a_i^\vee \mid i \in I\} \subset \h$ \emph{simple coroots} and denote with $\Delta^{re}$ the set of real roots (cf.\ \cite[§5.1]{Kac_1990} or \cite[Chapter 6.1]{Marquis_2018}).

\begin{Remark}\label{Rem:AffineGCM}
The original approach of Kac is to define that the sets $\Pi$ and $\Pi^\vee$ are linearly independent, see \cite[§1.1]{Kac_1990}. 
        In the general case, where the generalized Cartan matrix $\AA$ is not necessarily invertible (e.g.\ if $\AA$ is of affine type), the consequence is that $\Pi \subset \h^*$ and $\Pi^\vee \subset \h$ cannot be assumed to be linearly independent. 

To address this problem, one can extend the Cartan subalgebra $\h$, if necessary, to ensure that the $\a_i$ become linearly independent. For a concrete discussion of this, see \cite[Proof of Proposition 1.1]{Kac_1990} and \cite[Chapter 3.5 and Example 7.10]{Marquis_2018}. It turns out that for $\AA \in \ZZ^{n \times n}$, where $l$ denotes the rank of $\AA$, the extended Cartan subalgebra must have dimension $2n-l$ and is isomorphic to $\CC^{2n-l}$. The dual space is also extended and constructed in such a way that the simple roots still have the relation $\langle \a_i, \a_j^\vee \rangle = a_{ij} \in \AA$.\\
On this basis, we denote the following as $\D_{\Kac}^\AA$ the root datum where the elements $c_i$ and $h_i$ are linearly independent and denote by $\Lambda_\Kac$ the corresponding $\ZZ$-module, which has minimum rank with respect to these properties, i.e.\ the rank of $\Lambda_\Kac$ is equal to $2n - l$. Similarly, we use $\Lambda_\Kac^\vee$ to denote the $\ZZ$-dual. Moreover, denote by $\lbrace v_i \mid 1 \leq i \leq 2n-l \rbrace$ the $\ZZ$-basis of the extended $\Lambda_\Kac$, where $h_i = v_i$ for $1 \leq i \leq n$. We follow here the notation and idea of \cite[Example 7.10]{Marquis_2018}.

\end{Remark}

Following \cite[Chapter 7]{Remy_2002} or \cite[Definition 7. 13]{Marquis_2018}, one can define a \emph{Kac--Moody algebra $\g_\D$ of type $\D$} as a Lie algebra with generators $\h_\D$, $\{e_i\}_{i \in I}$ and $\{f_i\}_{i \in I}$, where $e_i$ and $f_i$ are to be symbols satisfying the following relations.
\begin{align*}
    [ \h_\D , \h_\D ] &= 0 \\
    [h, e_i] &= c_i(h) e_i \\
    [h, f_i] &= -c_i(h) f_i \\
    [e_i, f_j] &= - \delta^i_j h_i \\
    (\ad e_i)^{1- a_{ij}} e_j &= (\ad f_i)^{1-a_{ij}} f_j = 0 \; (i \neq j)
\end{align*}
for $h \in \h_\D$ and $i, j \in I$. Note that $\delta^i_j$ here means the Kronecker symbol. In the following, we denote with $\g(\AA)$ the Kac--Moody algebra of type $\D_{\Kac}^\AA$.

\begin{Convention} \label{convention3.4}
For the remainder of the article we consider the root datum $\D_{\Kac}^\AA$, where $\AA$ is an indecomposable, symmetrizable generalized Cartan matrix of size $n$ with rank $l$. 
\end{Convention}

\medskip

Let $\AA$ be a generalized Cartan matrix and consider $\Pi = \{\a_i \mid i \in I \}$, then define the \emph{root lattice} as 
\begin{equation}\label{Eq:RootLattice}
    Q \coloneqq \bigoplus_{i \in I} \ZZ \a_i
\end{equation}
and call an element $\a \neq 0 \in Q$ a \emph{root} if 
$$
\lbrace x \in \g(\AA) \mid \ad_h(x) = \a(h)x \; \forall h \in \h \rbrace \neq \{ 0 \}.
$$
In \cite[Chapter 3.4]{Marquis_2018} the concept of \emph{gradation} is explained, which describes the representation of a vector space as the direct sum of subspaces. Following this, this root lattice is called a \emph{$Q$-gradation} of the root space.\\
Denote by $\Delta$ all roots of $\g(\AA)$. Furthermore, a root $\a$ is called \emph{positive} if the coefficients of this representation are non-negative, and \emph{negative} otherwise. By \cite[Proposition 3.14 (5)]{Marquis_2018} follows that every root is either positive or negative.

\medskip

Next, recall the \emph{Weyl group} of $\g(\AA)$ from \cite[§3.7]{Kac_1990}, \cite[Chapter 4.2]{Marquis_2018} or \cite[7.1.3]{Remy_2002}: it is the group $\W < \GL(\h)$ generated by the set $S = \{ r_1^\vee , \ldots, r_n^\vee \}$ of reflections
$$
r_i^\vee \colon \h \to \h \, , \; h \mapsto h - \a_i(h)\a_i^\vee.
$$
A root $\a \in \Delta$ is called \emph{real}, if there exists a $w \in \W$, such that $w \a \in \Pi$. Denote by $\DD \coloneqq \W. \Pi$ the set of all real roots in $\g(\AA)$.

\subsection{The Tits Functor}\label{Sec:TitsFunctor}\label{Rem:ExpMap}

Let $R$ be a ring and recall that $\D_{\Kac}^\AA = \left( I, \AA, \Lambda_\Kac, (c_i)_{i \in I}, (h_i)_{i \in I} \right)$. Then a split torus scheme (cf.\ \cite[Section 7.4]{Marquis_2018}) is given by 
$$
T \colon \ZZ\text{-alg} \to \Grp \, , \; T(R) \coloneqq \Lambda_{\Kac}^\vee \otimes_\ZZ R^\times;
$$
alternatively, one can consider $T(R)$ as $\Hom_{\mathrm{grp}}(\Lambda_\Kac, R^\times)$:
$$
\Lambda^\vee \otimes R^\times \to \Hom_{\mathrm{grp}}(\Lambda_\Kac, R^\times) \, , \; h \otimes r \mapsto \left(\lambda \mapsto r^{\lambda(h)} \right).
$$
Since $\{v_i \}$ is a $\ZZ$-basis of $\Lambda_{\Kac}^\vee$ (see \Cref{Rem:AffineGCM}) and the rank is equal $2n-l$, the set 
$$
\{ r^{v_i} \mid r \in R^\times \, , \; 1 \leq i \leq 2n-l \}
$$
generates $T(R)$, cf.\ \cite[Example 7.55]{Marquis_2018}. If $\KK \in \{ \RR, \CC \}$, one can define an exponential map between $\h$ and $T$ as follows
$$
\KMExp \colon \h \to T(\KK) \, , \; (v_i \otimes r) \mapsto (v_i \otimes \exp(r)),
$$
where $\exp \colon \KK \to \KK^\times$ denotes the usual exponential map $\exp : (\RR,+) \to (\RR^\times, \cdot)$ or $\exp: (\CC,+) \to (\CC^\times, \cdot)$.
Since $T(\KK)$ is isomorphic to $(\KK^\times)^{2n-l}$ and $\h$ is isomorphic to $\KK^{2n-l}$, the exponential map $\KMExp$ corresponds to the natural exponential map $\exp \colon \KK^{2n-l} \to (\KK^\times)^{2n-l}$.

\bigskip
A \emph{basis of type $\D_{\Kac}^\AA$} (as in \cite[Section 2]{Tits_1987}, also \cite[Definition 7.77]{Marquis_2018}) is defined as a triple $\F \coloneqq \left( \G, (\p_i)_{i \in I}, \eta \right)$ where $\G \colon \ZZ\text{-alg} \to \Grp$ is a group functor, $(\p_i)_{i \in I}$ is a assortment of morphisms of functors $\p_i \colon \SL_2 \to \G$ and a morphism of functors $\eta \colon T \to \G$. One calls the group functor $\G$ \emph{Tits functor of type $\D$} if its satisfies the following axioms:

\begin{myitems}[label=\textbf{KMG.\arabic*}] 
    \item If $\KK$ is a field, $\G(\KK)$ is generated by $\p_i(\SL_2(\KK))$ and by $\eta(T(\KK))$.
    \item For every ring $R$, the homomorphism $\eta \colon T(R) \to \G(R)$ is injective.
    \item For $i \in I$ and $r \in R^\times$, one has 
    $$
    \p_i \begin{pmatrix} r & 0 \\ 0 & r^{-1} \end{pmatrix}  = \eta\left( r^{h_i} \right),
    $$
    where $r^{h_i}$ denotes the element $\lambda \mapsto r^{\langle \lambda, h_i \rangle}$ of the torus for each $h_i \in \Lambda^\vee$ and $\lambda \in \Lambda$.
    \item If $\iota \colon R \to \KK$ is an injective morphism of a ring $R$ in a field $\KK$, then $\G(\iota) \colon \G(R) \to \G(\KK)$ is injective.
    \item There is a homomorphism 
    $$\Ad \colon \G(\CC) \to \Aut([\g(\AA), \g(\AA)])$$
    whose kernel is contained in $\eta(\T_\Lambda(\CC))$, such that for $c \in \CC$ and $i \in I$,
    $$
    \Ad \left( \p_i  \begin{pmatrix} 1 & c \\ 0 & 1 \end{pmatrix}  \right) = \exp(\ad c e_i) , \;
    \Ad \left( \p_i  \begin{pmatrix} 1 & 0 \\ c & 1 \end{pmatrix}  \right) = \exp(\ad -c f_i) ,
    $$
    and for $t \in \T_\Lambda(\CC)$ and $i \in I$,
    $$ \Ad(\eta(t))(e_i) = t(c_i) \cdot e_i , \; \Ad(\eta(t))(f_i) = t(-c_i) \cdot f_i.$$
\end{myitems}
Note that, on fields, the Tits functor is unique up to isomorphism, see \cite[Theorem 7.82]{Marquis_2018}.

\begin{Convention} \label{convention3.5}
For the remainder of this article, let $\KK \in \{ \RR, \CC \}$ and denote by $G$ the evaluation of the Tits functor $\G$ with respect to the Kac--Moody root datum $\D_{\text{Kac}}^\AA$ on the field $\KK$; $G$ is called a \emph{split (minimal) Kac--Moody group over $\KK$ (of type $\AA$)} and $T$ its \emph{(standard) torus} of $G$.
\end{Convention}

\medskip

For any simple real root $\a_i \in \Delta^{re}$, $i \in I$, define the following subgroups of $G$ using the corresponding morphisms $\p_i$ from the definition of the Tits functor.
$$
U_{\a_i} \coloneqq \p_i \hspace{-2pt} \left( \left\lbrace \begin{pmatrix} 1 & r \\ 0 & 1 \end{pmatrix} \middle| \; r \in \KK \right\rbrace \right)
\; \text{ and } \; \; 
U_{-\a_i} \coloneqq \p_i \hspace{-2pt} \left( \left\lbrace \begin{pmatrix} 1 & 0 \\ -r & 1 \end{pmatrix} \middle| \; r \in \KK \right\rbrace \right).
$$
Now define the following group element for each simple root $\a_i$:
$$
\Tilde{s}_{i} \coloneqq \p_i \hspace{-2pt} \begin{pmatrix} 0 & 1 \\ -1 & 0 \end{pmatrix}.
$$
The group generated by these elements, $\EW \coloneqq \langle \Tilde{s}_{i} \mid i \in I \rangle$, is called the \emph{extended Weyl group}. A detailed discussion of the extended Weyl group and its properties can be found in \cite[Section 18]{GHK_2017}, \cite[Section 3.13]{FHHK_2017} or \cite[p. 172]{Kac-Peterson_1985}. Note that there is a canonical surjection from the extended Weyl group to the Weyl group given by 
\begin{equation}\label{Eq:Surjection_ExtWeyl-Weyl}
\EW \to \W \, , \; \Tilde{s}_{i} \mapsto r_i^\vee.    
\end{equation}
The kernel of the map is given by $\EW \cap T$, cf.\ \cite[Proposition 2.1]{Kac-Peterson_1985} or \cite[Section 3.13]{FHHK_2017}.\\
By the definition of a real root, each $\a \in \Delta^{re}$ can be written as $\a = w.\a_i$ for $\a_i \in \Pi$ and $w \in \W$. This translates, using the extended Weyl group, to the fact that for each real root $\a \in \Delta^{re}$ we can define the corresponding \emph{root (sub-)group} $U_\a$ as follows
$$
U_\a = U_{w.\a_i}\coloneqq \Tilde{w} U_{\a_i} \Tilde{w}^{-1},
$$
where $\Tilde{w} \in \EW$, see \cite[Section 3.13]{FHHK_2017}. 
Let $\a_i \in \Pi$, then one can describe an element of the root subgroup $U_{\a_i}$ as 
$$
x_{\a_i}(r) \coloneqq \p_i \hspace{-2pt} \begin{pmatrix} 1 & r \\ 0 & 1 \end{pmatrix} \in U_{\a_i},
$$
and an element of $U_{-\a_i}$ as
$$
x_{-\a_i}(r) \coloneqq \p_i \hspace{-2pt} \begin{pmatrix} 1 & 0 \\ -r & 1 \end{pmatrix} \in U_{-\a_i},
$$
for $r \in \KK$.\\
Recall from \Cref{Eq:RootLattice} that a root $\alpha$ is either positive or negative. Therefore, we can define $U_+ = \langle U_\a \mid \a \text{ is positive} \rangle$ as the subgroup only generated by the positive root subgroups. On the same way we define $U_-$ as the subgroup of $G$ only generated by the negative root subgroups.\\
Further define for $\a \in \DD$ the corresponding \emph{rank one subgroup} of $G$ as $G_\a \coloneqq \langle U_\a, U_{-\a} \rangle$. 
Note that by the axiom \textbf{KMG.3} the matrix 
$$
\begin{pmatrix}
    -1 & 0 \\
    0 & -1
\end{pmatrix}
$$
is mapped by $\p_i$ to $(-1)^{h_i}$, which is not the trivial element of the corresponding rank one subgroup, see \cite[Exercise 7.33]{Marquis_2018}, \cite[Chapter 8.2.1]{Remy_2002} and \cite[Chapter 8.4.2]{Remy_2002}. Therefore $G_\a \cong \SL_2(\KK)$, for details see \cite[Proposition 3.17]{FHHK_2017}.\\
Define the standard \emph{Borel subgroups} of $G$ by $B_\pm \coloneqq T \ltimes U_\pm$ (cf.\ \cite[Corollary 7.122]{Abramenko_2010}) and the conjugates of these groups are called Borel subgroups.

\begin{Convention}\label{Conv:Non-Degenerate-Tits-Functor}
In this article, we only consider Tits functors that are \emph{non-degenerate}, i.e.\ we assume that the morphisms $\p_i$, $i \in I$, from the definition of the Tits functor do not map $\SL_2(\KK)$, where $\KK$ is an arbitrary field, to the image of $U_+$ in the group given by the constructive Tits functor. The group given by the constructive Tits functor (cf. \cite[Definition 7.47]{Marquis_2018}) is a Kac--Moody group defined by generators and relations.\\
The assumption of non-degeneracy allows us to apply \cite[Theorem 7.82 (2)]{Marquis_2018}: Over fields, there is an isomorphism between the group given by the Tits functor and the group given by the constructive Tits functor. For details and a precise definition, see \cite[Definition 7.83]{Marquis_2018} and \cite[Definition 7.47]{Marquis_2018}.
\end{Convention}

For a split Kac--Moody group $G$ over $\KK$ of type $\AA$ one defines the \emph{Kac--Peterson topology} as the final group topology with respect to the maps $\p_i \colon \SL_2(\KK) \to G$ and $\eta \colon T \to G$, where $\SL_2(\KK)$ and $T$ are equipped with their natural Lie group topologies. For details of construction, see \cite[Section 7.5]{HKM_2013}.\\
This topology turns $G$ into a topological group with the following properties:
\begin{itemize}
    \item $G$ is Hausdorff (\cite[Proposition 7.21]{HKM_2013} and \cite[Definition and Remark 2.1]{HK_2021}). 
    \item The Kac--Peterson topology induces the unique Lie group topology on $T$ (\cite[Definition and Remark 2.1]{HK_2021} and \cite[Corollary 7.17]{HKM_2013}).
\end{itemize}
For more details, refer to \cite[Chapter 7]{HKM_2013}, \cite{HK_2021} and \cite{Harring_2020}. Also, many of these properties were announced by Kac and Peterson in \cite[Section~4G]{Kac-Peterson_1983} without proof.

\medskip

Now we prove \cite[Proposition 3.10]{FHHK_2017} in the context of this article. This allows us to endow the torus with a Lie group topology.

\begin{Proposition}\label{Prop:Torus-Exp-Topologisch}
Let $G$ be a split Kac--Moody group over $\KK$ of type $\AA$, where $\AA$ is an indecomposable, symmetrizable generalized Cartan matrix of size $n$ and rank $l$. Further denote by $\g(\AA)$ the corresponding Kac--Moody algebra. Equip $\h \cong \KK^{2n-l}$ with the vector space topology and $T \cong (\KK^\times)^{2n-l}$ with the Kac--Peterson subspace topology. Then $T$ is a Lie group homeomorphic to $(\KK^\times)^{2n-l}$ equipped with its standard topology.
\end{Proposition}

\begin{Proof}
The proof follows the same strategy as the proof of \cite[Proposition 3.10]{FHHK_2017}, namely to consider $G$ as a subgroup of a larger Kac--Moody group with trivial center and to apply \cite[Proposition 7.18]{HKM_2013}.\\
Recall that the root datum of $G$ is given by
$$
\D_{\Kac}^\AA = (I, \AA, \Lambda_{\Kac}, (c_i)_{i \in I}, (h_i)_{i \in I}),
$$
where a $\ZZ$-basis for $\Lambda_{\Kac}$ is given by $\{ u_j \}_{1 \leq j \leq 2n-l}$ and for $\Lambda_{\Kac}^\vee$ by $\{ v_j \}_{1 \leq j \leq 2n-l}$, such that $h_i = v_i$ for $1 \leq i \leq n$.\\
Now let $\mathbb{B} \in \ZZ^{2n-l \times 2n-l}$ be a generalized Cartan matrix with rank $2n-l$, such that $\AA$ is included in $\mathbb{B}$ as a principal submatrix. This can be done in a similar way to expanding a non-invertible generalized Cartan matrix to obtain the extended Cartan algebra, see \cite[Proposition 1.1]{Kac_1990} or \cite[Section 3.5]{Marquis_2018} (or \Cref{Rem:AffineGCM}). Consider now the corresponding root datum 
$$
\widetilde{\D} = (J, \mathbb{B}, \Lambda_{\Kac}, (u_i)_{i \in J}, (v_i)_{i \in J}),
$$
where $J = \{ 1, \ldots, 2n-l \}$ and by construction we have $\langle u_j, v_i \rangle = b_{ij}$ and that $\langle u_j , v_i \rangle = a_{ji}$ for $i, j \in I$. Denote by $\widetilde{G}$ the Kac--Moody group associated to $\widetilde{\D}$.\\
One can construct an embedding of $G$ into $\widetilde{G}$ as follows: For any $k \in \KK$, $r \in \KK^\times$ the embedding is given by
\begin{align*}
    \pi \colon G &\to \widetilde{G} \\
    x_{\pm \a_i}(k) &\mapsto x_{\pm \a_i}(k) \; \;  (1 \leq i \leq n) \\
    r^{v_i} &\mapsto r^{v_i} \; \; (i \in J) ,
\end{align*}
where $\a_i = c_i \otimes 1 $ for $i \in I$. Compare this to \cite[Proposition B.5]{HK_2021} or \cite[Appendix C]{Diss}, where it is done for arbitrary split Kac--Moody groups. 
As in \cite[Proposition 3.5.2]{Harring_2020}, the kernel of $\pi$ is a finite central subgroup, possibly trivial.\\
Both are Kac--Moody groups and since we assumed that the Tits functor is non-degenerate, it is already clear that the relations of a Kac--Moody group (see \cite[Definition 7.47]{Marquis_2018}) are preserved under the mapping $\pi$. Moreover, both groups can be equipped with the Kac--Peterson topology, and by construction is $G$ a subgroup of $\widetilde{G}$.
Based on the construction, $\mathbb{B}$ is invertible and thus $\widetilde{\D}$ has finite center (see \cite[Proposition 3.5.2]{Harring_2020}). This allows us to apply \cite[Proposition 7.18]{HKM_2013}. It follows from this proposition that the torus of $\widetilde{G}$ endowed with the Kac--Peterson topology is homeomorphic to $(\KK^\times)^{2n-l}/\ker(\pi)$ equipped with its standard Lie topology. Note, that $\ker(\pi)$ is a finite central subgroup by \cite[Proposition 3.5.2]{Harring_2020}. We conclude that the torus of $G$ is homeomorphic to $(\KK^\times)^{2n-l}$ as a Lie group, covering the Lie group $(\KK^\times)^{2n-l}/\ker(\pi)$.
\qed
\end{Proof}


\subsection{Twin Buildings}\label{Sec:Buildings}
This section provides a short overview of the basic notions of a (twin) building. For more details see \cite{Abramenko_2010}. Let
$$
\W \coloneqq \langle r_1 , \ldots, r_n \mid (r_i r_j)^{m_{ij}} = e \rangle,
$$
be a \emph{Coxeter group}, where $(m_{ij})_{1 \leq i, j \leq n} \subset \ZZ_{>0} \cup \{ \infty \} $ such that 
\begin{itemize}
    \item $m_{ii} = 1$,
    \item $m_{ij} = m_{ji} \geq 2 $ for $i \neq j$.
\end{itemize}
The pair $(\W, S)$ with generating set $S = \{ r_1, ..., r_n \}$ is a \emph{Coxeter system}. 

Note that any Weyl group of a Kac--Moody algebra forms a Coxeter system and the order of the elements are in relation with the entries of the generalized Cartan matrix, see \cite[Proposition 4.22]{Marquis_2018}. Not all Coxeter groups appear as Weyl groups, for instance the finite Coxeter group of type $H_3$ cannot be a Weyl group.

\begin{Definition}[cf.\ Definition 5.1 from \cite{Abramenko_2010}]
Given a Coxeter system $(\W,S)$,  a \emph{building $\Delta$ of type $(\W,S)$} is a pair $(\C, \delta)$, where $\C$ is a non-empty set whose elements are called \emph{chambers} together with a \emph{Weyl distance function} $ \delta \colon \C \times \C \to \W$ such that for all $C, D , E \in \C$ the following properties hold:
\begin{itemize}
    \item $\delta(C,D) = e \Leftrightarrow C = D$;
    \item if $\delta(C,D) = w \in \W$ and $\delta(E,C) = s \in S$, then $\delta(E,D) \in \{ sw, w \}$. Further, if $l(sw) = l(w)+1$, then $\delta(E,D) = sw$, where $l(.)$ denotes the length of the reduced word in $\W$;
    \item if $\delta(C,D) = w \in \W$, then for any $s \in S$ there exists a chamber $F \in \C$ such that $\delta(F,C) = s$ and $\delta(F,D) = sw$.
\end{itemize}
\end{Definition}

In order to associate a building to a group, we need the notion of a \emph{$BN$-pair}, which implies the \emph{Bruhat decomposition}, for details see \cite[Section 6.2.6]{Abramenko_2010}. It should also be noted that in this context of $BN$-pairs the abstract Coxeter group as we defined above becomes a Weyl group in the sense of a Kac--Moody group (for the definition see the end of \Cref{Sec:KM_Algebra}).\\
Now let $G$ be a group with subgroup $B$, $(\W, S)$ a Coxeter system and $C \colon \W \to B \backslash G / B$ the bijection given by 
\begin{equation}\label{Eq:C(w)}
C(w) = B \Tilde{w} B,
\end{equation}
where $\Tilde{w}$ is a representative of the preimage of $w$ in the extended Weyl group under the surjection $\EW \to \W$, see \Cref{Eq:Surjection_ExtWeyl-Weyl}. Note that $\W$ stands for a Weyl group in this context. This explicit description allows us to define a decomposition, as we will see later.\\
Next assume that for all $s \in S$ and $w \in \W$ the property $C(sw) \subseteq C(s)C(w) \subseteq C(sw) \cup C(w)$ is satisfied. This property allows us to define a building using the definition from above: consider the coset $G/B$ as the set of chambers and define the Weyl distance as follows 
\begin{align*}
    \delta \colon G/B \times G/B &\to B \backslash G / B \to \W \\
    (gB, hB) &\mapsto Bg^{-1}hB \mapsto C^{-1}(Bg^{-1}hB).
\end{align*}

This leads to the following equivalence, which allows to compute the Weyl distance
$$
\delta(gB, hB) = w \Longleftrightarrow g^{-1}h \in C(w).
$$
This construction satisfies the definition from above and thus leads to a building. In the following we denote by $\Delta_G$ a building which is associated to a group $G$.

The notion of a  $BN$-pair guarantees the existence of a Bruhat decomposition and, thus, a building to act on.

\begin{Definition}[cf.\ Definition 6.55 from \cite{Abramenko_2010}]
Let $G$ be a group with subgroups $B,N$. The pair $(B,N)$ is called \emph{$BN$-pair} for $G$ if the following is satisfied:
\begin{itemize}
    \item $G$ is generated by $B$ and $N$.
    \item The intersection $N \cap B$ is normal in $N$.
    \item The quotient $N/T$ admits a set of generating involutions $S$ such that
    \begin{itemize}
        \item for all $w \in N/T$ and $s \in S$ we have $sBw \subseteq BswB \cup BwB$, and
        \item $sBs^{-1}$ is not a subgroup of $B$ for all $s \in S$.
    \end{itemize}
\end{itemize}
\end{Definition}

We denote by $T$ the intersection $N \cap B$ and by $\W$ the quotient $N/T$. If a group $G$ has a $BN$-pair $(B,N)$, then the \emph{Bruhat decomposition} holds:
$$
G = \bigsqcup_{w \in \W}BwB,
$$
see \cite[Section 6.2.6]{Abramenko_2010}.
This decomposition combined with \Cref{Eq:C(w)} yield the desired bijection
\begin{align*}    
C \colon \W &\to B \backslash G / B \\
w &\mapsto C(w) \coloneqq BwB.
\end{align*}

\medskip

Given two buildings $\Delta_+ = (\C_+, \delta_+)$ and $\Delta_- = (\C_-, \delta_-)$ of type $(\W,S)$, where $\C_- \cap \C_+ = \emptyset$. A \emph{twin building of type $(\W,S)$} with \emph{codistance function $\delta^\star$} is a triple
$$
\Delta_\pm \coloneqq (\C_+, \C_-, \delta^\star)
$$
where the map
$$
\delta^\star \colon (\C_+ \times \C_-) \cup (\C_- \times \C_+) \to \W
$$
satisfies for any $C \in \C_{\epsilon}$, $D \in \C_{-\epsilon}$ and $\epsilon \in \{ +, - \}$ 
\begin{itemize}
    \item $\delta^\star(C,D) = \delta^\star(D,C)^{-1}$;
    \item if $\delta^\star(C,D) = w \in \W$ and if $\delta_{\epsilon}(E,C) = s \in S$ for $E \in \C_{\epsilon}$ and $l(sw) < l(w)$, then $\delta^\star(E,D) = sw$;
    \item for any $s \in S$ there exists a chamber $E \in \C_{\epsilon}$ such that $\delta_{\epsilon}(E,C) = s$ and $\delta^\star(E,D) =sw$.
\end{itemize}
 For more details, see \cite[Definition 5.133]{Abramenko_2010}.

By \cite[Definition 6.78]{Abramenko_2010} a building can be constructed from a group $G$ if it has a \emph{twin $BN$-pair} $(B_+, B_-, N)$, which is a $BN$-pair with respect to the subgroup $B_+$ and $B_-$. Similar to the case for a single building, a twin $BN$-pair guarantees the \emph{Birkhoff decomposition} (cf.\ \cite[Proposition 6.81]{Abramenko_2010}):
$$
G = \bigsqcup_{w \in \W} B_\epsilon w B_{-\epsilon},
$$
where $\epsilon \in \{ +, - \}$. Again, as in the case of a single building, this decomposition allows to define a codistance function on the set of chambers $G/B_\pm$. Note that $B_+, B_-$ and $N$ are subgroups of $G$ satisfying some properties, for details see \cite[Definition 6.78]{Abramenko_2010}.

Examples of groups that possess a (twin) $BN$-pair are the general linear group (\cite[Chapter 6.5]{Abramenko_2010}), the symplectic group (\cite[Chapter 6.6]{Abramenko_2010}), isotropic orthogonal groups (\cite[Chapter 6.7]{Abramenko_2010}) and isotropic unitary groups (\cite[Chapter 6.8]{Abramenko_2010}). 

\begin{Remark}\label{Rem:RGD}
The notion of an RGD system is very useful for studying Kac--Moody group. For instance, an RGD system allows one to construct a $BN$-pair and, thus, associate a building to a group.
Below we recall the definition of an RGD system, see also \cite[Section 8.6.1]{Abramenko_2010} and \cite[Definition 7.82]{Abramenko_2010}. In this article we will only need a few properties of an RGD system, which is why we will not go into the technical details here.\\
Let $G$ be a Kac--Moody group over $\KK$ of type $\AA$. Then, by \cite[Theorem 7.69]{Marquis_2018}, the triple $(G, (U_\a)_{\a \in \DD}, T)$ represents a (linear) RGD system. This means that this triple satisfies the following properties
\begin{myitems}[label=\textbf{RGD.\arabic*}]
\setcounter{myitemsi}{-1}
    \item For all $\a \in \DD$ the root group $U_\a$ is non-trivial.
    \item For all $\a, \beta \in \DD$ with $\a \neq \beta$ such that the set $\{ \a, \beta \}$ is prenilpotent, we have $[ U_\a, U_\beta ] \leq U_{(\a, \beta)}$, where $(\a, \beta) \coloneqq  [ \a, \beta ] \setminus \{ \a, \beta \}$. For details on the notation see \cite[Definition 8.41]{Abramenko_2010}.
    \item For every $s \in S$ there exists a function $m \colon U_{\a_s}\setminus \{e \} \to G$, such that for all $u \in U_{\a_s} \setminus \{ e \}$ and $\a \in \DD$ we have 
    \begin{align*}
        &m(u)\in U_{-\a_s} u U_{-\a_s} \\
        &m(u)U_\a m(u)^{-1} = U_{s\a}.
    \end{align*}
    Furthermore, for all $u,v \in U_{\a_s}\setminus \{ e \}$ follows $m(u)^{-1}m(v) \in  T$.
    \item For all $s \in S$ $U_{-\a_s}$ is not a subgroup of $U_+$.
    \item $G = T \langle U_\a \mid \a \in \DD \rangle$
    \item $ T \leq \bigcap_{\a \in \DD} N_G(U_\a)$
\end{myitems}

Following \cite[Chapter 8.8]{Abramenko_2010} one can define subgroups $B_\pm$ by
$$
B_\pm \coloneqq TU_\pm,
$$
see \cite[Remark 7.87]{Abramenko_2010}. Further in \cite[Chapter 8.8]{Abramenko_2010} it is shown, that $ N_G(U_\pm) = B_\pm $ and that $B_+ \cap B_- = T$ (\cite[Corollary 8.78]{Abramenko_2010}). Altogether, by \cite[Theorem 8.80]{Abramenko_2010} an RGD system leads to a $BN$-pair with bijection 
$$
N_G(T)/T \cong \W.
$$
By \cite[Corollary 8.79]{Abramenko_2010} follows that the last RGD axiom is an equality,
$$
T = \bigcap_{\a \in \DD} N_G(U_\a).
$$
\end{Remark}

\subsection{A Characteristic Subgroup}\label{Sec:G_Dagger}
A very useful subgroup for studying symmetric spaces associated to Kac--Moody group is the \emph{centered subgroup}, that is the subgroup generated only by the root groups.

\begin{Definition}\label{Def:G1}
Let $G$ be a split Kac--Moody group over $\KK$ of type $\AA$. Denote by $G^\dagger < G$ the subgroup of $G$ generated by the root subgroups of $G$, i.e.
$$ 
G^\dagger \coloneqq \langle U_\a \mid \a \in \Delta^{re} \rangle.
$$
This subgroup is called the \emph{centered subgroup}.
\end{Definition}

We can deduce that $G^\dagger$ is a normal subgroup of $G$: By \textbf{RGD.4} we can write 
$$
G = T G^\dagger.
$$
To deduce, that $G^\dagger$ is normal, consider an element $g \in G$ which we can also write as $g = t \widetilde{g}$. Now compute
$$
g G^\dagger g^{-1} = t \widetilde{g} G^\dagger \widetilde{g}^{-1} t^{-1}.
$$
By definition of $G^\dagger$ and by \textbf{RGD.5} follows, that $g G^\dagger g^{-1} \subseteq G^\dagger$. From \cite[Corollary 8.79]{Abramenko_2010} follows that the torus of a Kac--Moody group can be described as the intersection of all normalizers of the root groups. We use this fact to define a torus for $G^\dagger$:
$$
T^\dagger \coloneqq \bigcap_{\a \in \DD} N_{G^\dagger}(U_\a) = T \cap G^\dagger.
$$
Since the RGD axioms are mainly using root groups, we can deduce that $(G^\dagger, (U_\a)_{\a \in \DD}, T^\dagger)$ forms an RGD system. Moreover, according to \cite[Proposition 8.82]{Abramenko_2010}, the kernel of the action of $G$ on the associated twin building $\Delta_\pm$ is contained in the torus $T$, so that 
$$
N_{G^\dagger}(T^\dagger) / T^\dagger \cong \W.
$$
Note that we use the same arguments here as in a similar statement which can be found under the definition of an RGD system (see \Cref{Rem:RGD}). By \cite[Theorem 8.80]{Abramenko_2010}, an RGD system leads to a $BN$-pair with this bijection.\\
We conclude this paragraph with the observation that the type of the twin building belonging to $G^\dagger$ is the same as the type of the twin building belonging to $G$.

\begin{Lemma}\label{Lem:G1-is-charakteristisch}
Let $G$ be a split Kac--Moody group over $\KK$ of type $\AA$. Then $G^\dagger$ is the derived subgroup of $G$, in particular a characteristic subgroup, i.e.
$$
f(G^\dagger) \subseteq G^\dagger 
$$
for all $f \in \Aut(G)$.
\end{Lemma}

\begin{Proof}
By \textbf{RGD.4} one has $G = T G^\dagger$ (see \cite[Chapter 8]{Abramenko_2010}). Moreover, $G^\dagger$ is a normal subgroup of $G$, since $$T = \bigcap_{\a \in \Delta^{re}}{N_G(U_\a)},$$ cf.\ \cite[Corollary 8.79]{Abramenko_2010}. Since the natural map $T \mapsto G/G^\dagger$ is surjective, the quotient $G/G^\dagger$ is abelian and thus $[ G,G ] \subseteq G^\dagger$.\\
On the other hand, by definition $G^\dagger$ is generated by the root subgroups, in particular $G^\dagger$ is generated by the standard rank one subgroups 
$$
G_i = \langle U_{\a_i} , U_{-\a_i} \rangle \cong \SL_2(\KK).
$$
Due to the fact that $\SL_2(\KK)$ is perfect for $k \in \lbrace \RR, \CC \rbrace$, one deduces that $G^\dagger$ is perfect. Thus 
$$
[G,G] \subseteq G^\dagger = [G^\dagger, G^\dagger] \subseteq [G,G]
$$
which finishes the proof.
\qed
\end{Proof}
Note, that there is also an action of $\Aut(G)$ on $G^\dagger$ given by the homomorphism $\Aut(G) \to \Aut(G^\dagger)$ which  restricts any $f \in \Aut(G)$ on $G^\dagger$.

\begin{Remark}\label{Rem:G1-is-KM}
Note that this centered subgroup $G^\dagger$ is the Kac-Moody group used by the authors in \cite{FHHK_2017} to develop their theory of Kac--Moody symmetric spaces.
\end{Remark}

\subsection{The Chevalley Involution}\label{Sec:Chevalley}
The first goal of this article is to construct a symmetric space with respect to a Kac--Moody group. For this purpose, it is useful and important to define an involution. Therefore, let $G$ be a Kac--Moody group over $\KK \in \{ \RR, \CC \}$ of type $\AA$. Then $\Theta \colon G \to G$ denotes a continuous group involution which is given for $\KK = \CC$ by
$$
\Theta(x_{\a}(r)) = x_{-\a}(\overline{r}) \, , \; \Theta(r^{v_i}) = \overline{r}^{-v_i},
$$
for $r \in \CC$ and $\a \in \DD$. From this definition one can easily derive the following
$$
\Theta(t) = \overline{t}^{-1}
$$
for all $t \in T$.\\
The map $\Theta$ is continuous with respect to the Kac--Peterson topology on $G$. Since the complex conjugation acts trivially on the reals, it can be omitted in the case $\KK = \RR$. By \cite[§8.2]{Caprace_2009} this involution can be obtained by lifting the Chevalley involution of a Kac--Moody algebra $g(\AA)$ (cf.\ \cite[§1.3]{Kac_1990}) to a continuous involution on $G$. In particular, this means that $\Theta$ maps $U_\a$ to $U_{-\a}$ and stabilizes the torus $T$.\\
Denote by $K = G^\Theta$ the set of fixed points of $\Theta$ and define by $\tau \colon G \to G$, $g \mapsto g \Theta(g)^{-1}$ the \emph{twist map} of $G$. In the case that $\KK \in \{ \RR, \CC \}$, it can be derived using the definition of $\Theta$ that
$$
\tau(t) =
\begin{cases}
    t^2 & \text{for } \KK = \RR \\
    |t|^2 & \text{for } \KK = \CC,
\end{cases}
$$
for all $t \in T$, or simply $\tau(t) = |t|^2$.\\
In the following lemma, we define useful subgroups and analyze their properties. The statements are based on \cite[Lemma 3.24]{FHHK_2017} and \cite[Lemma 3.26]{FHHK_2017}.

\begin{Lemma}\label{Lem:Eigenschaften-A_M}
Let $G$ be a split Kac--Moody group over $\KK$ of type $\AA$, where $\KK \in \{ \RR, \CC \}$. Define the subgroups $A \coloneqq \tau(T)$ and $M \coloneqq K \cap T$. Then we have $T \cong M \times A$ and the following hold:
\begin{myitems}[label=\textbf{(A\arabic*)}, wide, labelindent=0pt]
    \item \label{itm:A0} For $g,h \in G$ one has $hK = gK \Leftrightarrow \tau(g) = \tau(h)$ and $\tau^{-1}(e) = K$.
    \item \label{itm:A1} The subgroup $A$ is isomorphic to $(\RR_{>0})^{2n-l}$.
    \item \label{itm:A2} $M \cap A = \{ e \} = K \cap A$.
    \item \label{itm:A3} For $\KK = \CC$ the subgroup $M$ is isomorphic to $(\SS^1)^{2n-l}$ and is the \emph{maximal compact subgroup}; in the case $\KK = \RR$, the $M$ is the \emph{torsion subgroup}, i.e.\ unique maximal finite subgroup of order $2^{2n-l}$.
    \item \label{itm:A4} $K \cap B_\pm = M$.
\end{myitems}
\end{Lemma}

\begin{proof}
For \ref{itm:A0} let $g,h \in G$ and it follows from $gK = hK$ that there is a $k \in K$ such that $g = hk$. For one direction one calculates
$$
gK = hK \Rightarrow \tau(g) = \tau(hk) = \tau(h),
$$
and for the other direction 
\begin{align*}
    \tau(g) = \tau(h) &\Rightarrow g \Theta(g)^{-1} = h \Theta(h)^{-1} \\
    &\Rightarrow \Theta(h^{-1}g) = h^{-1}g \\
    &\Rightarrow h^{-1}g \in K \\
    &\Rightarrow gK = hK.
\end{align*}
To see that $\tau^{-1}(e) = K$, let $g \in G$ be an element with $\tau(g) = e$. Then calculate
$$
\tau(g) = e \Leftrightarrow \Theta(g) = g \Leftrightarrow g \in K.
$$
The statement of \ref{itm:A1} follows directly from the definition of the maps, see the paragraph above. Concretely, this shows that $A = \tau(T) \cong (\RR_{>0})^{2n-l}$, because the dimension of $T$ is $2n-l$.\\
Next, we want to check \ref{itm:A2}. Let $g \in K \cap A \subseteq T$. Then by definition of $K$ and $\Theta$ follows that
$$
g = \Theta(g) = \overline{g}^{-1}.
$$
Thus $g \overline{g} = e$ and from the fact that $A = \tau(T) \cong (\RR_{>0})^{2n-l}$ one concludes that $g = e$. In particular, the statement that $M \cap A = \{ e \}$ is true follows then directly from $M \subset K$.\\
The statement of \ref{itm:A3} follows from \textbf{(A3)} and the polar decomposition
$$
T \cong 
\begin{cases}
    (\CC^\times)^{2n-l} \cong (\RR_{>0})^{2n-l} \times (\SS^1)^{2n-l} & \text{for } \KK = \CC \\
    (\RR^\times)^{2n-l} \cong (\RR_{>0})^{2n-l} \times \{ \pm 1 \}^{2n-l} & \text{for } \KK = \RR.
\end{cases}
$$
Finally, we prove the statement of \ref{itm:A4}. Therefore recall, that by the definition of $B_\pm$ one has $\Theta(B_\pm) = B_\mp$ and that $T = B_+ \cap B_-$. Note, that this is based on the fact that $U_+ \cap U_- = \{ e \}$. Thus we have $K \cap B_\pm = K \cap T$ = M, by definition.
\end{proof}

\subsection{The Iwasawa Decomposition}\label{Sec:Iwaswa}
In the theory of symmetric spaces, it is necessary to define flats because they are the most important objects for understanding actions and causal boundaries. To do so, one needs the topological Iwasawa decomposition, i.e.\ that the multiplication induces a homeomorphism 
$$
m_\pm \colon U_\pm \times A \times K \to G .
$$
In general, this is known for centered Kac--Moody groups, i.e.\ Kac--Moody groups generated only by their root groups without any additional torus, see \cite[Section 3]{Kac-Peterson_1985}.\\
First we will prove the algebraic statement that for a Kac--Moody group $G$ there is a decomposition $G = U_\pm A K$. After establishing the algebraic Iwasawa decomposition, we will use the same strategy as in the proof of \cite[Theorem 3.31 (ii)]{FHHK_2017} to prove the topological assertion for Kac--Moody groups over $\KK$ of arbitrary type $\AA$.\\
The following section is mainly based on \cite[Section 3.30]{FHHK_2017}. Denote by $\KK$ the field of real or complex numbers and with $\AA$ a symmetrizable, indecomposable generalized Cartan matrix.

\begin{Proposition}\label{Prop:AlgIwasawa}
Let $G$ be a Kac--Moody group over $\KK$ and of arbitrary type $\AA$. Then $G$ decomposes into
$$
G = U_\pm A K = K A U_\pm.
$$
Moreover, the maps,
$$
    m_\pm \colon U_\pm \times A \times K \to G
$$
given by multiplication, are bijective.
\end{Proposition}

\begin{Proof}
By \Cref{Rem:RGD} and \Cref{Sec:G_Dagger}, one can decompose $G$ into $T G^\dagger$, where $T$ is the torus and $G^\dagger$ is the centered subgroup of $G$. Since $\KK$ plays no role in the algebraic decomposition, and $G^\dagger$ is a centered Kac--Moody group, one can use \cite[Theorem 3.31 (ii)]{FHHK_2017}. Therefore denote by $K^\dagger$ the fixed point set $(G^\dagger)^\Theta$ and by $B_\pm^\dagger$ the (standard) Borel subgroup $B_\pm^\dagger = T^\dagger \ltimes U_\pm$. Then there is the decomposition
$$
G^\dagger = B_\pm^\dagger K_\pm^\dagger.
$$
Recall that $T = M A$, $M \subseteq K$ by \Cref{Lem:Eigenschaften-A_M} and that $T$ normalizes $U_\pm$. Also denote with $A^\dagger = \tau(T^\dagger)$ and with $M^\dagger = K \cap T$. All properties of \Cref{Lem:Eigenschaften-A_M} are also true for the corresponding subgroups in the context of $G^\dagger$. It follows
\begin{align*}
    G &= T G^\dagger \\
    &= T B_\pm^\dagger K^{\dagger} \\
    &= T U_\pm A^\dagger M^\dagger K^\dagger \\
    &= U_\pm (MA)A^\dagger K^\dagger\\
    &= U_\pm MA K^\dagger \\
    &= U_\pm AM K^\dagger \subset U_\pm AK
\end{align*}
since $M, K^\dagger \subset K$ (see \Cref{Lem:Eigenschaften-A_M} \ref{itm:A4}). The other order of the decomposition holds by symmetry (e.g., via group inversion).

\smallskip

To see that these maps $m_\pm$ are bijective, we need to check injectivity. To do this, we assume that $u_1 a_1 k_1 = u_2 a_2 k_2$. By reordering follows,
$$
K \ni k_1 k_2^{-1} = a_1^{-1} u_1^{-1} u_2 a_2 \in A U_\pm A \subseteq B_\pm.
$$
Hence $a_1^{-1} u_1^{-1} u_2 a_2 \in K$ and $B_\pm \cap K = M \subseteq T$, see \Cref{Lem:Eigenschaften-A_M}. And as a result
$$
u_1^{-1} u_2 \in U_\pm \cap T = \lbrace 1 \rbrace \Longrightarrow u_1 = u_2,
$$
(cf.\ \cite[Lemma 7.62]{Abramenko_2010}) and from this, that
$$
a_1^{-1} a_2 \in M \cap A = \lbrace 1 \rbrace \Longrightarrow a_1 = a_2.
$$
Therefore $k_1 = k_2$ and thus $m_\pm $ is injective.
\qed
\end{Proof}

In the following we consider only Kac--Moody groups of non-spherical type $\AA$. In the case where $\AA$ is of spherical type, the Kac--Moody group of type $\AA$ is a finite-dimensional Lie group and the topological Iwasawa decomposition is a well-known fact, see for example \cite[Theorem 6.46]{Knapp_2002}. Further, there are also proofs of the topological Iwasawa decomposition with respect to other groups, e.g.\ for $\operatorname{GL}_n(\RR)$ (\cite[Proposition 3.12]{Platonov_1993}) or in greater generality for reductive $\RR$-groups (\cite[Theorem 3.9]{Platonov_1993}).

\smallskip

Now for the topological version, the idea is to use the same strategy as in the proof of \cite[Theorem 3.31 (ii)]{FHHK_2017}, to be precise we will construct the inverse functions and show that they are open. The fact that we can use the same strategy as in \cite{FHHK_2017} is based on the statement that every Kac--Moody group has a strong topological building, as shown by Grüning and the first-named author of this article in \cite[Appendix B]{HK_2021} for centered Kac--Moody groups. The very technical generalization for arbitrary split Kac--Moody groups can be found in \cite[Appendix C]{Diss}. For an introduction to topological buildings see \cite[Section 3.1]{HKM_2013} and for strong topological buildings in particular see \cite[Definition 3.21]{HKM_2013}.

\begin{Lemma}\label{Lemma:Fasserung}
Let $G$ be a Kac--Moody group over $\KK$ and of type $\AA$, where $\AA$ is non-spherical. Then the fibration 
$$
p_\pm \colon \nicefrac{G}{AU_\pm} \to \nicefrac{G}{B_\pm}
$$
is trivial with compact fiber $M$, i.e., a topological direct product of the fiber and the base space.
\end{Lemma}

\begin{Proof}
First, note that $G$ leads to a non-spherical strong topological twin building, according to a result by Grüning and the first-named author of this article, see \cite[Corollary B.8]{HK_2021} (or \cite[Appendix C]{Diss}). Furthermore, note that $T \cong A \times M$ as topological groups and that $U_\mp \times T \times U_\pm$ is open. In particular the map $U_\mp \times T \times U_\pm \to B_\mp B_\pm$ is a homeomorphism (cf.\ \cite[Corollary B.7]{HK_2021}) and $B_\mp B_\pm$ is open (cf.\ \cite[Lemma 6.1]{HKM_2013}). Hence 
$$
\{ e \} \times M \times \{ e \} \to U_\mp \times T \times U_\pm \to U_\mp \times A \times U_\pm
$$
provides a local trivialization of the fibration 
$$ 
M \to \nicefrac{G}{AU_\pm} \to \nicefrac{G}{B_\pm},
$$
in particular, every point in $G/AU_\pm$ has a neighborhood that looks like an open set in $G/B_\pm$ times $M$.\\
Since the non-spherical twin building associated to $G$ is a strong topological twin building, it is contractible (cf.\ \cite[Theorem 5.13]{HKM_2013}) and thus the fibration $p_\pm$ is trivial.
\qed
\end{Proof}

For the next lemma we need the notion of \emph{$k_\omega$-spaces}: a Hausdorff topological space $X$ is called a $k_\omega$-space, if there is a family of compact subspaces $(X_n)_{n \in \NN}$, called $k_\omega$-sequence, such that
$$
X = \bigcup_{n \in \NN} X_n
$$
and any open subset $Y \subseteq X$ is open in $X$ if and only if for every $n \in \NN$ $Y \cap X_n$ is open in $X_n$. Call the pair $(X, (X_n)_{n \in \NN})$ \emph{$k_\omega$-pair}. For more details about $k_\omega$-spaces see \cite{k-Omega} or \cite[Definition 3.33 et seq.]{FHHK_2017}. Note that a Kac--Moody group over $\KK$ equipped with the Kac--Peterson topology is a $k_\omega$-group, for a justification and details see \cite[Proposition 7.10]{HKM_2013} and \cite[Definition and Remark 2.1]{HK_2021}.\\
Here we recall a lemma about $k_\omega$-spaces for the reader's convenience that we use in order to derive  the topological Iwasawa decomposition.

\begin{Lemma}[Lemma 3.34 from \cite{FHHK_2017}]\label{Lem:k-Omega-1}
Let $(X,(X_n))$ and $(Y, (Y_n))$ be $k_\omega$-pairs and let $f \colon X \to Y$ be a continuous bijection such that for all $n \in \NN$ there exists a $m \in \NN$ such that $Y_n \subset f(X_m)$.\\
Then $f(X_n)$ is a $k_\omega$-sequence for $Y$ and $f$ is a homeomorphism.
\end{Lemma}


Now, we can prove the following, which is a version of \cite[Proposition 3.36]{FHHK_2017} adapted to the context of this article.

\begin{Lemma}\label{Lemma:HilfeIwasawa}
Let $G$ be a Kac--Moody group over $\KK$ of type $\AA$, where $\AA$ is non-spherical and define the maps
\begin{align*}
\iota_\pm \colon K &\to G/AU_\pm \\
k &\mapsto k AU_\pm.
\end{align*}
These maps are homeomorphisms.
\end{Lemma}

\begin{Proof}
The goal is to use \cref{Lem:k-Omega-1}. According to the algebraic Iwasawa decomposition (cf.\ \Cref{Prop:AlgIwasawa}), the maps $\iota_\pm$ are bijections, since multiplication is continuous, these maps are also continuous.\\
Next, we recall that $G$ acts strongly transitively on the associated twin building (cf.\ \Cref{Rem:RGD}) and we can therefore use the Birkhoff decomposition according to \cite[Proposition 6.75]{Abramenko_2010}. Based on the Birkhoff decomposition we can define the following:
$$
G_l^\pm \coloneqq \bigcup_{w \in \W, \, l(w) \leq l} B_\mp w B_\pm,
$$
and denote with $\widetilde{D}_{l,\pm}$ the image of $G_l^\pm$ in $G/AU_\pm$ and with $D_{l,\pm}$ the image of $G_l^\pm$ in $G/B_\pm$. By \cite[Corollary 3.13]{HKM_2013} we have 
$$
\left( \nicefrac{G}{B_\pm} , \left( D_{l,\pm} \right) \right)
$$
is a $k_\omega$-pair. Note that the preimage of $D_{l, \pm}$ under the fibration $p_\pm$ is exactly $\widetilde{D}_{l, \pm}$ and that the fibration by \Cref{Lemma:Fasserung} is trivial with compact fiber $M$, so we observe that
$$
\left( \nicefrac{G}{AU_\pm} , \left( \widetilde{D}_{l,\pm} \right) \right)
$$
is also a $k_\omega$-pair.\\
In order to apply \Cref{Lem:k-Omega-1}, it remains to prove that $(K, (K_l^\pm))$ is a $k_\omega$-pair, where $K_l^\pm \coloneqq K \cap G_l^\pm$. Therefore by \cite[Corollary 7.11]{HKM_2013} one knows that
$$
G = \varinjlim G_l^\pm
$$
and in particular that
$$
K = \varinjlim K_l^\pm.
$$
To conclude that $(K, (K_l^\pm))$ is a $k_\omega$-pair, we must show that the subsets $(K_l^\pm)$ are compact. To this end, recall that $G$ can be written as $K B_\pm$ using \Cref{Prop:AlgIwasawa}, and hence $K$ acts chamber transitively on the twin building $\Delta_\pm$ associated to $G$.
In particular $$\Delta_\pm \quad\quad \text{and} \quad\quad K/(K \cap B_\pm) = K/(K \cap T)$$ are isomorphic as $K$-sets.
Since $G_l^\pm$ is a union of $B_\pm$-double cosets and, moreover, is invariant under the Cartan-Chevalley involution, we have 
$$
G_l^\pm = K_l^\pm  B_\pm.
$$

Denote by $\C_\pm$ the set of chambers of $\Delta_\pm$, then by using the Bruhat order (\cite[Exercise 3.59 (c)]{Abramenko_2010}) we can define the set 
$$
E_w(C) = \{ D \in \C_\pm \mid \delta_\pm(C,D) = w \} \subseteq \Delta_\pm,
$$
where $C \in \C_\pm$ denotes a chamber in $\Delta_\pm$, $w \in \W$ an element of the Weyl group and $\delta_\pm$ the Weyl distance in relation to the two halves of the twin building. The above isomorphism of $K$-sets between $\Delta_\pm$ and  $K/(K \cap B_\pm) = K/(K \cap T)$ induces a bijection between 
$$
\bigcup_{w \in \W \, , \; l(w) \leq l} E_w(C) \quad \quad \text{and} \quad \quad K_l^\pm/ (K_l^\pm \cap B_\pm).
$$
This is a finite union of compact sets, see \cite[Corollary 3.10]{HKM_2013}.
Recall also that the compact group $M$ is the maximal compact subgroup of the Lie group $T$ (cf.\ Proposition~\ref{Prop:Torus-Exp-Topologisch}) and that $K_l^\pm \cap B_\pm \subseteq M$. Thus $$
K \to K/(K \cap B_\pm)
$$
is a locally trivial fiber bundle by \cite{Palais-Bundles} with compact fibers, and so are
$$
\pi \colon K_l^\pm \to K_l^\pm/ (K_l^\pm \cap B_\pm).
$$
Thus $K_l^\pm$ is a compact set for every $l \in \NN$, therefore, the pair $(K, (K_l^\pm))$ is a $k_\omega$-pair. Together with the fact that
$$
\iota_\pm(K_l^\pm) = \widetilde{D}_{l, \pm}
$$
we can apply \Cref{Lem:k-Omega-1}.
\qed
\end{Proof}

At this point, we are able to prove the topological Iwasawa decomposition for Kac--Moody groups of arbitrary type over $\CC$. 

\begin{Theorem}\label{Thm:Iwasawa}
Let $G$ be a split Kac--Moody group over $\KK$ of type $\AA$, where $\AA$ is a symmetrizable, irreducible, non-spherical Cartan matrix. Then each of the multiplication maps
\begin{align*}
    m_\pm &\colon U_\pm \times A \times K \to G
\end{align*}
is a homeomorphism.
\end{Theorem}

\begin{Proof}
By the algebraic Iwasawa decomposition \Cref{Prop:AlgIwasawa} we know that $m_\pm$ is a bijection, so it follows from the continuity of multiplication in $G$ that $m_\pm$ is continuous with respect to the Kac--Peterson topology. Thus, all that remains is to construct an inverse continuous map as in the proof of \cite[Proposition 3.31]{FHHK_2017}.\\
Let $g \in G$ and define the two elements 
$$
k(g) \coloneqq \iota_\pm^{-1}\hspace{-2pt}\left( gAU_\pm \right) \; \text{ and } \; h(g) \coloneqq k(g)^{-1}g,
$$
where the maps $\iota_\pm$ are the same as in \Cref{Lemma:HilfeIwasawa}. Therefore we define the map
\begin{align*}
    n_\pm \colon G &\to K \times AU_\pm \\
    g = k(g)h(g) &\mapsto (k(g), h(g)),
\end{align*}
which are by construction inverse to $m_\pm$ and they are continuous, since $G$ is a topological group and $\iota_\pm$ is a homeomorphisms. It remains to check, that the multiplication maps $A \times U_\pm \to AU_\pm$ are open, but this follows from \cite[Corollary 7.30]{HKM_2013}, as the Levi decomposition is a semi-direct product of closed subgroups.
\qed
\end{Proof}

A useful application of the algebraic Iwasawa decomposition is the following lemma, which can be found for centered Kac--Moody groups over $\RR$ in \cite[Lemma 3.27]{FHHK_2017}.

\begin{Lemma}\label{Lem:NGT-tau(G)}
Let $G$ be a Kac--Moody group over $\KK$ of type $\AA$. Then $ N_G(T) \cap \tau(G) = A$.
\end{Lemma}

\begin{proof}
The following applies: $\tau(T) = A$ and $T \subseteq N_G(T)$, i.e.\ $A \subseteq N_G(T) \cap \tau(G)$. For the converse implication let $g \in N_G(T) \cap \tau(G)$. By \Cref{Prop:AlgIwasawa} there exist elements $u \in U_+$, $t \in T$ and $k \in K$, such that
$$
g = \tau(utk) = u \tau(t) \Theta(u)^{-1} \in U_+ \tau(t) U_-.
$$
Since the RGD system obtained from a Kac--Moody group is refined (see \cite[Theorem 1.5.4]{Remy_2002}), we can use the refined Birkhoff decomposition 
$$
G = \bigsqcup_{n \in N_G(T)} U_+ n U_-,
$$
see \cite[Theorem 5.2.3 (g)]{Kumar_2002}.\\
Since $g \in N_G(T)$ we conclude that $g = \tau(t)$ and hence $g \in A$. 
\end{proof}


\section{Kac--Moody Symmetric Spaces}

\subsection{A criterion for symmetric spaces}

We start with a technical building-theoretic lemma; for a brief introduction to buildings see \Cref{Sec:Buildings}. Let $\KK \in \{ \RR, \CC \}$ and $\AA$ be a symmetrizable, indecomposable generalized Cartan matrix.\\
An element $g \in G$ of a Kac--Moody group $G$ is called \emph{symmetric}, if
$$
\Theta(g) = g^{-1}.
$$
 Further we call a twin apartment $\Theta$-stable, if the twin building associated to $G$ is invariant as a set under the action of $\Theta$. Recall that $\Delta_\pm = G/B_\pm$, then we consider the action of $\Theta$ on the building as follows
$$
\Theta(g.B_\pm) = \Theta(g).B_\mp .
$$

\begin{Lemma}[Lemma 4.2 in \cite{Horn_2017}]\label{Lem:Equivalence}
Let $G$ be a Kac--Moody group over $\KK$ of type $\AA$ and denote by $\Delta_\pm$ the associated twin building. For a symmetric element $g \in G$ the statements below are equivalent:
\begin{myitems}
    \item \label{itm:fir} The element $g$ fixes a $\Theta$-stable twin apartment chamberwise. 
    \item \label{itm:sec} The element $g$ fixes a twin apartment chamberwise.
    \item \label{itm:thr} The element $g$ stabilizes a chamber.
    \item \label{itm:for} For all chambers $D$ in the twin building $\Delta_\pm$, the length of any minimal gallery from $D$ to a chamber in the set of $g$-orbits $\{ g^n.D \mid n \in \ZZ \}$, is bounded.
    \item \label{itm:fiv} For some chamber $D$ in the twin building $\Delta_\pm$, the length of any minimal gallery from $D$ to a chamber in the set of $g$-orbits $\{ g^n.D \mid n \in \ZZ \}$, is bounded.
    \item \label{itm:six} The element $g$ stabilizes a spherical residue in either half of the twin building.
\end{myitems}
\end{Lemma}

\begin{Remark}
    This statement was formulated and proven by Max Horn in \cite{Horn_2017}. We reproduce the proof here and, in fact, add some details for the reader's convenience.
\end{Remark}

In the following we use the notation of a twin apartment given in \cite[Definition 5.176]{Abramenko_2010}: Let $C$ and $C'$ be opposite chambers, i.e.\ $C \in \C_\epsilon$, where $\epsilon \in \{ +, - \}$ and $C' \in \C_{- \epsilon}$. Then $\Sigma(C, C')$ is defined as the set 
$$
\Sigma(C, C') \coloneqq \{ D \in \C_\epsilon \mid \delta_\epsilon(C,D) = \delta^\star(C',D) \}.
$$
For details on the notation, see \Cref{Sec:Buildings} or \cite[Chapter 5.8]{Abramenko_2010}.\\
Furthermore, a $J$-\emph{residue} is a certain equivalence class of chambers: two chambers $C, D$ in $\Delta_\epsilon$ are called $J$-equivalent if $\delta_{\epsilon}(C,D) \in W_J$, where $\W_J \coloneqq \langle J \rangle$ for a set $J \subseteq S$ w.r.t.\ Coxeter system $(\W, S)$ of $\Delta_\pm$. For details see \cite[Chapter 5.3]{Abramenko_2010}.

\begin{Proof} 
Let $g$ be a symmetric element of $G$. The implications $\ref{itm:fir} \Rightarrow \ref{itm:sec}$, $\ref{itm:sec} \Rightarrow \ref{itm:thr}$, and $\ref{itm:for} \Rightarrow \ref{itm:fiv}$ are clear.\\
The idea now is to prove that $\ref{itm:thr} \Rightarrow \ref{itm:fir}$, thus the first three points are equivalent. Then we show that 
$$
\ref{itm:thr} \Rightarrow \ref{itm:for} \Rightarrow \ref{itm:fiv} \Rightarrow \ref{itm:six} \Rightarrow \ref{itm:sec},
$$
which gives us the desired statement.
\begin{myitems}
    \item[$\ref{itm:thr} \Rightarrow \ref{itm:fir}$]
    Let $D$ be a chamber in $\Delta_\pm$ which is stabilized by $g$, i.e.\ $g.D = D$. If $g$ stabilizes $D$, then also $g^{-1}$ stabilizes $D$:
    $$
    D = (g^{-1} g).D = g^{-1}.(g.D) = g^{-1}.D.
    $$
    Hence
    $$
    \Theta(D) = \Theta(g.D) = g^{-1}. \Theta(D).
    $$
In particular, one concludes that $D = g.D$ implies $\Theta(D) = g.\Theta(D)$.
    Moreover, $g$ stabilizes the $\Theta$-stable twin apartment $\Sigma(D, \Theta(D))$, since $\Theta(D)$ is opposite of $D$. 
    \item[$\ref{itm:thr} \Rightarrow \ref{itm:for}$]
    First, recall that the building $\Delta_\pm$ associated with the Kac--Moody group $G$ is of type $(\W,S)$.\\
    Let $D' \in \Delta_+$ be a chamber in the positive half of the twin building which is stabilized by $g$. This implies that $\Theta(D') \in \Delta_-$ is stabilized by $g$ as well. Also let $D \in \Delta_+$ be an arbitrary chamber and $n \in \ZZ$. Then we know that
    $$
    \delta_+(D',D) = \delta_+(g^n.D', g^n.D) = \delta_+(D', g^n.D).
    $$
    Next, denote with $l \colon \W \to \NN$ the length function of the Weyl group $\W$ with respect to the generating set $S$. By \cite[Definition 5.1 (WD3)]{Abramenko_2010} one derives a triangle inequality, which is given for the length function $l$, so that one can compute
    $$
    l(\delta_+(D, g^n.D)) \leq l(\delta_+(D, D')) + l(\delta_+(D', g^n.D)) = 2l(\delta_+(D,D')).
    $$
    Hence the orbit $\{ g^n.D \mid n \in \ZZ \}$ is bounded and by the symmetry arguments of the to halves of a twin building, the claim follows.
    \item[$\ref{itm:fiv} \Rightarrow \ref{itm:six}$]
    This statement is exactly the Bruhat-Tits fixed point theorem applied on the CAT(0) realization of $\Delta_\pm$, see \cite[Corollary 12.67]{Abramenko_2010}.    
    \item[$\ref{itm:six} \Rightarrow \ref{itm:thr}$]
    As we have already seen, if $g$ stabilizes a chamber $D$, its inverse $g^{-1}$ also stabilizes $D$. The same calculation can be done to see, that if $g^{-1}$ stabilizes $D$, then also $g$ stabilizes $D$:
    $$
    D = (g g^{-1}).D = g.(g^{-1}.D) = g.D.
    $$
    Now let $R$ be a spherical residue, i.e.\ a residue where $\W_J$ is finite, stabilized by $g$ and by the computation
    $$
    \Theta(R) = \Theta(g.R) = g^{-1}. \Theta(R),
    $$
    and the previous comment, we have $\Theta(R)$ is also stabilized by $g$. Hence the symmetric element $g$ stabilizes the spherical residue in each half of $\Delta_\pm$.\\
    The idea is now, to show that $g$ fixes a chamber in $R$ and hence a twin apartment in $\Delta_\pm$. In particular, if $g$ fixes a chamber $D$ in $R$, $g$ fixes also $\Theta(D)$, which is an opposite chamber. This leads to the fact, that $g$ fixes then the twin apartment $\Sigma(D, \Theta(D))$ point-wise, i.e.\ $g$ stabilizes a twin apartment chamberwise.\\
    By \cite[Proposition 67.27]{Abramenko_2010}, one knows that the stabilizer of a spherical residue of type $J$ equals a spherical parabolic subgroup $P_J$, compare this with the discussion before \cite[Proposition 7.75]{Marquis_2018}. Since $g$ stabilizes $R$ and the opposite residue $\Theta(R)$, it is an element of 
    $$
    g \in P_J \cap \Theta(P_J).
    $$
    Using the fact, that a parabolic subgroup, in the context of a Kac--Moody group, has a Levi decomposition (see \cite[Theorem 6.2.2]{Remy_2002} or \cite[Proposition 7.75]{Marquis_2018}) one can analyze this intersection.
    Therefore, denote by
    $$
    P_J = G_J \ltimes U_J^+ \, \text{ and } \Theta(P_J) = G_J \ltimes U_J^-
    $$
    the Levi decomposition of $P_J$ and note that $\Theta(P_J)$ is obtained by applying the Chevalley involution. Note that in the convention of \cite[Chapter 6.2.2]{Remy_2002} $G_J$ is denoted by $M(F)$.\\
    Here we used the notion from \cite[Chapter 7]{Marquis_2018}, in particular the group $G_J$ stands for a subgroup of $G$, which is generated by the torus $T$ and the root groups $U_\a$, where $\a$ are the roots contained in the span of simple roots corresponding to $J$, see the discussion before \cite[Proposition 7.75]{Marquis_2018}. The second part of the Levi decomposition, $U_J^\pm$, is the normal subgroup of $U_\pm$ generated by $U_\a$, where $\a$ is in the complement of the roots used to generate $G_J$.\\    
    Since $g$ is an element of the intersection, we conclude that it can be written in two ways:
    $$
    x u_+ = g = y u_-,
    $$
    where $x,y \in G_J$, $u_+ \in U_J^+$ and $u_- \in U_J^-$. In fact, one realizes that
    $$
    g^{-1}.R = (u_-^{-1} y^{-1}).R = u_-^{-1}.R = R,
    $$
    since $y^{-1}$ is an element of the stabilizer (see the second statement in the theorem of \cite[Chapter 6.2.2]{Remy_2002}). From the proof of \cite[Lemma 6.2.1]{Remy_2002} (and for the sake of understanding with the notation from \cite[Theorem 6.2.2]{Remy_2002}) we obtain that $u_- \in P_J \cap U_J^- = \lbrace e \rbrace$, and analogously one can show with $\Theta(P_J)$ that $u_+$ must also be trivial.\\
    And so, the action of $g$ on the residue depends only on the part of $G_J$, and since $\W_J$ is finite, this is a Chevalley group, i.e.\ a linear algebraic group. By \cite[Proposition 1.10]{Borel_1991} we can consider this group as a closed subgroup of $\GL_n(\KK)$. To clarify the notation, we write $M_g$ for the corresponding matrix in the subgroup of $\GL_n(\KK)$.\\
    For a more detailed discussion of this embedding we refer to \cite[Chapter 3.2]{Platonov_1993} and in particular to \cite[Proposition 3.13]{Platonov_1993}. There it is explained that a spherical Kac--Moody group (i.e.\ in this concrete case $G_J$) can be embedded in $\GL_n(\KK)$. This is done by using an appropriate conjugation so that the Chevalley involution is translated into the transpose inverse composed with complex conjugation when $\KK = \CC$. Note that this embedding, where the Chevalley involution is translated in a pleasant way, must be chosen and is not given or natural.\\ 
    Taking a closer look on the property that $g$ is symmetric turns $M_g$ it into a Hermitian matrix, in particular
    $$
    \Theta(M_g) = M_g^{-1} \Longleftrightarrow  \left( \overline{M_g} \right)^{-T} = M_g^{-1} \Longleftrightarrow M_g^{T} = \overline{M_g}.
    $$
    Thus, $M_g$ is diagonalizable and hence it is conjugated to a diagonal matrix $D_g$. In particular, $g$ is conjugate to an element of a maximal torus $T_J$ of $G_J$. The torus $T_J$ embeds into a maximal torus of $G$, and so $g$ is conjugate to an element of the standard torus $T$ of the Kac--Moody group $G$. Since the standard torus fixes the pair of fundamental chambers $(C_+ , C_-)$ in $\Delta_\pm$, a conjugated element fixes a translated pair of the fundamental chambers. In detail if $T$ fixes $C_\pm$ and $aga^{-1} \in T$ denotes a conjugated element in $T$, one calculates
    $$
    T.C_\pm = C_\pm \Longrightarrow \left( aga^{-1} \right).C_\pm = C_\pm \Longrightarrow g.\left(a^{-1}.C_\pm \right) = a^{-1}.C_\pm.
    $$
    This completes the proof of this implication.
\end{myitems}
\qed
\end{Proof}

Next we use the preceding lemma to conclude that the homogenenous space $G/K$ for a Kac--Moody group $G$ is a symmetric space. This lemma is one of the key ingredients of the theory we develop in this article.\\
The version for non-affine centered Kac--Moody groups over $\RR$ is \cite[Proposition 3.38]{FHHK_2017}.

\begin{Lemma}\label{Lem:tau(G)-K-e}
Let $G$ be a split Kac--Moody group over $\KK$ of type $\AA$. Then
$$
\tau(G) \cap K = \{ e \}.
$$
\end{Lemma}

\begin{Proof}
Denote by $\Delta_\pm$ the associated twin building to the Kac--Moody group $G$. Further, let $g$ be an element of $\tau(G)$, i.e.\ $g = h \Theta(h)^{-1}$, for an element $h \in G$. Then one calculates
$$
\Theta(g) = \Theta(h \Theta(h)^{-1}) = \Theta(h) h^{-1} = g^{-1}.
$$
Thus, if $g \in K \cap \tau(G)$, one realizes that $g = \Theta(g) = g^{-1}$, from which we conclude that $g$ is symmetric and has order at most $2$. This leads to the fact that the orbits of $g$ are bounded as in point \ref{itm:for} of \Cref{Lem:Equivalence}. Hence by the proof of $\ref{itm:thr}$ to \ref{itm:fir} from \Cref{Lem:Equivalence} one knows, that $g$ fixes the $\Theta$-stable twin apartment $\Sigma(C, \Theta(C))$ chamberwise. Thus $g$ is contained in a corresponding torus $T'$ of $G$. From \Cref{Prop:AlgIwasawa} we know that the group $G$ can be decomposed as $G = KAU_\pm$. This can be turned into the decomposition 
$$
G = K B_\pm 
$$
by using $B_\pm = T \times U_\pm$ and the properties described in \Cref{Lem:Eigenschaften-A_M}. This decomposition yields that $K$ acts transitively on both halves of the associated twin building $\Delta_\pm$.\\
Hence there is a $k \in K$ such that $k T' k^{-1} = T$. Thus the element
$$
kgk^{-1} = k  \tau(h) \Theta(k^{-1}) = \tau(kh)
$$
is an element of $T$, $K$ and $\tau(G)$. By \Cref{Lem:Eigenschaften-A_M} and \Cref{Lem:NGT-tau(G)} we have
$$
kgk^{-1} \in T \cap K \cap \tau(G) = \{ e \}
$$
and thus $g = e$.
\qed
\end{Proof}

\subsection{Kac--Moody Symmetric Space}\label{Sec:KM-Sym-Space}
Loos' abstract definition of a symmetric space, see \Cref{Def:SymSp-Loos}, requires a set and a multiplication map $\mu$ that encodes point reflections within the symmetric space. It turns out that the Chevalley involution composed with the complex conjugation $\Theta$ (cf.\ \Cref{Sec:Chevalley}) is essential to define such a map $\mu$. Since the Kac--Moody group $G$ is a topological group, see \Cref{Sec:TitsFunctor}, we deduce that $K$ is a closed subgroup, since $\Theta$ is continuous and thus the quotient $G/K$ is a topological Hausdorff space. 

\begin{Proposition}\label{Prop:KM-SymSp}
Let $G$ be a split Kac--Moody group over $\KK$ of type $\AA$. Define the reflection map
\begin{align}
    \begin{split}
        \mu \colon G/K \times G/K &\to G/K \\
        (gK, hK) &\mapsto \tau(g)\Theta(h)K.
    \end{split}
\end{align}
The map $\mu$ is continuous, and if $\tau(G) \cap K = \{ e \}$, then the pair $(G/K, \mu)$ is a symmetric space and there is a natural action by automorphisms
$$
G \to \Sym(G/K) \, , \; g \mapsto (hK \mapsto ghK).
$$
\end{Proposition}

\begin{Proof}
The proof is a straightforward calculation, independent of the choice of $\AA$ and $\KK$, and can be found in \cite[Proposition 4.2]{FHHK_2017} or \cite[Proposition 4.3]{Diss}.
\qed
\end{Proof}

\begin{Definition}
The symmetric space $(G/K, \mu)$ is the \emph{Kac--Moody symmetric space} over $\KK$ of type $\AA$. 
\end{Definition}

\begin{Remark}
Note that this definition of the reflection map $\mu$ also holds for any topological group with an involution. Moreover, the first three properties in \Cref{Def:SymSp-Loos} can be proved without any assumption about the group, see the calculation in the proof of \cite[Proposition 4.2]{FHHK_2017}. Thus any topological group with an involution that satisfies the last property gives rise to a symmetric space. Here are two examples of such symmetric spaces.
\begin{myitems}
    \item The first example which comes to mind is the group given by $(\RR^n,+)$ with involution $\Theta \colon \RR^n \to \RR^n$, $\Theta(x) = -x$. Then one computes for $\tau(x) = x - \Theta(x) = 2x$ and hence the reflection of the symmetric space is given by $\mu(x,y) = \tau(x)+\Theta(y)+K = 2x-y$ and $K$ is equal to $\{ 0 \} \subseteq \RR^n$.

\newcommand{\SO}{\operatorname{SO}}

    \item Another example is the Lie group $\SL_n(\RR)$ with involution $\Theta(A) = A^{-T}$. Clearly the fixed point set $K$ equals the special orthogonal group $\SO(n)$. Note, that $\tau(A) = AA^T$ for a matrix $A \in \SL_n(\RR)$, hence the matrices in the image of the twist map are symmetric and thus diagonalizable over $\RR$. Moreover, the image $\tau(G)$ consists of positive semi-definite matrices.\\
    To check if $\tau(G) \cap K$ is trivial, we use that $\tau(A)$ can be diagonalized: let $S$ be an orthogonal real matrix such that
    $$
    S (AA^T) S^{-1} = D
    $$
    is a diagonal matrix. Assume now, that $AA^T \in K = \SO(n)$, then it is true that
    $$
    (AA^T)(AA^T)^T = (AA^T)^2 = \mathbb{E}_n,
    $$
    where $\mathbb{E}_n$ stands for the $n \times n$-identity matrix. Thus we can compute
    \begin{align*}
        \left( S (AA^T) S^{-1} \right)^2 &= D^2  \\
        S (AA^T)^2 S^{-1} &= D^2 \\
        S S^{-1} &= D^2
    \end{align*}
    and it can be deduced that the entries of $D$ are $\pm 1$. But by the positive semi-definite property only $+1$ is possible.
\end{myitems}
\end{Remark}

\subsection{The Group Model}\label{Sec:GroupModel}

So far, we can describe a symmetric space associated with a Kac-Moody group by $(G/K, \mu)$, which is called the \emph{coset model}. In this section we want to set up another model to describe the same symmetric space with $\tau(G)$ using the quadratic representation of an abstract symmetric space, see \Cref{Rem:QaudraticRepresentation}. The following section is based on \cite[Section 4.12]{FHHK_2017}.\\
In order to use the quadratic representation, we need a pointed symmetric space. Therefore choose $eK$ as the base point. Then, by \Cref{Rem:QaudraticRepresentation}, there is an isomorphism of symmetric spaces 
$$
(G/K, \mu) \cong (T(G/K, eK), \mu_T).
$$
Next, we take a closer look at $T(G/K, eK)$. Therefore, we first compute the kernel of the group action of $G$ on $G/K$, which is given by left multiplication:
$$
G \times G/K \to G/K \, , \; (g, hK) \mapsto ghK.
$$
If $g \in G$ is an element of the kernel, then $ghK = hK$ for all $h \in G$. In particular, this assertion also applies if $h \in K$ and in this case $ghK = K$, which implies $g \in K$.\\
Let $g$ now be an element of $K$, which acts trivially on $hK$ for an element $h \in G$. Then we can compute $ ghK = hK \Leftrightarrow h^{-1}gh \in K $ and by using \Cref{Lem:Eigenschaften-A_M} \ref{itm:A0} we have $\tau(h^{-1} g h) = e$. By a direct calculation this leads to $g \tau(h) g^{-1} = \tau(h)$, and thus the kernel of the action is contained in $C_K(\tau(G)) \subseteq K$.\\
For the converse consider an element $c \in C_K(\tau(G))$ and recall the following properties: The group $G$ can be decomposed into $G^\dagger$ and $T$, where $G^\dagger$ is generated by the root subgroups of $G$, see \Cref{Def:G1}. This allows us to apply \cite[Proposition 3.39]{FHHK_2017}, which gives us the following
$$
\left\langle \tau\hspace{-2pt}\left( G^\dagger \right) \right\rangle = G^\dagger.
$$
Now let $c \in C_K(\tau(G))$ and $h = tg \in G = TG^\dagger$ the we can calculate
\begin{align*}
    chK = ctgK = tcgK = tgcK = hcK = hK.
\end{align*}
We have therefore shown the reverse inclusion that $C_K(\tau(G))$ is contained in the kernel.


\smallskip

Next, denote by $\ell_g \colon G/K \to G/K$ the left multiplication with an element $g \in G$, i.e.\ $\ell_g(hK) = ghK$. It follows from the latter calculation that $\ell_g = \ell_h$ if $g$ and $h$ are in the same coset of $C_K(\tau(G))$. Then we can calculate 
\begin{align*}
    s_{gK}(hK) &= \tau(g) \Theta(h) K \\
    &= \ell_g \hspace{-2pt} \left( \Theta (g^{-1}h) K \right) \\
    &= \ell_g \hspace{-2pt} \left( s_{eK} \hspace{-2pt} \left( g^{-1}hK \right) \right) \\
    &= \left( \ell_g \circ s_{eK} \circ \ell_{g^{-1}} \right) \hspace{-2pt} (hK).
\end{align*}
Using this, we can rewrite a transvection:
\begin{align*}
    s_{gK} \circ s_{eK} &= \ell_g \circ s_{eK} \circ \ell_{g^{-1}} \circ s_{eK} = \ell_{\tau(g)}
\end{align*}
and thus
\begin{align*}
    T(G/K, eK) &= \{ s_{gK} \circ s_{eK} \mid gK \in G/K \} \\
    &= \{ \ell_{\tau(g)} \mid g \in G \} \\
    &= \{ \tau(g) C_K(\tau(G)) \mid g \in G \}.
\end{align*}
From this, one can define a map $\zeta \colon \tau(G) \to T(G/K, eK)$ by sending an element $\tau(g)$ of $\tau(G)$ to the coset $\tau(g) C_K(\tau(G))$. Then the previous observation shows that the map $\zeta$ is surjective. Furthermore, if an element is contained in the kernel of $\zeta$, it must be contained in $K$ by the previous calculation. Since $\tau(G) \cap K = \{ e \}$, see \Cref{Lem:tau(G)-K-e}, we conclude that the kernel of $\zeta$ is trivial. By transport of structure, the map establishes a symmetric space structure on $\tau(G)$, where the reflection map is given by $\mu_\tau(g,h) = gh^{-1}g$. 

\smallskip

Altogether there are the following isomorphism of symmetric spaces
$$
(G/K, \mu) \underset{\eta}{\cong} (T(G/K, eK), \mu_T) \underset{\zeta}{\cong} (\tau(G), \mu_\tau).
$$
The symmetric space $(\tau(G), \mu_\tau)$ is called \emph{group model} of the Kac--Moody symmetric space. The reflection map $\mu_\tau$ is the same as the reflection map given in \Cref{Ex:SymSpaces} \ref{itm:GroupModel} for a general group.

\begin{Remark}\label{Eq:Erzeugnis-Tau}
As in \cite[Proposition 4.11 (iii)]{FHHK_2017}, the transvection group of the Kac--Moody symmetric space can also be written as a quotient of the Kac--Moody group. But in the general case of arbitrary split Kac--Moody groups, we have to deal with the enlarged torus. We adopt the strategy used in the proof of \cite[Proposition 4.11]{FHHK_2017}, but we have to pay attention to the torus.\\
Define the following quotient, which, as we will see, acts faithfully on a Kac--Moody symmetric space:
$$
\eff{G^\dagger A} \coloneqq \nicefrac{G^\dagger A}{C_K(G^\dagger A)}.
$$
Recall that $G^\dagger$ is the derived subgroup of $G$, see the proof of \Cref{Lem:G1-is-charakteristisch}. Note that \emph{eff} in the index stands for \emph{effectively}, which means that the group acts faithfully on $G/K$.\\
Since $M \subseteq K$ and $\tau(T) = A$, see \Cref{Lem:Eigenschaften-A_M}, one concludes that 
$$
\langle \tau(G) \rangle = \langle \tau(G^\dagger T) \rangle \subseteq G^\dagger A.
$$
On the other hand, one knows that $G^\dagger = \langle \tau(G^\dagger) \rangle $ (cf.\ \cite[Proposition 3.39]{FHHK_2017}) as well as $\tau(T) = A$.
Note that \cite[Proposition 3.39]{FHHK_2017} is formulated for centered Kac--Moody groups with an RGD system, and $G^\dagger$ is a centered Kac--Moody group with an RGD system, see the discussion after \Cref{Rem:G1-is-KM}. Thus
$$
\langle \tau(G) \rangle = G^\dagger A.
$$
It follows from the previous observations that an element $g \in G$ acts trivially on a Kac--Moody symmetric space $G/K$ if and only if it centralizes $G^\dagger A$. In terms of the notation of left multiplication used above, $l(g_1) = l(g_2)$ only if $g_1$ and $g_2$ are in the same coset $C_K(G^\dagger A)$.

\smallskip

We now want to associate $\eff{G^\dagger A}$ with the transvection group of $G/K$, 
$$
\Trans(G/K) \coloneqq \langle s_{gK} \circ s_{hK} \mid gK, hK \in G/K \rangle.
$$
The strategy is the same as in the previous discussion and based on the proof of \cite[Proposition 4.11]{FHHK_2017}. Let $g \in G$ be an arbitrary element, $s_{gK}$ the corresponding point reflection and denote by $b = eK$ the base point of $G/K$. Then as done before, one calculates
\begin{align*}
\begin{split}
    s_{gK} \circ s_{b} = l(\tau(g)).
\end{split}
\end{align*}
Thus, for the generators of the transvection group, we observe the following presentation
\begin{align*}
    s_{gK} \circ s_{hK} &= \left( s_{gK} \circ s_{b} \right) \circ \left( s_{hK} \circ s_{b} \right)^{-1} \\
    &= l(\tau(g) \tau(h)^{-1}).
\end{align*}
Further, $\tau(h)^{-1} = \tau(\Theta(h))$, hence the transvection group is generated by elements contained in 
$$
\tau(G)^2 C_K \hspace{-2pt}\left(G^\dagger A\right).
$$
Together with \Cref{Eq:Erzeugnis-Tau} follows that
$$
\eff{G^\dagger A} = \Trans(\nicefrac{G^\dagger A}{K}) = \Trans(G/K).
$$
\end{Remark}

On the basis of this section, we have two isomorphic perspectives when we consider a symmetric Kac--Moody space. One is the coset model $G/K$, which is the natural way to construct the symmetric space according to Loos. The other is the group model $\tau(G)$, which allows us to view the symmetric space as a subset of the original Kac--Moody group $G$. This gives us an intrinsic, algebraic view of the symmetric space.

\subsection{Flats}
To analyze the structure of the Kac--Moody symmetric space of a Kac--Moody group $G$ over $\KK$ of type $\AA$, we first have to define flats. Euclidean flats are subsets of the symmetric space, which are isomorphic to the Euclidean symmetric space (cf.\ \Cref{Ex:SymSpaces} \ref{itm:EuclSymSpace}). In order to define Euclidean flats in the context of Kac--Moody symmetric spaces, we need to define a \emph{real form} of the Cartan subalgebra $\h = \sum_{1 \leq i \leq n} \CC \a^\vee_i + \sum_{n+1 \leq i \leq 2n-l} \CC h_i \cong \CC^{2n-l}$:
$$
\aa \coloneqq \erzeug{\RR}{\a^\vee_1, \ldots, \a^\vee_n, h_{n+1}, \ldots, h_{2n-l}} \cong \RR^{2n-l}.
$$

\begin{Proposition}\label{Prop:Flat-KM}
Let $G$ be a split Kac--Moody group over $\KK$ of type $\AA$ and equip the real form $\aa$ with its Euclidean symmetric space structure (cf.\ \Cref{Ex:SymSpaces} \ref{itm:EuclSymSpace}). 
Then the map 
\begin{align}
    \begin{split}
        \p_g \colon \aa &\to gAK \\
        X &\mapsto g \exp(X) K
    \end{split}
\end{align}
is an isomorphism of symmetric spaces for every $g \in G$. In particular, the subset $gAK$ of the Kac--Moody symmetric space $G/K$ is an Euclidean flat of dimension $\dim \aa = 2n-l$. 
\end{Proposition}

\begin{Proof}
First of all, it should be noted that by restricting the exponential mapping \Cref{Rem:ExpMap} to the real form $\aa$, $\KMExp(\aa)$ is bijective and the image is given by the subgroup $A$ of the torus. Furthermore, in this case $\KMExp$ corresponds to the natural exponential map and in particular to the group homomorphism between $(\RR^{2n-l},+)$ and $((\RR_{>0})^{2n-l}, \cdot)$.\\
Now, using that the action of $G$ on the symmetric space $G/K$ is given by an automorphism (cf.\ \Cref{Prop:KM-SymSp}), we can assume that $g = e$. Let $X, Y \in \aa$ and use that for all $t \in A$ we have that $\Theta(t) = t^{-1}$ (cf.\ \Cref{Sec:Chevalley}). Then one computes
\begin{align*}
\mu(\p_e(X), \p_e(Y)) &= \tau(\exp(X)) \Theta(\exp(Y)) K \\
&= \exp(X) \Theta(\exp(X))^{-1} \Theta(\exp(Y)) K \\
&= \exp(2X - Y) K = \p_e(\mu_{Eucl}(X,Y)).
\end{align*}
This shows, that $\p_g$ is an isomorphism of symmetric spaces.\\
To check whether the image $\mu(gaK, ga'K)$ is contained in $gAK$ for all $a, a' \in A$, we do a straightforward calculation
\begin{align*}
    \mu(gaK, ga'K) &= \tau(ga) \Theta(ga') K \\
    &= \tau(ga) \Theta(ga') K \\
    &= (ga) \Theta(ga)^{-1} \Theta(ga')K \\
    &= g a \Theta(a)^{-1} \Theta(g)^{-1} \Theta(g) \Theta(a') K \\
    &= g aa a'^{-1} K \in gAK.
\end{align*}

The last step is to check that $gAK \subset G/K$ is closed. \Cref{Thm:Iwasawa} tells that the multiplication $m \colon U_\pm \times A \times K \to G$ is a homeomorphism. This yields that the subspace $AK$ and $gAK$, for any $g \in G$, are closed in $G/K$ since the multiplication with an element induces a homeomorphism within the topological group.
\qed
\end{Proof}

\begin{Definition}
Call the map $\p_g$ from the previous proposition \emph{chart centered at $g$} and the image $\p_g(\aa) = gAK$ \emph{standard flat} for every $g \in G$ of the Kac--Moody symmetric space $(G/K, \mu)$.    
\end{Definition}

We conclude the section about flats by state a characterization. The following statement is based on \cite[Theorem 5.17]{FHHK_2017}. 

\begin{Theorem}\label{Thm:Characterization-Flats}
Let $(G/K, \mu)$ be the Kac--Moody symmetric space with respect to a Kac--Moody group over $\KK$ of type $\AA$, where $\AA$ is a symmetrizable, indecomposable generalized Cartan matrix. Then every weak flat is contained in a standard flat. Especially
\begin{itemize}
    \item standard flats are exactly the maximal flats;
    \item all weak flats are Euclidean;
    \item all weak flats are flats;
    \item $G$ acts transitively on maximal flats.
\end{itemize}
\end{Theorem}

\begin{Proof}
The proof of \cite[Theorem 5.17]{FHHK_2017} relies on the fact that a Kac--Moody group has an RGD system and that the coset model is isomorphic to the group model. The isomorphism between the group model and coset model is given in \Cref{Sec:GroupModel}. The arguments in the proof of \cite[Theorem 5.17]{FHHK_2017} based on RGD systems are independent of the type of the generalized Cartan matrix $\AA$ and of whether the underlying field $\KK$ of the Kac--Moody group is $\RR$ or $\CC$. Therefore, the proof given for \cite[Theorem 5.17]{FHHK_2017} is applicable, with the exception of claim $5$.\\
In claim $5$ a reductive split real Lie group is used and the fact that it is a subgroup of $\GL_m(\RR)$. For $\KK = \CC$ we also have an algebraic group which is a subgroup of $\GL_m(\CC)$ (cf.\ \cite[Proposition 1.10]{Borel_1991}. As we have seen in the proof of \Cref{Lem:Equivalence}, the involution $\Theta$ at the level of the matrix group can be understood as transpose, inversion, which is combined with complex conjugation (cf.\ \cite[Theorem 3.7 and Proposition 3.9]{Platonov_1993}). In the concrete calculation of the product $xy \cdot xy$, the upper left entry is equal to a sum of the squares of the absolute values of the entries in the first column. Now the arguments in claim $5$ work verbatim exactly as in the real case.
\qed
\end{Proof}

\subsubsection{Flats in the Group Model}
This section is based on \cite[Section 5.1]{FHHK_2017}.\\
In \Cref{Sec:GroupModel} we have defined two isomorphisms $\eta \colon G/K \to T(G/K)$ and $\zeta  \colon T(G/K) \to \tau(G)$. If we now combine them, we obtain the following isomorphism $\hat{\tau} \colon G/K \to \tau(G)$, $\hat{\tau}(gK) = \tau(g)$ between the coset and the group model and that it is equivariant under the action of $G$. We are now interested in the image of a standard flat in the group model under this isomorphism:
\begin{align*}
    \hat{\tau}(gAK) &= \tau(gA) \\
    &= g \tau(A) \Theta(g)^{-1}.
\end{align*}
In the following we use the notation
$$
g \star \tau(A) \coloneqq g \tau(A) \Theta(g)^{-1}.
$$
Thus, the image of the standard flat $gAK$ in the group model is given by $g \star \tau(A)$. Therefore, denote by $F[g] \coloneqq g \star A \subseteq \tau(G)$ the standard flat of the group model (cf.\ \cite[Proposition 5.4]{FHHK_2017}). In order to establish a relationship between the flats in the symmetric space $G/K$, the twin apartments in the corresponding building $\Delta_\pm$ and the maximal tori of $G$, we need the following lemma.

\begin{Lemma}[Remark 5.5 from \cite{FHHK_2017}]\label{Lem:N_KT=N_AT}
Let $G$ be a split Kac--Moody group over $\KK$ of arbitrary type. Then 
$$
N_K(A) = N_K(T).
$$    
\end{Lemma}
\begin{Proof}
Recall that in the case $\KK = \CC$ $T$ is connected and in the case $\KK = \RR$ $A$ is the identity component of $T$. This indicates that $N_K(T) < N_K(A)$, based on the fact that the conjugation is continuous.\\
For the reverse inclusion, consider an element normalizing $A$, then the conjugation with this element shifts $T$ to a maximal torus $T'$ of $G$. But the standard maximal torus $T$ containing $A$ is unique (a sufficiently large subgroup of $T$ fixes a unique twin apartment which in turn determines $T$; cf.\ \cite[Lemma 4.9]{Caprace_2009}). Therefore $T = T'$, thus this element normalizes $T$.
\qed
\end{Proof}

\begin{Proposition}[Remark 5.5 from \cite{FHHK_2017}]\label{Prop:Equi-Action-G-on-Flats}
There is a $G$-equivariant bijection between the set of maximal tori of $G$, the set of standard flats in $G/K$ and the twin apartments of $\Delta_\pm$, the corresponding twin building that belongs to $G$.\\
In fact, $G$ acts transitively on each of these objects.
\end{Proposition}

\begin{Proof}
Note that the strongly transitive action of $G$ on the twin building $\Delta_\pm$ implies the transitive $G$-action on the set of twin apartments; cf.\ \cite[Lemma~6.73]{Abramenko_2010}. The second claim is thus an immediate consequence of the first.

Denote by $ T_{max} $ the set of all maximal tori of $G$ and with $F_{std}$ the set of all standard flats in $G/K$. Then there is a $G$-equivariant map between these sets given by
\begin{align*}
    f \colon T_{max} &\to F_{std} \\
    gTg^{-1} &\mapsto gAK.
\end{align*}
To see that this map is well-defined, we first look at $N_G(T)$. We know by \cite[Lemma 3.27]{FHHK_2017} that $N_G(T) = A \ltimes N_K(T)$, since the proof of the quoted statement is given by \Cref{Lem:Eigenschaften-A_M} and the refined Birkhoff decomposition (see \cite[Theorem 5.2.3 (g)]{Kumar_2002} or proof of \Cref{Lem:NGT-tau(G)}). Together with \Cref{Lem:N_KT=N_AT} and \cite[Lemma 3.27]{FHHK_2017} follows that $N_G(T)$ stabilizes $AK$: Let $g \in N_G(T)$, then it can be written as $g = ak$ for $a \in A$ and $k \in N_K(A)$ and hence one calculates
$$
gAK = ak AK = a(kA^{-1}k)kK = aAK = AK.
$$
Together with the fact that $N_G(T)$ acts trivially by conjugation on $T_{max}$, we conclude that $G/N_G(T)$ can be identified with $T_{max}$ by mapping $gN_G(T) \mapsto gTg^{-1}$, in particular $f$ is surjective, see also the last point of \Cref{Thm:Characterization-Flats}.\\
For the injectivity, we use the group model and compute
\begin{align*}
    F[g]\Theta(F[g]) &= g \tau(A)\Theta(g)^{-1} \Theta(g) \Theta(\tau(A)) g^{-1} \\
    &= g \tau(A) \Theta(\tau(A)) g^{-1} \\
    &= g \tau(A) g^{-1} = gAg^{-1}.
\end{align*}
Since $A$ is contained in the standard maximal torus $T$ (see \cite[Lemma 4.9]{Caprace_2009}), $gAg^{-1}$ is contained in $gTg^{-1}$.\\
Thus $gTg^{-1}$ is uniquely determined by the flat $F[g]$, and the map $f$ is injective.\\
Note that $T$ is the stabilizer of the standard twin apartment and the group $G$ acts transitively on $T$ and the twin apartments. For the action of $G$ on $T$ the stabilizer of the standard torus is given by $N_G(T)$ and for the action of $G$ on the twin apartments is $N_G(T)$ the stabilizer of the standard twin apartment. Thus, for the twin apartments it can be deduced that the maximal tori of $G$ are exactly the chamberwise stabilizers of the twin apartments of $\Delta_\pm$.
\qed
\end{Proof}

\begin{Corollary}\label{Co:NG(T)=Stab(AK)}
Let $G$ be a split Kac--Moody group over $\KK$ of arbitrary type. Then 
$$
N_G(T) = \Stab_G(AK) \coloneqq \lbrace g \in G \mid g(AK) = AK \rbrace.
$$
\end{Corollary}

\begin{Proof}
The inclusion $ N_G(T) \subseteq \Stab_G(AK)$ follows from the proof of the theorem \Cref{Prop:Equi-Action-G-on-Flats}.\\
For the inverse inclusion, let an element $g$ of $\Stab_G(AK)$ be given, i.e.\ $ g(AK) = AK $. Then calculate
$$
gTg^{-1} = g(AM)g^{-1} = AM = T,
$$
where we have used the decomposition of $T$ and the fact that $M \subseteq K$, see \Cref{Lem:Eigenschaften-A_M}.
\qed
\end{Proof}


\section{Automorphism of Kac--Moody Symmetric Spaces}\label{Sec:Automorphism-KacMoody}
As a next step in developing the theory of symmetric spaces for Kac--Moody groups over $\KK$ of type $\AA$, one can analyze the automorphism group of the whole symmetric space and the local action on flats. For the non-affine real case, in chapter 6 of \cite{FHHK_2017}, the authors observe several statements for the global automorphism group on the symmetric space as well as for the local automorphism group with respect to the flats.

\subsection{Global Action}\label{Sec:GlobalAction}
From the characterization of flats one can deduce that the Kac--Moody group acts transitively on pointed maximal flats. In the following we call such an action of a Kac--Moody group $G$ on the corresponding symmetric space \emph{strongly transitive}, see \cite[Definition 2.32]{FHHK_2017}.

\begin{Proposition}\label{Prop:G-Strong-Transitive-Action-GK}
Let $G$ be Kac--Moody group over $\KK$ of arbitrary type. Then $G$ acts strongly transitively on the symmetric space $G/K$.
\end{Proposition}
\begin{Proof}
We are going to use \cite[Proposition 2.33]{FHHK_2017}. This statement is shown for abstract symmetric spaces, hence we can apply it here. In order to use this statement, we will check the required conditions. In detail, it has to be checked that $\Trans(G/K) < G < \Aut(G/K)$, furthermore that $G$ acts transitively on maximal flats and that all maximal flats are Euclidean. The right inclusion of the first condition is obviously satisfied, and the left inclusion is shown in \Cref{Eq:Erzeugnis-Tau}.\\
From \Cref{Prop:Equi-Action-G-on-Flats} it follows that $G$ acts transitively on standard flats. Since all standard flats are maximal flats and all maximal flats according to \Cref{Thm:Characterization-Flats} are Euclidean in our context, all assumptions of \cite[Proposition 2.33]{FHHK_2017} are satisfied.
\qed
\end{Proof}

A central assumption in \cite{FHHK_2017} which allows to decompose an automorphism, is that the Kac--Moody group is centered (cf.\ \Cref{Rem:G1-is-KM}), i.e.\ the Kac--Moody group is generated only by its root groups. To establish a similar results for non-centered Kac--Moody groups $G$ of type $\AA$, we have to take care of the extended torus. In order to do this, we consider in the following section $G^\dagger \trianglelefteq G$.\\
By \Cref{Rem:G1-is-KM} we know, that the subgroup $G^\dagger$ is a centered Kac--Moody group, thus using \cite[Theorem 2.4]{CM_2005}, we can write any automorphism of $G^\dagger$ as a product of an inner automorphism, the Chevalley involution, a diagonal automorphism, a graph automorphism and a field automorphism. Following, we give a brief overview of these automorphism, for details see \cite[Section 8.2]{CM_2005}.\\
To do this, we need to recall some facts from \Cref{Sec:TitsFunctor}. Firstly, the map $\eta$ from the Tits functor; secondly, the fact that an element of the root group $U_\a$ is given by $x_\a(r)$ for $r \in \KK \in \{ \RR, \CC \}$ and $\a \in \Delta^{re}$; and, lastly, that the torus of a Kac--Moody group $T = \T_\Lambda(\KK)$ is isomorphic to $\Hom_{grp}(\Lambda, \KK^\times)$. Hence we can define the map
$$
\iota(t) \colon \Lambda \to \KK^\times \, , \; \lambda \mapsto t(-\lambda) = (t(\lambda))^{-1}
$$
for each $t \in T$. Now, for a given centered split Kac--Moody group $H$ of type $\AA$, $r \in \KK$, $\a_i \in \Pi$, $\a \in \Delta^{re}$ and $t \in T$, one can describe the types of the automorphisms as follows:
\begin{myitems}[itemsep=0pt]
    \item An \emph{inner automorphism} $c_g \colon H \to H$ for an element $g \in H$ is given by $c_g(h) = g h g^{-1}$. 

    \item A \emph{sign automorphism} is a map $s \colon H \to H$ such that 
    $$
    s(x_{\a_i}(r)) = x_{-\a_i}(r) \, , \; s(x_{-\a_i}(r)) = x_{\a_i}(r) \, \text{ and } \; s(\eta(t)) = \eta(\iota(t)).
    $$

    \item A \emph{diagonal automorphism} is the identity on the torus and acts on $x_{\a}(r)$ by multiplying $r$ with an element of the field depending on $\a \in \Delta^{re}$.

    \item Let $\sigma$ be a permutation of the index set $I$, then a \emph{diagram automorphism} is given as a map $\gamma$ such that
    $$
    \gamma(x_{\a_i}(r)) = x_{\a_{\sigma(i)}}(r) \, , \; \gamma(x_{-\a_i}(r)) = x_{-\a_{\sigma(i)}}(r) \, \text{ and } \; \gamma(\eta(r^{h_i})) = \eta(r^{h_{\sigma(i)}}).
    $$    

    \item A \emph{field automorphism} $f \colon \KK \to \KK$ acts on $H$ using the fact that $H$ is defined by a functor. In particular for any $t \in T$ it is $f \circ t$ and for any $\a \in \Delta^{re}$ and $r \in \KK$ it is $f(x_{\a}(r)) = x_{\a}(f(r))$.
\end{myitems}

Note that we choose $H$ as $G^\dagger$ for our purpose.

\begin{Remark}\label{Rem:Theta=Sign}
From the discussion in \cite[Section 8.2]{CM_2005} follows that the sign automorphism is unique and from \cite[Section 8.2.1]{Caprace_2009} follows that it is related to the Chevalley involution. This leads to the conclusion that for $\KK = \RR$ one has
$$
s = \Theta
$$
and for $\KK = \CC$ one has
$$
s \circ f = \Theta,
$$
where $f$ denotes the field automorphism of the complex conjugation.
\end{Remark}

With this concrete decomposition, we can now rewrite \cite[Proposition 6.5]{FHHK_2017} in the context of the centered subgroup $G^\dagger$ of a Kac--Moody group $G$ of type $\AA$ over $\RR$ or $\CC$.

\begin{Proposition}\label{Prop:AutG1-Zerlegung}
Let $G$ be a split Kac--Moody group of type $\AA$ over $\KK \in \{ \RR, \CC \}$ and $G^\dagger \trianglelefteq G$ the centered subgroup. Then one can decompose the automorphism group, i.e.\
$$ \Aut(G^\dagger) \cong \left( \Ad(G^\dagger) \rtimes \left( D \times \langle \Theta \rangle \right) \right) \rtimes \Aut(\Gamma),$$
where $\Ad(G^\dagger) \coloneqq G^\dagger / Z(G^\dagger)$, $D$ is a group of diagonal automorphisms and $\Aut(\Gamma)$ denotes the diagram automorphisms of the corresponding Dynkin diagram with respect to $\AA$.
\end{Proposition}
Note that this explicit statement about the structure of the automorphism group is mainly based on the result of \cite[Theorem 2.4]{Caprace_2009}.
\begin{Proof}
The case $\KK = \RR$ is proven in \cite[Proposition 6.5]{FHHK_2017}, where it is formulated for any generalized Cartan-Matrix $\AA$.\\
For the complex case, we recall the decomposition of an automorphism of $G^\dagger$ (\cite[Theorem 2.4]{CM_2005}): Denote with $c_g$ an inner automorphism w.r.t.\ an element $g \in G$, with $\omega$ the Chevalley involution, with $d$ a diagonal automorphism, with $\gamma$ a diagram automorphism and with $f$ a field automorphism. Then any automorphism $\a \in G^\dagger$ can be written as 
$$ \a = c_g \circ \omega \circ d \circ \gamma \circ f.$$
In the complex case, there is only one non-trivial continuous field automorphism $f$, the complex conjugation. By the definition of $\gamma$ and $d$ we know that $f$ commutes with them and we can deduce that
$$
\a = c_g \circ \omega \circ f \circ d \circ \gamma = c_g \circ \Theta \circ d \circ \gamma,
$$
where we have used that $\Theta$ is defined as $\omega \circ f = \Theta$.\\
Since the varioua subgroups generated by these automorphisms intersect trivially, it follows that the product of the generated subgroups is semi-direct.
\qed
\end{Proof}

\begin{Remark}\label{Rem:Embedding-Aut}
\begin{itemize}
    \item[] 
    \item Since we are interested in statements about Kac--Moody symmetric spaces, we also need to define a symmetric space with respect to the centered subgroup $G^\dagger$. In particular, let $K^{\dagger} \coloneqq G^\dagger \cap K$, and it follows directly $(G^\dagger )^\Theta = K^{\dagger}$. Now define the corresponding Kac--Moody symmetric space 
    $$
    \left( \nicefrac{G^\dagger}{K^{\dagger}} , \res{G^\dagger}{\mu} \right)
    $$
    and denote it by $X^\dagger$. Note that this is the unreduced Kac--Moody symmetric space, which is also treated in part in the work on which this is based, see \cite[Definition 4.5 (i)]{FHHK_2017}.

    \item Since $G^\dagger$ is a centered Kac--Moody group which corresponds to the Kac--Moody group treated in \cite{FHHK_2017}, we can apply \cite[Proposition 6.4]{FHHK_2017}. The proposition states that 
    $$
    \Aut(G^\dagger) \hookrightarrow \Aut(\Delta)
    $$
    which allows us to interpret automorphisms of $G^\dagger$ as automorphisms of the corresponding twin building $\Delta_\pm$ from the RGD system associated to the Kac--Moody group $G$. By an automorphism of $\Delta$ we mean here a map of the chamber set $\Delta = \Delta_+ \cup \Delta_-$ into itself, preserving the adjacency relation and the opposition relation.
\end{itemize}
\end{Remark}

We conclude this section by pointing out the following statement, which describes the automorphism group of a Kac--Moody symmetric space coming from the centered subgroup. 

\begin{Theorem}[Theorem 6.12 (i) \cite{FHHK_2017}]\label{Thm:X1-Isom}
Let $G$ be a split Kac--Moody group of type $\AA$ over $\KK \in \{ \RR, \CC \}$ and $G^\dagger \trianglelefteq G$ its centered subgroup. Then 
$$
\Aut \left( X^\dagger \right) \cong \AutEff(G^\dagger) \coloneqq \left( \nicefrac{G^\dagger}{C_{K^{\dagger}}(G^\dagger)} \rtimes (D \times \langle \Theta \rangle) \right) \rtimes \Aut(\Gamma) ,
$$
where $C_{K^{\dagger}}(G^\dagger)$ is the centralizer of $G^\dagger$ in $K^{\dagger}$.
\end{Theorem}

\begin{Proof}
The proof is analogous to the proof for \cite[Theorem 6.12 (i)]{FHHK_2017}. Note that, by \Cref{Rem:Embedding-Aut}, the symmetric space $X^\dagger$ is the same symmetric space considered in \cite{FHHK_2017}. Additionally, the underlying field and type of Kac--Moody groups are irrelevant as they utilize \Cref{Prop:AutG1-Zerlegung}, a basic statement applicable to all cases.
\qed
\end{Proof}

\section{Local Action}\label{Sec:LocalAction}
One can define local transformations of flats for general symmetric spaces (cf.\ \cite[Section 2.31]{FHHK_2017}) and then one can transfer these definitions into the context of Kac--Moody symmetric spaces and analyze them.\\
In this section we denote by $G$ a split Kac--Moody group over $\KK \in \{ \RR, \CC \}$ of type $\AA$, where $\AA$ denotes an indecomposable, symmetrizable generalized Cartan matrix of arbitrary type. Further denote by $\BB$ the bilinear form on the Cartan subalgebra $\h$ induced by the symmetrizable generalized Cartan matrix. For details of the construction of this bilinear form, see \cite[§2.1]{Kac_1990} (or \cite[Section 3.7]{Marquis_2018} or \cite[A.22]{FHHK_2017}).\\
This section is primarily based on \cite[Section 6.14]{FHHK_2017}. We will sometimes refer to statements from \cite{FHHK_2017} without any justification when the statements we are referring to do not depend on the underlying field or the type of the Kac--Moody group. We will point out statements that require further explanation before being applicable into our context.

\bigskip

For an abstract symmetric space according to Loos one can define the following terms, compare them with \cite[Definition 2.32]{FHHK_2017} and \cite[Definition 2.35]{FHHK_2017}: Let $X$ be an abstract symmetric space containing a maximal Euclidean flat $F$ of rank $m$ and denote by $p$ a point of $F$. 
\begin{itemize}
    \item The pair $(p,F)$ is called a \emph{pointed maximal flat}. 
    \index{Flat! Pointed Maximal}
    
    \item Call a point $q \in F$ \emph{singular} with respect to $p$, if there exists a second maximal flat $E \subseteq X$ such that $p,q \in E$ and $E \neq F$. If not, call $q$ \emph{regular}. Further, let $\sing{F}(p)$ denote the subset of singular points of $(p,F)$ and $\reg{F}(p)$ the set of regular points. Hence the flat $F$ can be written as $F = \reg{F}(p) \sqcup \sing{F}(p)$.
    \index{Point!Singular} \index{Point!Regular}

    \item The isomorphism of symmetric spaces $\p \colon \RR^{m} \to F$, see the discussion before \cite[Proposition 2.21]{FHHK_2017}, is called a \emph{chart of $F$ centered at $p$} if $p = \p(0)$.
    
    \item Call a map $f \colon F \to F$ \emph{linear at $p$} if for a chart of $F$ centered at $p$ the following is true 
    $$
    \p^{-1} \circ f \circ \p \in \GL_{m}(\RR).
    $$
    \index{Linear Map}

    \item A map $f \colon F \to F$ linear at $p$ is called \emph{local transformation} of a pointed flat $(p,F)$ if it preserves the decomposition 
    $$
    F = \reg{F}(p) \sqcup \sing{F}(p).
    $$
    The set of all local transformations of $(p,F)$ forms a group and is denoted by $\GL(p,F,\sing{F})$.
    \index{Local Transformation $\GL(p,F,\sing{F})$}
\end{itemize}

\bigskip

The idea in the following is to use these terms in the context of Kac--Moody symmetric spaces, since we have a good understanding of flats there. In \cite[Section 6.14]{FHHK_2017} the authors give a description of singular points in terms of the real form $\aa$ and connect it to the notion of singular points of abstract symmetric spaces. Afterwards they analyze the group of local transformations and compare the Weyl group of the Kac--Moody symmetric space with the Weyl group of the Kac--Moody algebra.\\
The difficulty in applying this description to the general case for all generalized Cartan matrices $\AA$ is that in \cite{FHHK_2017} the semisimple adjoint quotient of a Kac--Moody group is used (cf. \cite[Definition 3.8]{FHHK_2017}) and for Kac--Moody groups of affine type the quotient prevents a representation of the Weyl group from being faithful and the corresponding bilinear form on $\h$ from being non-degenerate (cf. \cite[Section A.22]{FHHK_2017}). In the following we give a general definition of the representation which is faithful for real or complex split Kac--Moody groups of arbitrary type.\\
Recall that the Weyl group $\W$ of the Kac--Moody algebra $\g(\AA)$ is a subgroup of $\GL(\h)$. As the reflections depend on the generalized Cartan matrix $\AA$, which has integer entries, the action stabilizes the real form $\aa$, so we can define the \emph{Kac--Moody representation}\index{Kac--Moody Representation $\rho_{KM}$}
\begin{align}\label{Eq:Rho-KM}
\begin{split}
    \rho_{KM} \colon \W \to \GL(\aa).
\end{split}
\end{align}

According to \cite[Proposition 3.12]{Kac_1990}, this representation is faithful and the invariant symmetric bilinear form $\BB$ on $\h$ is non-degenerate, see \cite[Lemma 2.1 b)]{Kac_1990} (or \cite[Theorem 2.2]{Kac_1990} or \cite[Proposition A.23]{FHHK_2017}). Note that the generating reflections $r_i^\vee \in \W$ can also be defined using the bilinear form $\BB$, see \cite[A.22]{FHHK_2017}.\\
Since the elements of $\W$ act as linear reflections on $\aa$ under the representation $\rho_{KM}$, we can define a \emph{reflection hyperplane} with respect to a real root $\a \in \Delta^{re}$ as
$$
H_\a \coloneqq \Fix(\rho_{KM}(r_\a^\vee)) \subseteq \aa.
$$
From the non-degeneracy of $\BB$ on $\h$ we conclude that $\BB$ is also non-degenerate on $\aa$. Therefore, the reflection $r_\a^\vee$ under $\rho_{KM}$ is the only $\BB$-orthogonal reflection on the hyperplane $H_\a$. 
In particular, this leads to a one-to-one correspondence between the real roots $\a \in \Delta^{re}$ and the hyperplanes $H_\a$.\\
Using these concepts, we can define singular points in terms of the real form $\aa$ as the union of the reflection hyperplane:
$$ 
\sing{\aa} \coloneqq \bigcup_{\a \in \Delta^{re}} H_\a.
$$

\smallskip

As we saw at the beginning of this section, we can also define singular points on $\aa$ by using the notion of singular points in the context of abstract symmetric spaces. To be able to compare the two terms, we need to make them compatible. Therefore, let $(e,AK)$ be the standard flat of the symmetric space $G/K$. Using the chart $\p_e$ (for the definition see \Cref{Prop:Flat-KM}), one can move the singular points of the flat $AK$ into the real form $\aa$, i.e.
$$
\p_e^{-1}\left( \sing{(AK)}(e) \right),
$$
where $\sing{(AK)}(e)$ denotes the set of singular points of $(e,AK)$.\\
Now we can state that the descriptions of singular points on the real form are identical. The statement for non-affine real Kac--Moody symmetric spaces can be found in \cite[Proposition 6.16]{FHHK_2017}.

\begin{Proposition}\label{Prop:Sing-Points-Equal}
Let $G/K$ be a Kac--Moody symmetric space. Then
$$
\bigcup_{\a \in \Delta^{re}} H_\a = \sing{\aa} = \p_e^{-1}\left( \sing{(AK)}(e) \right).
$$
\end{Proposition}
\begin{Proof}
We can use the proof of \cite[Proposition 6.16]{FHHK_2017} verbatim, since it starts with an assertion that follows in our general case from \Cref{Prop:G-Strong-Transitive-Action-GK}, and then uses only building-theoretic arguments that are available since each Kac--Moody group has an RGD system over $\KK \in \{ \RR, \CC\}$, see \Cref{Rem:RGD}. Note that the type of $\AA$ does not matter.
\qed
\end{Proof}

We are using this statement to analyze the group of local transformations. Therefore, define the group
$$
\GL(\aa,\sing{\aa}) \coloneqq \{ f \in \GL(\aa) \mid f(\sing{\aa}) = \sing{\aa} \},
$$
which preserves the decomposition of $\aa$ into singular and regular points. By \Cref{Prop:Sing-Points-Equal} we can deduce that there is an isomorphism 
\begin{align*}
    \GL(e, AK, \sing{(AK)}) &\to \GL(\aa, \sing{\aa}) \\
    f &\mapsto \p_e^{-1} \circ f \circ \p_e.
\end{align*}
Since $G$ acts strongly transitively, see \Cref{Prop:G-Strong-Transitive-Action-GK}, any pointed maximal flat $(p,F)$ can be mapped to the standard pointed flat $(e, AK)$. Furthermore, the $G$ action preserves the decomposition into regular and singular points. Thus we observe an isomorphism
\begin{align*}
    \GL(p, F, \sing{F}) &\to \GL(\aa, \sing{\aa}) \\
    f &\mapsto \p_p^{-1} \circ f \circ \p_p,
\end{align*}
where $\p_p$ stands for a chart centered at $p$ and by \Cref{Prop:Sing-Points-Equal} the chart identifies $\sing{\aa}$ with $\sing{F}(p)$.

\begin{Convention}
In order to be able to use some corollaries from the appendix of \cite{FHHK_2017}, we need an additional assumption in the rest of this section: the size $n$ for the symmetrizable generalized Cartan matrix $\AA \in \ZZ^{n \times n}$ must be at least $2$.
\end{Convention}

By Corollary A.37 from \cite{FHHK_2017} follows that we can write any element of $\GL(\aa, \sing{\aa})$ as a product of a $\BB$-orthogonal linear transformation and a homothety.\\
Note that \cite[Corollary A.37]{FHHK_2017} is applicable here, since the arguments and statements that the authors used in the proof of this corollary are based on the fact that the bilinear form $\BB$ coming from the generalized Cartan matrix $\AA$ is non-degenerate and that $\rho_{KM}$ is faithful. All these properties are also satisfied in our context and thus we can apply the corollary.\\
Hence, by the isomorphism between $\GL(p, F, \sing{F})$ and $\GL(\aa, \sing{\aa})$, we can use this decomposition of $\GL(\aa, \sing{\aa})$ to rewrite the group of local transformations.\\
Therefore given a pointed flat $(p,F)$ and a local transformation $f \in \GL(p,F, \sing{F})$, then $f$ is a \emph{local automorphism}\index{Local Automorphism $\Aut(p,F)$} if $ \p^{-1} \circ f \circ \p $ is $\BB$-orthogonal for a chart $\p$ centered at $p$. We denote the subgroup of local automorphism by $\Aut(p,F) < \GL(p,F,\sing{F})$. By the previous discussion, we can decompose any local transformation into a product of a homothety and a $\BB$-orthogonal transformation.
Thus we can decompose the group of local transformation as follows
\begin{align}\label{Eq:GL=Aut(pF)}
\GL(p,F,\sing{F}) \cong \RR_{>0} \times \Aut(p,F),
\end{align}
where $\RR_{>0}$ represents the action of the homotheties on the flat $F$. 

\begin{Proposition}[cf.\ Corollary 6.20 in \cite{FHHK_2017}]\label{Prop:GL-Isom}
Let $G$ be a Kac--Moody group over $\KK$ of type $\AA$, where $\AA \in \ZZ^{n \times n}$ is non-spherical and symmetrizable. For any pointed flat $(p,F)$ one has
$$ 
\Aut(p,F) \cong (\W \rtimes \Aut(\W,S)) \times \ZZ/2\ZZ
$$
and hence
$$
\GL(p,F,\sing{F}) \cong \RR_{>0} \times (\W \rtimes \Aut(\W,S)) \times \ZZ/2\ZZ,
$$
where $\Aut(\W,S)$ denotes the group of automorphisms of $\W$ that preserve $S$ as a set. 
\end{Proposition}
\begin{Proof}
Recall that the Weyl group $\W$ of the Kac--Moody group with generating set $S$ is a Coxeter system $(\W, S)$ and that a Coxeter system induces a simplicial Coxeter complex $\Sigma$ (cf.\ \cite[Chapter 3]{Abramenko_2010} and \cite[Section A.1]{FHHK_2017}). Denote in the following with $\Aut(\Sigma)$ the group of simplicial automorphisms. By \cite[Lemma A.31]{FHHK_2017} one knows that
$$
\Aut(\Sigma) = \W \rtimes \Aut(\W,S)
$$
and there the action of $\Aut(\Sigma)$ on $\aa$ is also described.\\
Next denote by $\operatorname{O}(\aa, \BB)$ the group of orthogonal linear maps with respect to the bilinear form $\BB$ and define the group
$$
\operatorname{O}(\aa,\sing{\aa}) \coloneqq \operatorname{O}(\aa,\BB) \cap \GL(\aa, \sing{\aa}).
$$
By \cite[Proposition A.32]{FHHK_2017} and \cite[Remark A.33 (ii)]{FHHK_2017}, one derives that 
$$
\Aut(\Sigma) \times \ZZ/2\ZZ \cong \operatorname{O}(\aa, \sing{\aa})
$$
and hence that
$$
\Aut(\Sigma) \times \ZZ/2\ZZ \cong \Aut(p,F).
$$
Note that the quoted statements apply to Kac--Moody groups of arbitrary type, since they require a faithful representation $\rho_{KM}$ and a non-degenerate bilinear form $\BB$ on $\aa$. Both are given in the general context.\\
Together with \cite[Lemma A.31]{FHHK_2017}, we obtain the first isomorphism. The second isomorphism follows directly from \Cref{Eq:GL=Aut(pF)}. 
\qed
\end{Proof}

\subsection{Comparing the Weyl Groups}
As at the end of \cite[Section 2] {FHHK_2017}, one can define a Weyl group for an abstract symmetric space. In this section we will relate this Weyl group of the symmetric space to the Weyl group of the Kac--Moody group. The section follows \cite[Section 2.31]{FHHK_2017}.

\smallskip

Let $(X, \mu)$ be a symmetric space containing a maximal Euclidean flat, let $G < \Aut(X)$ be a group acting transitively on maximal flats of $X$, and let $(p,U)$ be a maximal pointed flat. Then define the following terms:
$$
\Stab_G(p,U) \coloneqq \lbrace g \in G \mid g(U) = U \, , \; g(p) = p \rbrace
$$
and
$$
\Fix_G(p,U) \coloneqq \lbrace g \in G \mid \forall u \in U : g(u) = u \rbrace.
$$
Note that $\Fix_G(p,U)$ is the same as $\Fix_G(U)$, but adding the $p$ in the notation makes it clear that we are considering an pointed flat here.

Then one can define the \emph{geometric Weyl group}\index{Weyl! Geometric Group} by 
$$
W(G \curvearrowright X) \coloneqq \nicefrac{\Stab_G(p,U)}{\Fix_G(p,U)},
$$
see \cite[Definition 2.37 (i)]{FHHK_2017}.\\
Note that by \cite[Proposition 2.36 (i)]{FHHK_2017} the definition of the geometric Weyl group is independent of the choice of the pointed flat up to conjugation.

\bigskip

Now we are able to compare the Weyl group of a Kac--Moody symmetric space with the Weyl group of the underlying Kac--Moody group. Therefore, denote in the following with $G$ a real Kac--Moody group of arbitrary type $\AA$. The reason why we consider a real Kac--Moody group goes back to the structure of the subgroup $M$ of the torus $T$ (see \Cref{Lem:Eigenschaften-A_M}). In the following discussion we will see that it is necessary that the elements of $M$ have order $2$, which is true for real Kac--Moody groups.\\
Recall from the discussion after \Cref{Def:G1} that the centered subgroup $G^\dagger \trianglelefteq G$ leads to a twin building of the same type as $G$ and hence that the Weyl groups are equal, i.e.
\begin{align}\label{Eq:G_dagger-Cong-W}
    \nicefrac{\N{G^\dagger}{T^\dagger}}{T^\dagger} \cong \W \cong \nicefrac{N_G(T)}{T}.
\end{align}

Recall that the extended Weyl group $\EW \coloneqq \langle \Tilde{s}_{i} \mid i \in I \rangle$ (cf.\ \Cref{Sec:TitsFunctor}) surjects on $\W$.\\
Due to \cite[(2.4)]{Kac-Peterson_1985} (see also \cite[Exercise 7.49]{Marquis_2018}), the elements $\Tilde{s}_i^2 = (-1)^{h_i}$ are contained in the torus and thus they generate a subgroup of $\EW \cap T$ which has order $2^n$. In addition, according to \cite[Corollary 2.3 (a)]{Kac-Peterson_1985} the elements of the intersection $\EW \cap T$ have an order less than or equal to $2$.\\
Note that the standard torus of $G^\dagger$ is denoted by $T^\dagger$ and can be decomposed as a topological group into 
$$
T^\dagger = M^\dagger \times A^\dagger, 
$$
where $A^\dagger$ and $M^\dagger$ are defined as in \Cref{Lem:Eigenschaften-A_M}. By the construction of $G^\dagger$ follows that $T^\dagger \cong (\RR^\times)^n$, where $n$ is the size of the Cartan matrix $\AA$. Hence it follows by the construction of $M^\dagger$ (cf.\ \Cref{Lem:Eigenschaften-A_M}) that $M^\dagger$ is the unique maximal finite subgroup of elements of order $2$. This implies the following  
$$
\EW \cap T = M^\dagger.
$$
Further, by the definition of $\EW$ it follows that the extended Weyl group is a subgroup of the centered subgroup $G^\dagger$. In \Cref{Sec:TitsFunctor} we have seen that $\EW$ surjects through $\Ad(-)$ restricted to $\h$ onto $\W$. Moreover, the kernel of this surjective map is exactly the intersection of $T$ with $\EW$. Thus we obtain the isomorphism
$$ 
\nicefrac{\EW}{(\EW \cap T)} \cong \W.
$$
Alternatively, we can also use \cite[Proposition 2.1]{Kac-Peterson_1985} to obtain the same statement.\\
Altogether we conclude that 
\begin{equation}\label{Eq:EW-M-Isom}
    \nicefrac{\EW}{M^\dagger} \cong \W. 
\end{equation}

To compare the Weyl groups we need the following helpful lemma, which is a generalization of \cite[Corollary 3.28]{FHHK_2017}.
\begin{Lemma}[cf.\ Corollary 3.28 from \cite{FHHK_2017}]\label{Lem:EW=NK(T)}
Let $G$ be a split Kac--Moody group of arbitrary type $\AA$. Then the extended Weyl group satisfies
$$
\EW M = N_K(T).
$$
\end{Lemma}
\begin{Proof}
Consider the canonical projection
$$
\pi \colon N_G(T) \to \nicefrac{N_G(T)}{T},
$$
and recall from \Cref{Lem:Eigenschaften-A_M} that $M \subset K$. In the proof of \cite[Lemma 3.27 (ii)]{FHHK_2017}, using the Iwasawa decomposition and the refined Birkhoff decomposition, it is shown that $N_G(T) = N_K(T)T$ holds. Note that the refined Birkhoff decomposition can be used since the $BN$-pair obtained from $G$ is refined according to \cite[Theorem 1.5.4]{Remy_2002}. The restriction $\res{N_K(T)}{\pi}$ is therefore still surjective. From \cite[Lemma 3.27 (ii)]{FHHK_2017} it now follows that 
$$
\nicefrac{N_G(T)}{T} = \nicefrac{A \rtimes N_K(T)}{(A \rtimes M)} \cong \nicefrac{N_K(T)}{M},
$$
and therefore the kernel of $\res{N_K(T)}{\pi}$ is equal to $M$. Thus 
$$
\nicefrac{N_K(T)}{M} \cong \nicefrac{N_G(T)}{T} \cong \W.
$$
On the other hand, one calculates
$$
\nicefrac{\EW M}{M} \cong \nicefrac{\EW}{(\EW \cap M)} = \nicefrac{\EW}{M^\dagger} \cong \W.
$$
For the first isomorphism the isomorphism theorem is used, for the second the definition of $M^\dagger$ with the fact that $\EW \cap T = M^\dagger$ is used, and for the last isomorphism \Cref{Eq:EW-M-Isom} is used.\\
From this follows
$$
\nicefrac{\EW M}{M} \cong \nicefrac{N_K(T)}{M},
$$
which shows the desired statement.
\qed
\end{Proof}

Finally, we can give the statement which compares the Weyl group of a real Kac--Moody group with the Weyl group of the corresponding symmetric space.

\begin{Proposition}\label{Prop:GeomW-AlgW-Isom}
Let $G$ be a split Kac--Moody group over $\RR$ and $(p,F)$ a pointed flat in the symmetric space $G/K$. Then $\Stab_G(p,F) \cong N_K(T)$ and $\Fix_G(p,F) \cong M$, moreover 
$$
    W \hspace{-2pt} \left( G \curvearrowright G/K \right) \cong \W .
$$
\end{Proposition}
\begin{Proof}
Since the definition of the Weyl group of a symmetric space is independent of the choice of the pointed flat, we can choose the standard flat $(e, AK)$.\\
Now we analyze the stabilizer and fixator of $G/K$ with respect to the action of $G$ on $G/K$. Observe that the fixator of $e$ is given by $K$ and $N_G(T)$ is the stabilizer of $AK$ (see \Cref{Co:NG(T)=Stab(AK)}). Thus $\Stab_G(e,AK) = N_K(T)$ follows and by \Cref{Lem:EW=NK(T)} follows that $\EW M = N_K(T)$.\\
Recall from \Cref{Lem:Eigenschaften-A_M} \ref{itm:A2} that $T = M \times A$, in particular $M$ centralizes $A$. Hence $M$ is a subgroup of $\Fix_G(e,AK)$.
Using \Cref{Lem:EW=NK(T)} and the isomorphism \Cref{Eq:EW-M-Isom} we conclude that 
$$
\nicefrac{\Stab_G(e,AK)}{M} = \nicefrac{N_K(T)}{M} = \nicefrac{\EW M}{M} \cong \nicefrac{\EW}{M^\dagger} \cong \W.
$$
So far we have established the isomorphism $\Stab_G(e,AK)/M \cong \W$. Since the representation $\rho_{KM} \colon \W \to \GL(\aa)$ (cf.\ \Cref{Eq:Rho-KM}) is faithful, each element different from the identity of $\W$ acts non-trivially on $\aa$. The isomorphism $\nicefrac{\Stab_G(e,AK)}{M} \cong \W$ from above hence implies  that $M = \Fix_G(e,AK)$. 
\qed
\end{Proof}

Note that from these isomorphisms used in the previous proof it can be deduced that any non-trivial coset of $\Stab_G(e,AK)/M$ acts non-trivially on the flat $AK$.

\subsection{A Local to Global Approach}
After we have discussed the local and global action of a Kac--Moody group on the corresponding symmetric space, we can try to connect them. In fact, one can do this using the Weyl groups and the fact that they are isomorphic, see \Cref{Prop:GeomW-AlgW-Isom}. The following statement comes from a theorem given for non-affine real split Kac--Moody groups, see \cite[Theorem 6.25]{FHHK_2017}. Here we generalize it to Kac--Moody groups of arbitrary type, and since we use \Cref{Prop:GeomW-AlgW-Isom}, our given assumption is that $G$ is a split real Kac--Moody group of arbitrary type $\AA$, where $\AA$ is non-spherical, indecomposable and symmetrizable. This enables us to use all the local action statements.
Before we state the theorem, it should be noted that a Coxeter diagram can be obtained from a Dynkin diagram $\Gamma$ by forgetting the double edges. Therefore, each diagram automorphism of $\Gamma$ induces an automorphism of the corresponding Coxeter diagram.\\
In summary, the group $\Aut(\Gamma)$ embeds in $\Aut(\W,S)$ and automorphisms of the embedded subgroup must preserve $S$ as a set.

\begin{Theorem}
Let $G$ be a split real Kac--Moody group of arbitrary type $\AA$, where $\AA$ is non-spherical, indecomposable and symmetrizable. Further denote by $G^\dagger$ the centered subgroup of $G$ and with $X^\dagger = G^\dagger/K^\dagger$ the corresponding Kac--Moody symmetric space and let $(p,F)$ be a pointed maximal flat in $X^\dagger$. Then the following diagram commutes:
\begin{center}
\tikzset{every picture/.style={line width=0.75pt}} 
\begin{tikzpicture}[x=0.75pt,y=0.75pt,yscale=-1,xscale=1]
\draw (42,22) node    {$W(\Aut(X^\dagger) \curvearrowright X^\dagger)$};
\draw (294.67,22) node    {$\GL\left( p,F,\sing{F}(p)\right)$};
\draw (35,120) node    {$\left( \W \rtimes \Aut( \Gamma )\right) \times \mathbb{Z} /2\mathbb{Z}$};
\draw (294.5,120) node    {$\mathbb{R}_{ >0} \times \left( \W \rtimes \Aut( \W ,S)\right) \times \mathbb{Z} /2\mathbb{Z}$,};

\draw (135,107) node [anchor=north west][inner sep=0.75pt]    {$\kappa$};
\draw (165,10) node [anchor=north west][inner sep=0.75pt]    {$\iota$};

\draw (41,62.4) node [anchor=north west][inner sep=0.75pt]    {$\cong $};
\draw (298,62.4) node [anchor=north west][inner sep=0.75pt]    {$\cong $};
\draw    (294.64,37.67) -- (294.64,102.5) ;
\draw [shift={(294.64,104.5)}, rotate = 270.1] [color={rgb, 255:red, 0; green, 0; blue, 0 }  ][line width=0.75]    (10.93,-3.29) .. controls (6.95,-1.4) and (3.31,-0.3) .. (0,0) .. controls (3.31,0.3) and (6.95,1.4) .. (10.93,3.29)   ;

\draw    (110,22) -- (225,22) ;
\draw [shift={(228,21.93)}, rotate = 539.9200000000001] [color={rgb, 255:red, 0; green, 0; blue, 0 }  ][line width=0.75]    (10.93,-3.29) .. controls (6.95,-1.4) and (3.31,-0.3) .. (0,0) .. controls (3.31,0.3) and (6.95,1.4) .. (10.93,3.29)   ;

\draw    (57.33,36) -- (57.33,106) ;
\draw [shift={(57.33,108)}, rotate = 270] [color={rgb, 255:red, 0; green, 0; blue, 0 }  ][line width=0.75]    (10.93,-3.29) .. controls (6.95,-1.4) and (3.31,-0.3) .. (0,0) .. controls (3.31,0.3) and (6.95,1.4) .. (10.93,3.29)   ;

\draw    (115,120) -- (173,120) ;
\draw [shift={(176,120)}, rotate = 540] [color={rgb, 255:red, 0; green, 0; blue, 0 }  ][line width=0.75]    (10.93,-3.29) .. controls (6.95,-1.4) and (3.31,-0.3) .. (0,0) .. controls (3.31,0.3) and (6.95,1.4) .. (10.93,3.29)   ;
\end{tikzpicture}
\end{center} 
Here $\kappa$ and $\iota$ denote canonical inclusion. Further, any local automorphism extends to a global automorphism if and only if $\Aut(\Gamma) = \Aut(\W,S)$.
\end{Theorem}

\begin{Proof}
The map $\iota$ is the natural inclusion since the Weyl group of the symmetric space takes values in the group of local transformations. The isomorphism of the right arrow follows directly from \Cref{Prop:GL-Isom}.\\
Now we are going to prove the isomorphism of the left arrow and there we use the same strategy as in the proof of Theorem 6.25 from \cite{FHHK_2017}. Therefore, recall that 
$$
W(\Aut(X^\dagger) \curvearrowright X^\dagger) \coloneqq \nicefrac{\Stab_{\Aut(X^\dagger)}(p,F)}{\Fix_{\Aut(X^\dagger)}(p,F)}.
$$
Since this definition of the Weyl group corresponding to the symmetric space is independent of the pointed flat (see \cite[Proposition 2.36]{FHHK_2017}, which is proved for abstract symmetric spaces due to Loos), one can consider the standard pointed flat $(e,A^\dagger K^\dagger)$. 

Using the isomorphism given in \Cref{Thm:X1-Isom}, one can use the decomposition to check how the single parts acts on the standard flat. For this purpose, recall the explanation of the different automorphism types given in \Cref{Sec:GlobalAction} or see \cite[Section 8.2]{CM_2005}. Furthermore, \Cref{Thm:X1-Isom} yields the isomorphism
$$
\Aut(X^\dagger) \cong \AutEff(G^\dagger) = \left( \nicefrac{G^\dagger}{C_{K^{\dagger}}(G^\dagger)} \rtimes (D \times \langle \Theta \rangle) \right) \rtimes \Aut(\Gamma).
$$
Now we use this decomposition of $\Aut(X^\dagger)$ to determine the stabilizer and the fixator by considering each factor separately.\\
Since the diagram automorphism $\Aut(\Gamma)$ stabilizes the standard flat $(e,A^\dagger K^\dagger)$ it remains in the decomposition and based on the fact that it is a subgroup of $\Aut(\W,S)$, one has the inclusion $\kappa$ in this component.
The diagonal automorphism fixes the standard flat pointwise and hence has no effect, i.e.\ the factor vanishes.\\  
The Chevalley involution $\Theta$ acts on the standard flat by inversion and so it corresponds to $\ZZ/2\ZZ$. Note that the semidirect product with $\Theta$ vanishes in the decomposition of the stabilizer in $\AutEff(G^\dagger)$ as one considers the action on $A^\dagger K^\dagger$ and by the proof of \cite[Proposition 6.5]{FHHK_2017} it is known that the Chevalley involution commutes with $\Aut(\Gamma)$ and it can be concluded that $\Theta$ also commutes with the inner automorphisms. 
In fact, since one considers the stabilizers of the standard flat $A^\dagger K^\dagger$, one knows that $\Theta$ commutes with the diagram automorphism and normalizes the left multiplication. If we denote the left multiplication with $g \in \Stab_{G^\dagger}\hspace{-2pt}(e,A^\dagger K^\dagger)$ by $l(g)$, then the following is true
$$
\Theta \circ l(g) = l(\Theta(g)) \circ \Theta.
$$
After this, it remains to analyze the first factor in the decomposition of $\Aut(X^\dagger)$. Note that the inner automorphism of $\AutEff(G^\dagger)$ corresponds to $G^\dagger / C_{K^\dagger}(G^\dagger)$, cf.\ \Cref{Thm:X1-Isom}. Further, by \Cref{Prop:GeomW-AlgW-Isom} one can consider the stabilizer and fixator of $A^\dagger K^\dagger$ as subgroups of $G^\dagger$:
$$
\Stab_{G^\dagger}\hspace{-2pt}(e,A^\dagger K^\dagger) \cong N_{K^{\dagger}}(T^\dagger) \quad \text{and} \quad \Fix_{G^\dagger}\hspace{-2pt}(e,A^\dagger K^\dagger) \cong M^\dagger.
$$
Using \Cref{Prop:GeomW-AlgW-Isom} and the isomorphism theorem, the following can be calculated
\begin{align*}
    \nicefrac{\Stab_{G^\dagger}\hspace{-2pt}(e,A^\dagger K^\dagger)/C_{K^{\dagger}}\hspace{-2pt}(G^\dagger)}{\Fix_{G^\dagger}\hspace{-2pt}(e,A^\dagger K^\dagger)/C_{K^{\dagger}}\hspace{-2pt}(G^\dagger)} & \cong \nicefrac{\Stab_{G^\dagger}\hspace{-2pt}(e,A^\dagger K^\dagger)}{\Fix_{G^\dagger}\hspace{-2pt}(e,A^\dagger K^\dagger)} \\
    & \cong \W.
\end{align*} 
So we have shown the isomorphisms of the vertical arrows, the commutativity of the diagram follows from the fact that the horizontal arrows in fact are inclusions.
\qed
\end{Proof}


\section{Comparing Symmetric Spaces} \label{comparing}
In the preceeding sections we constructed symmetric spaces associated to Kac--Moody groups over $\KK \in \{ \RR, \CC \}$ and with respect to arbitrary symmetrizable generalized Cartan matrices. Now we can give an summary and compare it with the symmetric space developed in \cite{FHHK_2017}, which was set up for non-affine Kac--Moody groups over $\RR$.\\
In \cite{FHHK_2017} the authors consider only real centered Kac--Moody groups with respect to the root datum $\D_{sc}^\AA$ (see \cite[Example 7.11]{Marquis_2018}). This choice implies that their coroots are linearly independent, but this does not mean that the roots are linearly independent. For this reason, the authors in \cite{FHHK_2017} work with the derived Kac--Moody algebra (see \cite[Definition 7.13]{Marquis_2018}), and use the Cartan subalgebra of this Kac--Moody algebra to construct the standard flat of the associated symmetric space.\\
If we now consider Kac--Moody groups with repsect to arbitrary generalized Cartan matrices, one problem we encounter is that Cartan matrices do not have to be invertible. In the case of a non-invertible generalized Cartan matrix, there is a non-trivial kernel of bilinear form the Cartan subalgebra of the derived Kac--Moody algebra, which is equal to the center (\cite[Lemma 2.1 a)]{Kac_1990}).\\
To work around this problem for non-affine Cartan matrices, the authors of \cite{FHHK_2017} use the semisimple adjoint quotient of a Kac--Moody group: this is the quotient of the group by the center of the Kac--Moody algebra that is lifted onto the group, cf. \cite[Definition 3.8]{FHHK_2017}. This quotient allows them to prove, for example, the useful statements \cite[Proposition 6.4]{FHHK_2017} and \cite[Theorem 6.12 (ii)]{FHHK_2017}.\\
If we also consider Kac--Moody symmetric spaces over $\CC$, we need a topological Iwasawa decomposition to define flats. At the time \cite{FHHK_2017} was created, no complex version of the topological Iwasawa decomposition was known. 

\smallskip

Here, we are talking about Kac--Moody groups over $\RR$ and $\CC$ with respect to the root datum $\D_{\Kac}^\AA$, where $\AA$ is an arbitrary generalized Cartan matrix. In this case the linear independence of the roots is given by enlarging the Cartan subalgebra, which leads to an enlarged torus at the group level.\\
This larger torus turns out to be a problem, because for the analysis of the global and local action of the Kac--Moody group on the corresponding symmetric space it is important that the two automorphism groups are related (see \Cref{Thm:X1-Isom}). The reason for this problem is Caprace's theorem: It is the basis for understanding the automorphisms of a particular type of Kac--Moody group. But it only describes the automorphisms of a Kac--Moody group generated by its root groups, see \cite[Theorem 4.2]{Caprace_2009}. On the extended torus one essentially asks for the automorphism group of $(\RR,+)$. This group is too big to be classified, so it is not possible to say anything about the automorphisms on the torus outside of $G^\dagger$.\\
Therefore, in order to study the local and global action, this paper is restricted to the subgroup $G^\dagger$ of a Kac--Moody group $G$, which is generated only by its root groups. Caprace's theorem applies to this centered subgroup $G^\dagger$ and leads to \Cref{Thm:X1-Isom}.\\
Regarding the case of generalizing the topological Iwasawa decomposition to complex split Kac--Moody groups, there is now a way to do this thanks to \cite[Corollary B.8]{HK_2021}. The quoted statement allows us to use the same methods as in \cite{FHHK_2017} to prove the topological Iwasawa decomposition in the complex case, see \Cref{Thm:Iwasawa}.\\
So far, we have achieved a generalization of the Kac--Moody symmetric space described in \cite{FHHK_2017}.


\section{Almost Split Kac--Moody Groups}
After developing a theory of symmetric Kac--Moody spaces over $\CC$ of arbitrary type $\AA$, where $\AA$ is an indecomposable, symmetrizable generalized Cartan matrix, it is of interest to ask about symmetric spaces of real forms of a Kac--Moody group. Thanks to the work of B. R\'{e}my, \cite{Remy_2002}, there is a theory of real forms in the context of Kac--Moody groups. So in this section we first give a brief overview of the theory and the definitions of the necessary notions and then describe how to develop a theory of almost split Kac--Moody symmetric spaces. This section is mainly based on \cite{Remy_2002}.

\subsection{Real Forms of Kac--Moody Groups}\label{Sec:RR-Forms}
Following \cite[Chapter 11]{Remy_2002}, we define several types of a real form. 

\begin{Definition}[cf.\ Definition 11.1.2 from \cite{Remy_2002}]
\begin{itemize}
    \item[] 
    \item Let $\KK \in \{ \RR, \CC \}$. Call a functor $F \colon \KK\text{-alg} \to \Grp$ a \emph{$\KK$-functor}.
    \item An $\RR$-functor $F$ is called a \emph{$\RR$-form} of a $\CC$-functor $H$, if there is a functorial isomorphism $\res{\CC}{F} \cong H$, i.e.\ the functor $F$ applied to the field $\CC$ is functorial isomorphic to the functor $H$.
\end{itemize}
\end{Definition}

The Tits functor $\G$ is clearly an example of an $\RR$-form of itself.\\
Let $F$ be an $\RR$-form of a $\CC$-functor $H$. Now, we want to define an action of the Galois group $\Gamma \coloneqq \operatorname{Gal}(\CC/\RR)$ on $H(\CC)$ using the real form $F$. Therefore denote by $R_\RR$ an $\RR$-algebra and by $R_\CC = R_\RR \otimes_\RR \CC$ the corresponding $\CC$-algebra given by scalar extension. Then for all $\sigma \in \Gamma$ the map
$$
1_{R_\RR} \otimes \sigma \colon R_\CC \to R_\CC
$$
is an isomorphism of $\RR$-algebras. By the fact, that the real form $F$ equals the $\CC$-functor $H$ over $\CC$, the isomorphism $1_{R_\RR} \otimes \sigma$ induces an automorphism $f_\sigma$ on $H(R_\CC)$. This leads to the following commutative diagram.

\begin{center}
\begin{tikzcd}[row sep = huge, column sep = huge]
F(R_\CC) \arrow{r}{F(1_{R_\RR} \otimes \sigma)} \arrow[swap, "=",d, leftrightarrow] & F(R_\CC) \arrow[d, leftrightarrow, "="] \\
H(R_\CC) \arrow[swap]{r}{f_\sigma} & H(R_\CC) 
\end{tikzcd}
\end{center}

To simplify the notation, we write $\sigma$ instead of $f_\sigma$. Note, that this action is a natural transformation, for more details see \cite[11.1.2]{Remy_2002}.

\begin{Definition}
Let $G$ be a Kac--Moody group over $\CC$ of type $\AA$, where $\AA$ is a indecomposable, symmetrizable generalized Cartan matrix. Further denote by $\G$ the Tits functor of type $\D_{\Kac}^\AA$, i.e.\ $\G(\CC) = G$. Then call an $\RR$-form $\F$ of the Tits functor $\G$ a \emph{functorial $\RR$-form of the Kac--Moody group $G$}.
\end{Definition}

Recall, that a Kac--Moody algebra $\g(\AA)$ has an associated universal enveloping algebra $\U_{\g(\AA)}$, see \cite[Chapter 3.1]{Marquis_2018} or \cite[7.3]{Remy_2002}. Following \cite[7.4.5]{Remy_2002} we denote by $\U$ a corresponding $\ZZ$-form of $\U_{\g(\AA)}$, i.e.
$$
\U \otimes_\ZZ \CC \cong \U_{\g(\AA)}.
$$
Denote by $\U_\RR$ the scalar extension of $\U$ with respect to $\RR$, i.e.\ $\U \otimes_\ZZ \RR$. Note, that the universal enveloping algebra can be seen as a functor, which is very useful to define the notion of a prealgebraic $\RR$-form:\\
The pair $(\F, \U_\RR)$, where $\F$ is a functorial $\RR$-form of a Kac--Moody group $G$ and $\U_\RR$ is an $\RR$-form of the universal enveloping algebra $\U_{\g(\AA)}$, is called \emph{prealgebraic $\RR$-form of $G$} under the following conditions:
\begin{itemize}
    \item The adjoint representation $\Ad \colon \G \to \Aut(\U)$, where $\Aut(\U) \colon \ZZ\text{-alg} \to \Grp$ is a functor, is equivariant with respect to the Galois action.
    \item For each subfield $\KK$ of $\CC$, $\F(\KK)$ is a subgroup of $\G(\CC) = G$.
\end{itemize}
For more details about prealgebraic $\RR$-froms see \cite[11.1.3]{Remy_2002} and for more details about the adjoint representation see \cite[9.5]{Remy_2002}.

\smallskip

From \cite[Proposition 11.2.2]{Remy_2002} follows that the Galois action on a Kac--Moody group $G$ stabilizes the set of tori in $G$ and that for every Galois action $\sigma$ there is an element $g \in G$, so that one can \emph{rectify} the action in the sense that the standard torus is fixed. Denote with $c_g(x)$ the map that conjugates an element $x \in G$ with $g$. Then we can state that the action
$$
c_{g^{-1}} \circ \sigma
$$
is well-defined and fixes the standard torus $T \leq G$. In the following denote by $\os$ the \emph{rectification of $\sigma$}.\\
Since a root $\a \in \DD$ is defined with respect to the Cartan subalgebra $\h$, see \Cref{Sec:KM_Algebra}, one can speak of a root $\a$ \emph{relative to} $\h$. In this context, we can also call a root group $U_\a$ \emph{relative to $T$} if the root $\a$ is relative to the Cartan subalgebra $\h$ underlying $T$.

\begin{Convention}[cf.\ Hypothesis 11.2.3 from \cite{Remy_2002}.]
Denote by $\SGR$ the condition, that a rectified Galois action of a Kac--Moody group maps a root group relative to $T$ to a root group of the same type.  
\end{Convention}

Condition $\SGR$ induces an action of the rectified Galois action on the set of real roots $\DD$. To make use of this action, we give the following definition, cf.\ \cite[11.2.5]{Remy_2002}.\\
Recall that $\U_\RR$ is the universal enveloping algebra and therefore satisfies the same grading as the roots of the underlying Kac--Moody algebra with respect to $Q$. For details of the grading of the universal enveloping algebra see \cite[Chapter 3.1]{Marquis_2018} and for the grading of the roots with respect to $Q$ see \Cref{Eq:RootLattice}.

\begin{Definition}
A prealgebraic $\RR$-form $(\F, \U_\RR)$ of a Kac--Moody group $G$ over $\CC$ of type $\AA$, where $\AA$ is an indecomposable, symmetrizable generalized Cartan matrix, is called an \emph{algebraic $\RR$-form} of $G$ if 
\begin{itemize}
    \item the condition $\SGR$ is satisfied,
    \item for every rectified Galois action $\os$ (i.e.\ $\os(T) = T$) the following is true
    \begin{itemize}
        \item the action $\os$ respects the decomposition into $Q$-homogeneous components of $\U_\RR$ (i.e.\ it respects the $Q$-gradation) and the induced permutation on $Q$ satisfies the following homogeneity condition:\\
        For all $\a \in \DD$ and for all $n \in  \NN$ we have $\os(n\a) = n \os(\a)$.
        \item the action $\os$ stabilizes $\Lambda $ and $\Lambda^\vee$.
    \end{itemize}
\end{itemize}
Furthermore, if the Galois action stabilizes the conjugacy class of the positive Borel subgroups of $G$, call the algebraic $\RR$-form an \emph{almost split Kac--Moody $\RR$-group}.
\end{Definition}

Note that by \cite[Lemma 11.3.2]{Remy_2002} one can rectify a Galois automorphism $\sigma$ in such a way that the sign of the Borel subgroup is preserved. Moreover, this rectification also stabilizes the standard torus $T$ and is well-defined up to $T$. Thus, this rectification is just a special case of the rectified Galois automorphism we defined above. In particular, this action induces an action on the Weyl group $\W$ and stabilizes the set of generators $S$, see \cite[Proposition 11.3.2]{Remy_2002}. Since we are only interested in almost split Kac--Moody $\RR$-groups and for simplicity, we denote this rectification by $\os$. For more details about almost split Kac--Moody $\RR$-groups see \cite[11.3.1]{Remy_2002}.

\smallskip

Before we discuss Galois descent, we give useful properties of the rectified Galois action. Let $G$ be a Kac--Moody group over $\CC$ of type $\AA$, where $\AA$ is indecomposable and symmetrizable. In addition, let $(\F, \U_\RR)$ be an algebraic $\RR$-form of $G$. Then any rectified Galois automorphism $\os$ has the following properties:
\begin{myitems}
    \item \label{itm:Wurzel_Wirkung} For any real root $\a \in \DD$ and any parameter $t \in \CC$, there is the equation
    $$
    \os(x_\a(t)) = x_{\os(a)}(k_\a \sigma(t)),
    $$
    where $k_\a \in \CC^\times$ is an element defined by the Galois action on a generator of the Kac--Moody algebra.

    \item The reflection $r_\a \in \W$ is mapped to $r_{\os(\a)}$.

    \item \label{itm:T_Wirkung} For any real root $\a \in \DD$ the element $h_\a \in \Lambda^\vee$ is mapped to $h_{\os(\a)}$.
\end{myitems}
For details, we refer to \cite[Proposition 11.2.5]{Remy_2002}.

\subsection{The Galois Descent}\label{Sec:GaloisDescent}
In linear algebra it is well known that a real form of a complex vector space can be identified with the fixed point set of an action given via Galois automorphism. R\'{e}my uses the same idea to develop further the theory of almost split Kac--Moody $\RR$-groups.\\
In \cite[Chapter 12]{Remy_2002}, he uses the fact that one can associate to a Kac--Moody group a twin building on which he can also define a Galois action which is compatible with the Galois automorphism action on the group. In detail, he uses a geometric realization to study the Galois action and for example to check if there exists at all fixed points, for the realization see \cite[Chapter 5]{Remy_2002} and for the action on the realization of the twin building see \cite[12.1.2]{Remy_2002}.\\
The theory R\'{e}my develops in Chapter 12 is very technical and uses a lot of new notions, which we do not need in detail, hence, we only give a short overview. 
First let $G$ be a Kac--Moody group over $\CC$ of type $\AA$, where $\AA$ is an indecomposable, symmetrizable generalized Cartan matrix. Assume that $G$ has an algebraic $\RR$-form $(\F, \U_\RR)$ which satisfies $\SGR$, i.e.\ $G$ is an almost split Kac--Moody $\RR$-group. Then call 
$$
G(\RR) \coloneqq \F(\RR)
$$
the \emph{rational points of $G$}. Furthermore, $G$ satisfies the \emph{descent condition \eqref{Eq:DCS}} if 
\[ \label{Eq:DCS}
G(\RR) = G^\Gamma, \tag{\DCS} 
\]
where $\Gamma$ denotes the group of the induced Galois automorphism on $G$, see \Cref{Sec:RR-Forms}. Call a subgroup $H$ of $\F(\CC) = G$ an $\RR$-subgroup, if it is stable under the Galois action. The rational points of an $\RR$-subgroup $H \leq G$, denoted by $H(\RR)$ are given by the Galois fixed points of $H$.\\
To make the following terms well-defined, we need to define a triple $(A,C,-C)$, where $A$ is an apartment and $C,-C$ are opposite chambers in $A$. This is similar to defining a standard torus and defines the geometric realization. In our context, it is only necessary that we obtain the following when applying the theory of \cite[Chapter 12]{Remy_2002}:
\begin{itemize}
    \item The group of rational points possesses a root system, which is not necessarily reduced. A \emph{real $\RR$-root} is a restricted root in the context of the geometric realization. To be precise, a root (there it is called half-apartment) in the geometric realization is given by a hyperplane on which a reflection is defined. In the process of doing the Galois descent, one restricts the hyperplane to the fixed point set. For more details see \cite[Chapter 5]{Remy_2002}. There the geometric realization and the relative terminology is defined.\\
    Denote by $\Delta^{\rel}_{\mathrm{re}}$ the set of all real $\RR$-roots with respect to $(A, C, -C)$ (cf.\ \cite[5.5.1]{Remy_2002}). Denote by $\root{\Pi}$ the set of simple roots of $\Delta^{\rel}_{\mathrm{re}}$.
    
    \item In this setting, one can define $\RR$-root groups as follows: Let $\a^{\rel} \in \Delta^{\rel}_{\mathrm{re}}$ be a real $\RR$-root, then denote by $U_\a^{\rel}$ the $\RR$-root group with respect to $\a^{\rel}$. The group $U_\a^{\rel}$ is defined as the group generated by $\{ U_\beta \mid \beta \in \Delta^{\rel}_\a \}$, where $\Delta^{\rel}_\a \coloneqq \{ \beta \in \DD \mid \exists \lambda \in \{ 1,2 \} \colon \beta^{\rel} = \lambda \a^{\rel} \} $. Denote by $\rootgrp{\a}$ the rational points of this subgroup, i.e.\ $\rootgrp{\a} \leq G(\RR)$, and call it \emph{relative root group} (cf.\ \cite[Definition 12.3.3]{Remy_2002}). By \cite[Lemma 12.3.4]{Remy_2002} follows, that any relative root group with respect to the real $\RR$-root $\a^{\rel}$ is non-trivial.

    \item By \cite[Theorem 12.6.3]{Remy_2002} it follows, that the group of rational points $G(\RR)$ of an almost split Kac--Moody $\RR$-group $G$ possess an RGD system.
\end{itemize}

\subsection{Relative Kac--Peterson Topology}

As the Kac--Peterson topology is defined for a Kac--Moody group $G$, see \Cref{Sec:TitsFunctor} for the definition, we can define a \emph{relative Kac--Peterson topology} on the rational points $G(\RR) \leq G$ of an almost split Kac--Moody $\RR$-group. The question is whether the relative Kac--Peterson topology coincides with the subspace topology induced by the surrounding Kac--Moody group. This section is based on \cite[8.2.3]{Diss}.\\
Let $G$ be an almost split Kac--Moody $\RR$-group which satisfies \eqref{Eq:DCS} and denote by $G(\RR)$ the rational points of $G$. In order to define the relative Kac--Peterson topology, we need that the relative root groups $\rootgrp{\a}$ (cf.\ \Cref{Sec:GaloisDescent}) have a Lie group structure and particularly a Lie group topology. From \cite[Lemma 12.5.4]{Remy_2002} follows that relative root groups are isomorphic to root groups of semisimple $\RR$-groups. This gives us, that a relative root group is a Zariski closed subgroup of the general linear group and hence is a Lie group. Next, we need to define \emph{relative rank one subgroups}
$$
G^{\rel}_{\a_i} = \left\langle \rootgrp{\a_i} , \rootgrp{-\a_i} \right\rangle,
$$
where $\a_i \in \root{\Pi}$.\\
Even the torus $Z(\RR) \leq G(\RR)$, defined in \cite[Chapter 12.5]{Remy_2002}, is a subgroup of the torus of the surrounding split Kac--Moody group and splits over $\RR$, see \cite[12.5.1]{Remy_2002} and \cite[12.5.3]{Remy_2002}. This means, that $Z(\RR)$ also carries a Lie group topology.\\
Now we have all the ingredients to define the relative Kac--Peterson topology in the same way as the Kac--Peterson topology, comparing to \Cref{Sec:TitsFunctor}.

\begin{Definition}
Let $G$ be an almost split Kac--Moody $\RR$-group satisfying \eqref{Eq:DCS} and denote by $G(\RR)$ the rational points of $G$. Then the \emph{relative Kac--Peterson topology} $\alsp{\topo{KP}}$ on $G(\RR)$ is the finest group topology on $G(\RR)$ such that the embeddings of the relative rank one groups and the embedding of $Z(\RR)$ are continuous with respect to their induced Lie topologies.
\end{Definition}

In the proof of \cite[Proposition B.3]{HK_2021} bounded subgroups are used to define a topology on a split Kac--Moody group $G$ which equals the Kac--Peterson topology. This can also be done for rational points of an almost split Kac--Moody $\RR$-group $G(\RR)$. To define the notion of a restricted subgroup in $G(\RR)$, we recall what is a parabolic subgroup in the context of a split Kac--Moody group $G$ over $\CC$, see \cite[p.150]{Marquis_2018}:\\
Denote by $\W_J$ the subgroup of the Weyl group $\W$ given by $\W_J = \langle r_i \mid i \in J \rangle$ for a subset $J \subseteq I$, where $I = \{ 1, \ldots, n \}$ , see \cite[Definition 2.12]{Abramenko_2010}. Then a parabolic subgroup $P_J$ is defined by
$$
P_J = B^+ \W_J B^+.
$$
If $\W_J$ is finite, call $P_J$ spherical. In the language of buildings, a parabolic subgroup is the stabilizer in $G$ of the standard residue $J$ in the positive half of the twin building $\Delta_\pm$. 

\begin{Definition}
Let $G$ be a split Kac--Moody group over $\KK \in \{ \RR, \CC \}$ of type $\AA$, where $\AA$ is an indecomposable, symmetrizable generalized Cartan matrix. Call a subgroup $M \leq G$ \emph{bounded} if it is contained in the intersection of two spherical parabolic subgroups with opposite signs.    
\end{Definition}

Since bounded subgroups are subgroups contained in spherical parabolic subgroups, they carry a Lie group structure and thus a Lie group topology. 

\begin{Definition}\label{Def:MBS-Topologie}
Let $G$ be a split Kac--Moody group over $\KK \in \{ \RR, \CC \}$ of type $\AA$, where $\AA$ is an indecomposable, symmetrizable generalized Cartan matrix. Denote by $\mbs{\KK}$ the set of maximal bounded subgroups. The \emph{maximal bounded subgroup topology} $\topo{MB}$ is the finest group topology, so that the embeddings of the maximal bounded subgroups in $G$ are continuous.
\end{Definition}

Note, that \cite[Proposition B.3]{HK_2021} is formulated for centered split Kac--Moody groups, thus we now prove a version for arbitrary split Kac--Moody groups.

\begin{Proposition}\label{Prop:B3-Paula}
Let $G$ be a split Kac--Moody group over $\KK \in \{ \RR, \CC \}$ of type $\AA$, where $\AA$ is an indecomposable, symmetrizable generalized Cartan matrix. Equip the standard torus with the Lie topology (see discussion after \Cref{Conv:Non-Degenerate-Tits-Functor} or \Cref{Prop:Torus-Exp-Topologisch}), then the Kac--Peterson topology $\topo{KP}$ equals the maximal bounded subgroup topology $\topo{MB}$.
\end{Proposition}

\begin{Proof}
We follow the strategy of the proof given for \cite[Proposition B.3]{HK_2021}. Since the cited proposition is given for centered Kac--Moody groups and in general a Kac--Moody group $G$ can be decomposed in $G = T G^\dagger$, see \Cref{Rem:RGD}, where $G^\dagger$ is a centered Kac--Moody group, see \Cref{Rem:G1-is-KM}, we have to take a closer look at the torus.\\
Note, that by definition the torus is a bounded subgroup,
$$
T = B_+ \cap B_-,
$$
and carries a Lie group topology, see discussion after \Cref{Conv:Non-Degenerate-Tits-Functor} or \Cref{Prop:Torus-Exp-Topologisch}. Thus the torus embeds as a closed subgroup into a maximal bounded subgroup and by the closed-subgroup theorem (\cite[Theorem 9.3.7]{HN_2012}) follows that the subspace topology of the embedded torus equals its Lie group topology.\\
As a next step, we apply \cite[Lemma 4.3]{Marquis_2012}, which states that the subspace topology induced by $\topo{KP}$ on the maximal bounded subgroups is equal to the Lie group topology. Although the statement is given for centered split Kac--Moody groups, it also applies to all split Kac--Moody groups, since the statement is based on building theory arguments.\\
Recall, that the Kac--Peterson topology induces the Lie group topology on the rank one subgroups $G_\a \cong \SL_2(\KK)$, $\a \in \DD$, see \cite[Definition and Remark 2.1]{HK_2021}. Furthermore, the rank one subgroups are bounded subgroups and a rank one subgroup embeds as a closed subgroups into a maximal bounded subgroup. By the closed-subgroup theorem (\cite[Theorem 9.3.7]{HN_2012}) follows that the subspace topology of the embedded rank one subgroups in a maximal bounded subgroup equals its Lie group topology.\\
As a result we have 
$$
\topo{MB} \subseteq \topo{KP}.
$$
To see the other direction, consider the maximal bounded subgroups in $(G, \topo{KP})$. Here, each maximal bounded subgroup has the Lie group topology induced by the subspace topology (cf.\ \cite[Lemma 4.3]{Marquis_2012}). 
The embedded relative rank one subgroups and the embedded torus also have the induced subspace Lie topology and one concludes that the embedding of the maximal bounded subgroups in $G$ is continuous with respect to the topology $\topo{KP}$. Therefore, the topology $\topo{KP}$ is coarser than $\topo{MB}$.
\qed
\end{Proof}

A parabolic $\RR$-subgroup of $G(\RR)$ is also a parabolic subgroup in the surrounding split Kac--Moody group $G$, see \cite[Definition 12.1.3]{Remy_2002} and \cite[Remark (1) 12.1.3]{Remy_2002}). Thus a maximal bounded $\RR$-subgroup embeds in a maximal bounded subgroup of $G$. By \cite[6.2]{Remy_2002} and \cite[Section 3.1]{Caprace_2009} follows that maximal bounded subgroups are (algebraic) Lie groups. 
To keep the notation simple, we denote the set of maximal bounded subgroups in $G(\RR)$ by $\mbs{\RR}$ and in $G$ by $\mbs{\CC}$.

\smallskip

Recall that a topological group is called \emph{almost connected} if the quotient of the group by the connected component of the neutral element is compact, see \cite{Marquis_2012}.\\
The action of an almost connected Lie group $H$ on the Davis realization of a building has a global fixed point, see \cite[Corollary B]{Marquis_2012}. Further, recall that $H/H_0$, where $H_0$ is the identity component of $H$, is compact and therefore has finitely many connected components. A relative rank one subgroup is a real Lie group and from \cite[Section 3.2, Corollary 1]{Platonov_1993} follows, that it has finitely many connected components. From the discussion above, we conclude that a relative rank one subgroup fixes a point in the Davis realizations of both halves of the twin building associated to the surrounding split Kac--Moody group $G$, i.e.\ $G^{\rel}_{\a} \subseteq \mbs{\CC}.$

From the closed-subgroup Theorem (cf.\ \cite[Theorem 9.3.7]{HN_2012}) and the fact that a relative rank one subgroup is a closed subgroup in a maximal bounded subgroup of $\mbs{\CC}$ because the Galois action on $G$ is continuous, we conclude that a relative rank one subgroup is an embedded Lie group. Moreover, the Lie group topology coincides with the subspace topology from $\mbs{\CC}$. Note that the subspace topology induced by $\topo{KP}$ on a maximal bounded group of $\mbs\CC$ is equal to the Lie group topology, see the proof of \Cref{Prop:B3-Paula}.\\
Now, we consider the topology that can be defined by the embedding of the maximal bounded subgroups. Recall the definition of $\topo{MBS}$, see \Cref{Def:MBS-Topologie}.

\begin{Definition}
Let $G$ be an almost split Kac--Moody $\RR$-group satisfying \eqref{Eq:DCS} and denote by $G(\RR)$ the rational points of $G$. Define the topology $\alsp{\topo{MB}}$ on $G(\RR)$ as the finest group topology so that the embeddings of the maximal bounded $\RR$-subgroups are continuous.
\end{Definition}

\begin{Lemma}
Let $G$ be an almost split Kac--Moody $\RR$-group satisfying \eqref{Eq:DCS} and denote by $G(\RR)$ the rational points of $G$. Then
\begin{align}\label{Eq:TopoMB}
    \alsp{\topo{MB}} = \res{G(\RR)}{\topo{MB}},
\end{align}
where $ \res{G(\RR)}{\topo{MB}}$ means the subspace topology induced by $\topo{MB}$.
\end{Lemma}

\begin{Proof}
Both topologies are defined as the finest topologies, moreover they coincide on the $\mbs{\RR}$:\\
Every parabolic $\RR$-subgroup of $G(\RR)$ is a parabolic subgroup of $G$. From the previous discussion, we know that the maximal bounded subgroups are Lie groups. This leads to the fact that $\mbs{\RR}$ are also Lie groups, and in particular they are equal to $\mbs{\CC} \cap G(\RR)$, equipped with the trace topology.
\qed
\end{Proof}

In the last step, we want to compare the topologies. To do this, we need a different point of view for the Kac--Peterson topology:\\
The universal topology $\topo{Prod}$, defined in \cite[Definition 7.19]{HKM_2013} using colimits, is equal to the Kac--Peterson topology (cf.\ \cite[Proposition 7.21]{HKM_2013}) and by transport of structure, it is a $k_\omega$ topology, for details see \Cref{Sec:Iwaswa} or \cite{k-Omega}. The authors of the article \cite[Corollary 6.7]{Glockner_2008} state that $G(\RR)$ is $k_\omega$ with respect to the induced subspace topology coming from the surrounding split Kac--Moody group.

\smallskip

Putting all things from this section together like a puzzle, we observe that by \Cref{Prop:B3-Paula} $\topo{KP}$ is equal to $\topo{MB}$ for split Kac--Moody groups. Moreover, by \cite[Proposition 7.21]{HKM_2013}, one knows that $\topo{Prod} = \topo{KP}$. In particular 
\begin{align}\label{Eq:SplitTopoEqual}
    \topo{MB} = \topo{KP} = \topo{Prod}.
\end{align}

This can be now used to compare the different structures on $G(\RR)$.
\begin{Theorem}\label{Prop:SubspaceTopo-Equal}
Let $G$ be an almost split Kac--Moody $\RR$-group satisfying \eqref{Eq:DCS} and denote by $G(\RR)$ the rational points of $G$. Then the following holds for $G(\RR)$:
\begin{myitems}[label=\textbf{\Roman*.}, itemsep=0pt]
    \item $\alsp{\topo{MB}} = \res{G(\RR)}{\topo{MB}} = \res{G(\RR)}{\topo{Prod}} = \res{G(\RR)}{\topo{KP}}$,
    \item $G(\RR)$ is a $k_\omega$-space with respect to the topology $\alsp{\topo{MB}}$, and
    \item $\alsp{\topo{KP}} = \alsp{\topo{MB}}$.
\end{myitems}
\end{Theorem}

\begin{Proof}
\begin{myitems}[label=\textbf{\Roman*.:}, itemsep=0pt]
    \item[] 
    \item This statement follows from the discussion above with \Cref{Eq:TopoMB}.
    
    \item This statement follows from the previous discussion. 
    
    \item Denote by $(\mbs{\RR}_i)_{i \in I}$ the family of maximal bounded $\RR$-subgroups in $G(\RR)$ and consider the topological spaces $(G(\RR), \alsp{\topo{KP}})$ and $(G(\RR), \alsp{\topo{MB}})$. A relative rank one subgroup embeds as a Lie subgroup into one of the $\mbs{\RR}_i$. Moreover the continuous embeddings of the relative rank one subgroups into $G(\RR)$ provide fewer conditions than the embeddings of the maximal bounded subgroups. Hence the map from $(G(\RR), \alsp{\topo{KP}})$ to $(G(\RR), \alsp{\topo{MB}})$ is continuous, i.e.\ $ \alsp{\topo{MB}} \subseteq \alsp{\topo{KP}}.$\\
    For the other direction consider the maximal bounded subgroups in the topological group $(G(\RR), \alsp{\topo{KP}})$. Every $\mbs{\RR}_i$ carry the Lie group topology and the embedded relative rank one subgroups carry the induced Lie subspace topology. Thus one deduces that the embedding of $\mbs{\RR}_i$ into $G(\RR)$ with respect to the topology $\alsp{\topo{KP}}$ is continuous. To be more precise, the topology $\alsp{\topo{KP}}$ comes from the embedding of relative rank one subgroups equipped with the Lie topology and by the closed-subgroup theorem (cf.\ \cite[Theorem 9.3.7]{HN_2012}) the maximal bounded subgroups carries the same Lie group topology. This means that the topology $\alsp{\topo{KP}}$ is coarser than $\alsp{\topo{MB}}$. 
\end{myitems}
\qed
\end{Proof}


\section{Symmetric Spaces of Almost Split Kac--Moody Groups} \label{almostsplit}
In this section, we are going to define a symmetric space for almost split Kac--Moody $\RR$-groups.\\
Therefore, let $G$ be a Kac--Moody group over $\CC$ of type $\AA$, where $\AA$ is an indecomposable, symmetrizable generalized Cartan matrix of size $n$ with rank $l$. Denote by $\Gamma$ the group of Galois automorphisms induced by $\operatorname{Gal}(\CC/\RR)$ on $G$ and assume that $G$ possesses an algebraic $\RR$-form $(\F,\U_\RR)$ satisfying the \eqref{Eq:DCS} condition, i.e.\ 
$$
G(\RR) = G^\Gamma,
$$
where $G(\RR)$ denotes the rational points of the almost split Kac--Moody $\RR$-group $G$. Denote by $\os$ the rectified Galois automorphism which stabilizes the standard torus $T$ of $G$, i.e.\ $\os(T) = T$.

\begin{Lemma}\label{Lem:Sigma-Theta-Tauscht}
Let $G$ be an almost split Kac--Moody $\RR$-group. Any rectified Galois automorphism $\os$ can be conjugated by a diagonal automorphism such that it commutes with $\Theta$.
\end{Lemma}

\begin{Proof}
Similar to the proof of the algebraic Iwasawa decomposition, we use the decomposition
$$
G = T G^\dagger.
$$
This allows us to prove the statement separately for $T$ and $G^\dagger$. Note, that by \Cref{Sec:GlobalAction} every automorphism of $\Aut(G^\dagger)$ can be decomposed into an inner automorphism, a sign automorphism, a diagonal automorphism, a diagram automorphism and a field automorphism.\\
By \Cref{Sec:RR-Forms} \ref{itm:Wurzel_Wirkung} one knows how $\os$ acts on an element of a root group $U_\a$, $\a \in \DD$: 
$$
\os(x_\a(t)) = x_{\os(\a)}(k_a \overline{t}),
$$
where $x_\a(t)$, $t \in \CC$, denotes an element $U_\a$, see \Cref{Sec:TitsFunctor}. The operation to add $k_\a$ can be understood as the action of a diagonal automorphism. Hence $\os$ can be rectified by conjugation with a diagonal automorphism $d$ so that $k_\a$ becomes $1$. In detail denote by
$$
\widehat{\sigma} \coloneqq d \circ \os \circ d^{-1}
$$
the rectified Galois automorphism, where $d$ acts by $ d(x_\a(t)) = x_\a(d_\a t) $ such that $ k_\a = \frac{\overline{d_\a}}{d_{\os(a)}} $. This results in
\begin{align*}
    \widehat{\sigma}(x_\a(t)) &= (d \circ \os \circ d^{-1})(x_\a(t)) \\
    &= (d \circ \os)(x_\a(d_\a^{-1} t)) \\
    &= d (x_{\os(\a)}(k_\a \overline{d_\a^{-1}} \overline{t})) \\
    &= x_{\os(\a)}(d_{\os(\a)} \overline{d_\a^{-1}} k_\a \overline{t}) \\
    &= x_{\os(\a)}(\overline{t}).
\end{align*}
By \cite[p.384]{CM_2005} one can always choose such a diagonal automorphism with respect to $d_\a$ and $d_{\os(\a)}$. Moreover, $\widehat{\sigma}(T) = T$, since the diagonal automorphism acts as the identity on $T$ by definition (cf.\ \cite[p.384]{CM_2005}).\\
This results in the observation, that $\Theta$ and $\widehat{\sigma}$ commutes on $G^\dagger$:
\begin{align*}
    \left( \Theta \circ \widehat{\sigma} \right)\hspace{-2pt}(x_\a(t)) &= \Theta(x_{\os(\a)}(\overline{t})) \\
    &= x_{-\os(\a)}(t),
\end{align*}
and 
\begin{align*}
    \left( \widehat{\sigma} \circ \Theta \right)\hspace{-2pt}(x_\a(t)) &= \widehat{\sigma}\left( x_{-\a}(\overline{t}) \right) \\
    &= x_{-\os(\a)}(t).
\end{align*}

It remains to show, that $\Theta$ and $\widehat{\sigma}$ commutes on $T$ and since every diagonal automorphism is the identity on $T$, we conclude that $\widehat{\sigma}$ acts as $\os$ on $T$. The action of $\os$ on $T$ is given by \Cref{Sec:RR-Forms} \ref{itm:T_Wirkung}. Now given a generating set $ \lbrace r^{v_i} \mid r \in \CC^\times \, , \; 1 \leq i \leq 2n-l \rbrace $ of $T$. Then we calculate for any $\a \in \DD$
$$
\os\left(  r^{v_i} \right) = \overline{r}^{v_{\os(i)}}.
$$
In \Cref{Sec:Chevalley} we have defined the action of $\Theta$ on $T$ and hence we can compute 
\begin{align*}
    \left( \Theta \circ \os \right)\hspace{-2pt}\left( r^{v_i} \right) &= \Theta\left( \overline{r}^{v_{\os(i)}} \right) = r^{-v_{\os(i)}}, \\
\intertext{and on the other hand}
     \left( \os \circ \Theta \right)\hspace{-2pt}\left( r^{v_i} \right) &= \os\left( \overline{r}^{-v_i} \right) = r^{-v_{\os(i)}}.
\end{align*}
\qed
\end{Proof}

Based on the previous lemma, we observe that
$$
\Theta\hspace{-2pt}\left(G^{\Gamma} \right) = \Theta(G)^{\Gamma},
$$
where we denote with $\Gamma$ the group of rectified and conjugated Galois automorphism $\widehat{\sigma}$. Note that this shows that $G(\RR)$ is closed under the action of $\Theta$.\\
Therefore $\Theta$ is a continuous involution on the topological group $G(\RR)$ and we define
$$
K(\RR) \coloneqq K^{\Gamma} = \left( G^\Theta \right)^\Gamma = \left( G^\Gamma \right)^\Theta.
$$

Similar to the approach for split Kac--Moody groups, we can define the twist map on $G(\RR)$ by
\begin{align*}
    \alsp{\tau} \colon G(\RR) &\to G(\RR) \\
    g &\mapsto g \Theta \left( g ^{-1} \right),
\end{align*}
and also a reflection map by
\begin{equation}
    \begin{split}
        \alsp{\mu} \colon G(\RR)/K(\RR) \times G(\RR)/K(\RR) &\to G(\RR)/K(\RR) \\
        (gK(\RR), hK(\RR)) &\mapsto \alsp{\tau}(g)\Theta(h)K(\RR).
    \end{split}
\end{equation}

Note, that $\alsp{\tau}$ equals $\tau$ since $G(\RR)$ is a subgroup of $G$. Moreover, since $G(\RR)$ is closed with respect to $\Theta$ and $\alsp{\tau}$, the map $\alsp{\mu}$ satisfies all properties of an abstract symmetric space (cf.\ \Cref{Def:SymSp-Loos}). Even the fourth property $\textbf{S.4}_{\text{global}}$ is true, because 
$$
\alsp{\tau}(G(\RR)) \cap K(\RR) \subseteq \tau(G) \cap K = \{ e \}. 
$$

This results in the following theorem, which is a generalization of the split case. 

\begin{Theorem}\label{Thm:Galois-Sym-Sp}
Let $G$ be an almost split Kac--Moody $\RR$-group satisfying \eqref{Eq:DCS} and denote by $G(\RR)$ the rational points of $G$. The quotient 
$$
X(\RR) \coloneqq \nicefrac{G(\RR)}{K(\RR)}
$$
is a symmetric space in the sense of Loos with the reflection map $\alsp{\mu}$. Moreover, there is a natural action of $G(\RR)$ on the symmetric space given by automorphisms
\begin{align*}
    G(\RR) &\to \Sym(X(\RR)) \\
    g &\mapsto (hK(\RR) \mapsto ghK(\RR)).
\end{align*}
\end{Theorem}

Since $\F(\CC) = G$ if $G$ is a complex (split) Kac--Moody group, it can be deduced that the theorem is in this particular case just \Cref{Prop:KM-SymSp} for $\KK = \CC$. Furthermore, if we consider a real split Kac--Moody group in the theorem, the following theorem corresponds to \Cref{Prop:KM-SymSp} for $\KK = \RR$.

\begin{Proof}
Since the properties of an abstract symmetric space due to Loos are satisfied, we need to check that 
$$
G(\RR) \to \Sym(X(\RR)) \, , \; g \mapsto (hK(\RR) \mapsto ghK(\RR))
$$
is an action. Thus, let $g,h,d \in G(\RR)$ and do the following calculation
\begin{align*}
    \alsp{\mu}(ghK(\RR), gdK(\RR)) &= \alsp{\tau}(gh) \Theta(gd) K(\RR) \\
    &= gh \Theta(gh)^{-1} \Theta(g) \Theta(d) K(\RR) \\
    &= gh \Theta(h)^{-1} \Theta(g)^{-1} \Theta(g) \Theta(d) K(\RR) \\
    &= g ( \alsp{\tau}(h) \Theta(d) K(\RR)) \\
    &= g \alsp{\mu}(hK(\RR), dK(\RR)).
\end{align*}
\qed
\end{Proof}

\subsection{Flats in Almost Split Kac--Moody Symmetric Spaces}
After establishing the notion of a symmetric space for almost split Kac--Moody groups, we can define flats. As we have seen in the proof of \Cref{Prop:Flat-KM}, we define a symmetric Euclidean space in the fixed point set of the Cartan subalgebra and use the exponential map to obtain an isomorphism.\\
Let $G$ be an almost split Kac--Moody $\RR$-group. Then, based on \Cref{Sec:RR-Forms} \ref{itm:T_Wirkung}, we know that the Galois action on $G$ also acts on $\h$. Further, a rectified Galois action stabilizes the standard torus $T$ of $G$ and since $G$ satisfies the condition \eqref{Eq:DCS}, we can define 
$$
\root{\aa} \coloneqq \aa \cap \h^\Gamma.
$$
Note, that $\root{\aa}$ is exactly the fixed point set of the real form $\aa$ under the Galois action. Recall, that $G(\RR) \leq G$ and by $\root{\aa} \subseteq \aa$ we can use the exponential map given in \Cref{Sec:TitsFunctor}: restricting the map onto $\root{\aa}$ leads to 
$$
\res{\root{\aa}}{\KMExp} \colon \root{\aa} \to A.
$$
The next step is to analyze the image of $\KMExp(\root{\aa})$. Therefore, let 
$$
h = \sum_{i=1}^{2n-l} r_i v_{i} \in \h = \Lambda^\vee \otimes_{\ZZ} \CC
$$
be an element of the Cartan algebra, where $r_i \in \CC$ are scalars and $v_i \in \Lambda^\vee$ are the $\ZZ$-basis of $\Lambda^\vee$ (cf.\ \Cref{Sec:TitsFunctor}). Now we can compute
\begin{align*}
    \os(\KMExp(h)) &= \os \hspace{-2pt} \left( \prod_{i=1}^{2n-l}\left(e^{r_i}\right)^{v_i}\right) = \prod_{i=1}^{2n-l} \sigma (e^{r_i})^{v_{\os(i)}} = \prod_{i=1}^{2n-l} (e^{\overline{r_i}})^{v_{\os(i)}} \\
    \KMExp(\os(h)) &= \KMExp \hspace{-2pt} \left(\sum_{i=1}^{2n-l} \overline{r_i} v_{\os(i)} \right) = \prod_{i=1}^{2n-l} (e^{\overline{r_i}})^{v_{\os(i)}}.
\end{align*}
The image of the exponential map restricted to $\root{\aa}$ is given by $A(\RR) \coloneqq A \cap G(\RR)$, since the exponential map is Galois-equivariant.\\
Note that $\root{\aa} \cong \RR^{m}$, where $m < 2n-l$, but the actual size of $m$ is not clear. It depends on how exactly the torus of $G(\RR)$ looks and splits.

\begin{Proposition}\label{Thm:Galois-Flats}
Let $X(\RR) = G(\RR)/K(\RR)$ be the symmetric space attached to the group of rational points of an almost split Kac--Moody $\RR$-group $G$. Moreover, equip $\root{\aa}$ with the Euclidean symmetric space structure. Then for any $g \in G(\RR)$ there is an isomorphism of symmetric spaces given by 
\begin{align}
    \begin{split}
        \alsp{\p}_g \colon \root{\aa} &\to gA(\RR)K(\RR) \\
        X &\mapsto g \left(\res{\root{\aa}}{\KMExp}\hspace{-4pt}(X)\right) K(\RR).
    \end{split}
\end{align}
The subset $gA(\RR)K(\RR)$ of $X(\RR)$ is a Euclidean flat of dimension $\dim \root{\aa}$. Call for any $g \in G(\RR)$ the set $gA(\RR)K(\RR)$ a \emph{standard $\RR$-flat}.
\end{Proposition}

\begin{Proof}
Since $\alsp{\p}_g$ is obtained from the map of \Cref{Prop:Flat-KM} by restricting the exponential map to a subset of $\aa$, the assertion that $\alsp{\p}_g$ is an isomorphism of symmetric spaces follows with the same arguments as in the proof of \Cref{Prop:Flat-KM}.\\
Assume that $g = e$, then it remains to check that $A(\RR)K(\RR)$ is a closed subset, since the multiplication is continuous. To do this, we can use the Iwasawa decomposition on $G$ and the fact, that $AK$ is a closed subset in the symmetric space associated to $G$. By \Cref{Prop:SubspaceTopo-Equal}, the induced subspace topology on $G(\RR)$ is the same as the relative Kac--Peterson topology on $G(\RR)$. Hence, if $AK \cap G(\RR) = A(\RR)K(\RR)$, the desired statement follows.\\
Recall from \Cref{Thm:Iwasawa} that an element $g \in AK$ can be written as $ak$ for a unique $a \in A$ and a unique $k \in K$. Assume that $ak \in AK \cap G(\RR) \subseteq G$ and let $\os$ a rectified Galois automorphism. Then we have $\os(A) \subseteq A$ and $\os(K) \subseteq K$. From the condition \eqref{Eq:DCS} it further follows that
$$
\os(ak) = ak = \os(a) \os(k).
$$
By the uniqueness of the decomposition, $AK \cap G(\RR) = A(\RR)K(\RR)$ follows. Hence the set $A(\RR)K(\RR)$, and so $gA(\RR)K(\RR)$ for all $g \in G(\RR)$, is closed in $X(\RR)$.
\qed
\end{Proof}

As for symmetric spaces of split Kac--Moody groups, we can give the characterization of flats. This is possible, since almost split Kac--Moody $\RR$-groups posses an RGD system and the root groups has a good algebraic group structure.\\
In the proof of the characterization, one uses that the group model of the symmetric spaces is isomorphic to the coset model. But, the isomorphism between the coset model and the group model established in \Cref{Sec:GroupModel} is independent of the Kac--Moody group. It is essentially based on the isomorphism specified in \Cref{Rem:QaudraticRepresentation} for abstract symmetric spaces according to Loos and on the existence of the twist map $\tau$. Hence, as a result we have
$$
(X(\RR), \alsp{\mu}) \cong (\alsp{\tau}(G(\RR)), \alsp{\mu}_\tau).
$$

\begin{Theorem}
Let $G$ be almost split Kac--Moody $\RR$-group and denote by $G(\RR)$ the rational points of $G$. Denote by $X(\RR) = G(\RR)/K(\RR)$ the symmetric space attached to $G$. Then every weak flat is contained in a standard flat. Especially 
\begin{myitems}[itemsep=0pt]
    \item standard flats are exactly the maximal flats;
    \item all weak flats are Euclidean;
    \item all weak flats are flats;
    \item $G(\RR)$ acts transitively on maximal flats.
\end{myitems}
\end{Theorem}

\begin{Proof}
As in the general split case, the proof works verbatim as the proof of \cite[Theorem 5.17]{FHHK_2017}, except claim $5$. There it is stated that the stabilizers of opposite spherical residues are reductive Lie groups. In the case of the group of rational points of an almost split Kac--Moody $\RR$-group we get the same:\\
From \cite[Proposition 6.27]{Abramenko_2010} follows that the stabilizer of a spherical residue is equal to a spherical parabolic subgroup. From \cite[Definition 12.1.3]{Remy_2002} follows that a parabolic subgroup in $G(\RR)$ is equal to a parabolic subgroup in $G$. Using the same arguments as in the last step of the proof of \Cref{Lem:Equivalence}, it can be be interpreted as a subgroup of $\GL_m(\RR)$, $m \in \NN$.\\
As in the proof of \Cref{Thm:Characterization-Flats}, we have to embed the involution $\Theta$ into the matrix group $\GL_m(\RR)$. According to \Cref{Lem:Sigma-Theta-Tauscht} we know that $\Theta$ commutes with every rectified Galois automorphism $\widehat{\sigma}$ and therefore we can embed $\Theta$ on these parabolic $\RR$-groups as transpose inverses. Thus, the arguments of claim $5$ from the proof of \cite[Theorem 5.17]{FHHK_2017} apply verbatim.
\qed
\end{Proof}


\bibliographystyle{alpha}
\bibliography{main}

\Addresses

\end{document}